\documentclass[11pt]{article}
\usepackage{amsmath,amsthm,amssymb,mathrsfs,amscd,amsthm,amscd,amsfonts,graphicx}
\usepackage{eucal}
\usepackage[T1]{fontenc}
\usepackage[utf8]{inputenc}
\usepackage{indentfirst,url,xypic}
\usepackage{lipsum}
\usepackage{faktor}
\usepackage[all]{xy}
\usepackage{url}
\usepackage[colorlinks,plainpages]{hyperref}
\usepackage{verbatim}
\usepackage{fancyvrb}
\usepackage{listings}
\usepackage{mathtools}
\usepackage{algorithm}
\usepackage[noend]{algorithmic}

\usepackage[colorlinks,plainpages]{hyperref}
\hypersetup{
	colorlinks=true,
	linkcolor=blue,
	filecolor=magenta,      
	urlcolor=cyan,
	citecolor=blue
}
\renewcommand{\qed}{\phantom{l}\hfill\qedsymbol}

\DeclareMathOperator{\id}{id}

\DeclareMathOperator{\Act}{Act}

\DeclareMathOperator{\Nil}{Nil}

\def\NRed{\mbox{NRed}}
\def\Red{\mbox{Red}}
\def\NFlat{\mbox{NFlat}}

\def\NWor{\mbox{NWor}}
\def\NNor{\mbox{NNor}}

\def\nga{\widetilde}

\def\Coker{\mbox{Coker}}
\def\NNor{\mbox{NNor}}

\newcommand{\td}{\Longleftrightarrow}

\newcommand{\gt}{\overline}
\newcommand{\tru}{\setminus}
\newcommand{\sr}{\longrightarrow}
\newcommand{\Sr}{\Longrightarrow}
\newcommand{\mtn}{\rightarrow}

\newtheorem{Definition}{ Definition}[section]
\newtheorem{Theorem}[Definition]{Theorem}
\newtheorem{Remark}[Definition]{Remark}
\newtheorem{Proposition}[Definition]{Proposition }
\newtheorem{Corollary}[Definition]{Corollary}
\newtheorem{Example}[Definition]{Example}

\newtheorem{Proof*}{Proof}
\newtheorem{Lemma}[Definition]{Lemma}

\newcommand{\B}{{\mathbb B}}
\newcommand{\Z}{{\mathbb Z}}

\newcommand{\C}{{\mathbb C}}
\newcommand{\Q}{{\mathbb Q}}
\renewcommand{\P}{{\mathbb P}}

\newcommand{\kb}{{\mathcal B}}
\newcommand{\kc}{{\mathcal C}}
\newcommand{\kd}{{\mathcal D}}

\newcommand{\kf}{{\mathcal F}}
\newcommand{\kg}{{\mathcal G}}
\newcommand{\kh}{{\mathcal H}}

\newcommand{\km}{{\mathcal M}}
\newcommand{\kn}{{\mathcal N}}
\newcommand{\ko}{{\mathcal O}}
\newcommand{\kp}{{\mathcal P}}

\newcommand{\kx}{{\mathcal X}}

\DeclareMathAlphabet{\mathsc}{U}{rsfs}{m}{n}

\newcommand{\kt}{\mathsc{T}}

\DeclareMathOperator{\Sing}{Sing}
\DeclareMathOperator{\Spec}{Spec}

\newcommand{\es}{\text{es}}

\DeclareMathOperator{\Ann}{Ann}

\DeclareMathOperator{\Def}{\kd\!\!\:\it{ef}}

\DeclareMathOperator{\SimNor}{SimNor}

\DeclareMathOperator{\depth}{depth}

\DeclareMathOperator{\Ker}{Ker}

\DeclareMathOperator{\mt}{mt}

\DeclareMathOperator{\ord}{ord}

\DeclareMathOperator{\pr}{pr}

\newcommand{\red}{\text{\it red}}
\DeclareMathOperator{\Reg}{Reg}

\begin{document}

\title{Equisingular and Equinormalizable Deformations of Isolated Non--Normal Singularities}
\author{Gert-Martin Greuel\\[1.0ex]{\em Dedicated to Henry Laufer on the occasion of his 70th Birthday}}


\maketitle

\begin{abstract}
We present new results on equisingularity and equinormalizability of families with isolated non--normal singularities (INNS) of arbitrary dimension. We define a $\delta$--invariant and a $\mu$--invariant for an INNS and
prove necessary and sufficient numerical conditions for equinormalizability and weak equinormalizability using $\delta$ and $\mu$. For families of generically reduced curves, we investigate the topological behavior of the Milnor fibre and characterize topological triviality of such families. 
Finally we state some open problems and conjectures.
In addition we give a survey of classical results about equisingularity and equinormalizability so that the article may be useful as a reference source.
\end{abstract}

\renewcommand{\contentsname}{Table of Contents}
\tableofcontents


\addcontentsline{toc}{section}{Introduction}


\vspace{1.0cm}

\noindent \textbf{\Large Introduction}\\

The main purpose of this article is to prove necessary and sufficient numerical conditions for equisingularity and equinormalizability of families of isolated non--normal singularities. In order to put these into perspective, we start with  a survey of classical results for families of reduced curve singularities, partly from a non--classical point of view. The main new results concern deformations of isolated non--normal singularities. We define a $\delta$--invariant and a $\mu$--invariant for these singularities und use them to prove necessary and sufficient numerical conditions for equinormalizability and weak equinormalizability. We determine also the number of connected components of the Milnor fibre. For families of generically reduced curves, we investigate the topological behavior of the Milnor fibre and characterize topological triviality of such families using $\delta$ and $\mu$. 
Finally we state some open problems and conjectures.
\medskip

A germ $(Z,z)$ of an arbitrary complex space $Z$ is said to have an isolated non--normal singularity (INNS for short) if $ U \smallsetminus\{z\}$ is normal for some neighbourhood $U$ of $z$, e.g. if $z$ is an isolated singularity of $Z$. We allow in particular that $Z$ is not reduced at $z$, which is a natural assumption since this happens necessarily for fibres of families where the total space has low depth (see \cite{St15} for interesting examples of deformations of isolated non--normal singularities). 
\medskip

For $(Z,z)$ an INNS, we introduce a $\delta$--invariant $\delta(Z,z)$ and a $\mu$--invariant $\mu(Z,z)$, which take care of the embedded component of $Z$ at $z$. These invariants generalize the definitions given in \cite{BuG80} for reduced and in \cite{BrG90} for non--reduced curve singularities and we continue to call $\mu(Z,z)$ the Milnor number of $(Z,z)$ if $(Z,z)$ is a curve singularity. 
\medskip

The main general result is that a deformation $f:(X,x)\to (\C,0)$ of an INNS $(X_0, x)$ admits a simultaneous normalization iff, for a good representative $f:X\to T$, $\delta(X_t)$ is constant for $t \in T$ and if the Milnor fibre $X_t, t\neq 0$, has no isolated points. This generalizes the well--known result of Teissier for deformations of reduced curve singularities and of Br\"ucker and the author for non--reduced isolated curve singularities.
Moreover, we study the problem of simultaneous weak normalization (where the fibres admit a simultaneous weak normalization), which has not been considered so far. We prove, among others, that a family of reduced isolated non--normal singularities admits a simultaneous weak normalization iff $\mu(X_t) - \delta(X_t)$ is constant.
\medskip

We apply the general results to families of curves and study several equisingularity conditions, in particular topological triviality. It is known by \cite{BuG80} that in a family of reduced curve singularities with section $\sigma$ the constancy of $\mu(X_t)$ along $\sigma$ is equivalent to topological triviality. We show by an example that  this is no longer the case for families of non--reduced curve singularities. The reason is that, in a family where the special fibre has an embedded component, the Milnor fibre may become disconnected, even if  $\mu(X_t)$ is constant. We provide an algebraically computable formula for the number of connected components for arbitrary families of isolated non--normal singularities. The nice thing about $\mu$ is, that in a family of generically reduced curves $\mu(X_t)$ is constant iff the topological Euler characteristic of $X_t$ is constant. Moreover, we prove that a family with section $\sigma$ is topologically trivial iff $\delta(X_t,\sigma(t))$ and the number of branches $r(X_t,\sigma(t))$ are  constant and that this is equivalent to $X_t$ being connected and $\mu(X_t,\sigma(t))$ being constant.
\bigskip

We give now an overview of the different sections of the paper. Section 1 reviews the classical  characterizations of equisingular families of plane curve singularities. These studies go back to Oskar Zariski's studies in equisingularity in the 1960ties in connection with his attempts to prove resolution of singularities by induction on the dimension. We review the basic results by Oskar Zariski,  Bernard Teissier and  L\^{e} D\~{u}ng Tr$\acute{\text{a}}$ng.
\medskip

In Section 3 we study equinormalizable deformations of a reduced plane curve singularity $(C,0) \subset (\C^n,0)$ and, in connection with this, the $\delta$--constant stratum in the semiuniversal deformation $\Phi$ of $(C,0)$. The non--classical approach here is that we consider deformations of the parametrization of $(C, 0)$, instead of deformations of the equation of $(C, 0)$. Deformations of the parametrization are much easier and the main result (taken from \cite{GLS07}) is, that the semiuniversal deformation of the parametrization is obtained by pulling back $\Phi$  to the normalization of the $\delta$--constant stratum of $\Phi$ (c.f. Theorem \ref{Theo4.1}). 
\medskip

In Section 4 we reconsider equisingular deformations (es-deformations) of a reduced plane curve singularity $(C,0)$ by using again deformations of the parametrization, as developed in \cite{GLS07}. This leads to a description of the semiuniveral es-deformation of the parametrization, which is as simple and explicit computable as the semiuniveral of $(C,0)$ itself (Theorem \ref{Theo4.2}). The parametric approach provides also an easy proof of Jonathan Wahl's result that the $\mu$--constant stratum in the semiuniversal deformation of $(C,0)$ is smooth (Theorem \ref{Theo4.6}).
\medskip

The presentation of our new results starts in Section 5. We introduce the $\varepsilon$--invariant (which measures the embedded component), the $\delta$--invariant and the $\mu$-invariant for isolated non--normal singularities (Definition \ref{Def5.2}). We describe the weak normalization of an INNS and characterize the class of weakly normal INNS by the minimum of the $\delta$--invariant. Moreover, we give some formulas which are useful for the computation of the new invariants and elaborate this in an example, including the {\sc Singular}--code (Example \ref{Ex5.6}).
\medskip

Section 6 contains the general concepts of simultaneous normalizations and equinormalizability and we prove some consequences for maps that admit a simultaneous normalization.  Moreover, we review the more recent results by Chiang--Hsieh, Lipman and Koll$\acute{\text{a}}$r, who treat simultaneous normalizations over  an arbitrary normal (resp. weakly normal) base space. Hung--Jen Chiang--Hsieh and Joseph Lipman reconsidered in 2006 Teissier's results and, in addition, families of projective varieties, replacing the constancy of delta by that of the Hilbert polynomial (Theorem \ref{Theo10}). In 2009 J$\acute{\text{a}}$nos Koll$\acute{\text{a}}$r proved the existence of a fine moduli space for the functor of simultaneous normalization of a projective morphism  (Theorem \ref{Theo6.13}), thus generalizing the results of Chiang--Hsieh and Lipman. 
\medskip

The main new outcomes  are presented from Section 7 on and later. In Section 7 we characterize simultaneous normalization of 1--parametric families $f:(X,x)\to (\C,0)$ of isolated non--normal singularities, mentioned at the beginning of this introduction. We need partial normalizations of the restriction of $f$ to subspaces of $X$, which we call moderations. The introduction of moderations unifies the proof of the main results of this section. The first main result is a local and global generalization to families of INNS of the classical $\delta$--constant criterion:

\begin{Theorem}{(cf. Theorem \ref{theorem1}, Theorem \ref{theo6.16}):}
\begin{enumerate}
\item [(1)]  Let $f : (X,x) \mtn (\C,0)$ be flat with $(X_0,x)$ an INNS of dimension $\geq 1$ and $f:X \to T$ a good representative.
Then $f$ is equinormalizable if and only if $f$ is $\delta$--constant and $X_t$ has no isolated points for $t \neq 0$.

\item [(2)] Let $f: X\to T$ be a flat morphism of complex spaces with $T$ a $1$--dimensional complex manifold, such that the non--normal locus of $f$ is finite over $T$. 
Then $f$ is equinormalizable iff $\delta(X_t)$ is constant on $T$ and the $1$--dimensional part of $X$  is smooth and does not meet the higher dimensional components of $X$.
\end{enumerate}
\end{Theorem}

We like to mention that we do not need to assume that the fibres are equidimensional (as this is the case in \cite{CL06} and \cite{Ko11}), the only necessary topological condition is that the fibres have no isolated points.
\bigskip

In Section 8 we prove a connectedness result for the Milnor fibre of an arbitrary morphism $f : (X,x) \mtn (\C,0)$:

\begin{Theorem}(cf. Theorem \ref{Prop7.4}, Corollary \ref{Cor8.4}):

Let $f:(X,x)\to (\C, 0)$ be a morphism of complex germs, $f: X \to T$ a good representative and $F := X_t, t\neq 0$, the Milnor fibre of $f$.

\begin{enumerate}
\item [(1)] If $(X_0,x)$  is reduced then F is connected.
\item [(2)] If $(X,x)$  is irreducible and if there exist points $y \in X_0$ arbitrary close to $x$ such that  $(X_0,y)$ is reduced, then $F$ is irreducible.
\end{enumerate}
\end{Theorem}

This a joint result with Bobadilla and Hamm; we prove (1) under a weaker assumption. In Proposition \ref{Prop7.5} we use this to determine an algebraically computable formula for the number of connected components of the Milnor fibre of a family of INNS.
\medskip

Moreover, in this section we study $\mu$--constant families of INNS and give in particular a numerical criterion for simultaneous weak normalization:

\begin{Theorem}(cf. Theorem \ref{Theorem7.9}):

Let $f:(X,x)\to (\C, 0)$ be flat with $(X_0,x)$ a reduced INNS of dimension $\geq 1$, $f: X \to T$ a good representative and $\omega: \widehat{X}\to X$ the weak normalization of $X$.
Then 
$$ \mu(X_t) - \delta(X_t) \text{ is constant} \Leftrightarrow \omega \text{ is a weak normalization of } f $$ 
i.e. $ \widehat{f} = f \circ \omega:  \widehat{X}\to T$ is flat with weakly normal fibres and the restriction $\omega_t: \widehat{f}^{-1}(t) \to X_t$ is the weak normalization of $X_t$ for $t \in T$.
 \end{Theorem}
\bigskip

The application of the general results of Section 7 and 8  to families of generically reduced curves and  to the study of several equisingularity conditions, in particular topological triviality, is presented in Section 9. The main results are the following:

\begin{Theorem}(cf. Theorem \ref{Theorem7.10}, Corollary \ref {Cor7.12}):

Let $f:(X,x)\to (\C,0)$ be flat with $(X_0, x)$ a generically reduced curve singularity and $f:X\to T$ a good representative. Then the following holds:
\begin{enumerate}
\item [(1)]  $\mu(X_0)-\mu(X_t)=1-\chi(X_t)$ for $t\in T$.
\item [(2)] The following are equivalent:
  \begin{enumerate}
\item [(i)] 	$\mu(X_t)-b_0(X_t)$ is constant,
\item [(ii)] 	$b_1(X_t)=0$,
\item [(iii)]   each connected component of $X_t$ is contractible.
  \end{enumerate}
\end{enumerate}
\end{Theorem}

This shows that the definition of the Milnor number $\mu$ for generically reduced curves (with arbitrary non--reduced embedded points) has good topological properties. If the family has a singular section, we can say more:

\begin{Theorem}(cf. Theorem \ref{Theorem7.14}, Corollary \ref {Cor7.16}):

Let $f:(X,x)\to (\C,0)$ be flat, $(X_0, x)$ a generically reduced curve singularity and $\sigma:(\C,0) \to (X,x) $ a section of $f$ such that $X_t\smallsetminus\sigma(t)$ is smooth. 
\begin{enumerate}
\item [(1)] The following are equivalent:
\begin{enumerate}
\item [(i)]  $f:X\to T$ is topologically trivial,
\item [(ii)]  $X_0$ and $X_t$ are embedded topologically equivalent for $t\in T$,
\item [(iii)]  $(X,x)$ is pure $2$--dimensional and $f:X\to T$ admits a weak simultaneous resolution,
\item [(iv)]  $(X,x)$ is pure $2$--dimensional and $\delta(X_t, \sigma(t))$ and $r(X,\sigma(t))$ are constant,
\item [(v)]  $\mu(X_t)$ is constant and $X_t$ is connected for $t\in T$.
\end{enumerate}

\item [(2)] $f$ admits a strong simultaneous resolution $\Leftrightarrow$ (i) \ldots (v) hold and mt($X_t,\sigma(t))$ is constant.
\end{enumerate}
\end{Theorem}
\bigskip

In Example \ref{Ex7.5} we give examples to show that $\mu$ is in general not semicontinuous and that $\mu$--constant does not imply topological triviality if $(X_0,x)$ is not reduced.
\medskip

We finish the article with some comments, open problems and conjectures in Section 10.
\medskip

The appendix in Section 11 contains notations and definitions about deformations of complex spaces and maps and which are particularly useful for understanding the relation between deformations of the equation and deformations of the parametrization. Moreover, we fix the notion of a good representative of a morphism of germs, to be used throughout this article.
\medskip

For all results we give either full proves or precise references, hoping that this (rather long) article may be useful as a reference source. 
\medskip

{\bf Acknowlegment:} I would like to thank J. Bobadilla and H. Hamm for the cooperation in proving Theorem \ref{Prop7.4} and  both as well as J.J. Ballesteros, L\^e D\~ung Tr\'ang, J. Lipman and J. Snoussi for additional references and fruitful discussions.

\medskip
\smallskip

\section{Equisingularity for plane curve singularities}

In a series of three papers ''Studies in Equisingularity I, II, III'' (cf. \cite{ZaI65}, \cite{ZaII65}, \cite{ZaIII68}), published in 1965 and 1968, Oscar Zariski initiated the study of equisingular families of algebroid hypersurfaces over an algebraically closed field of characteristic $0$. The motivation behind this study was the idea to prove resolution of singularities by projection and by induction on the dimension (see also \cite{Za77} and \cite{Li00}).
\medskip

The rough and much simplified idea is as follows. Let $X$ be an $r$--dimen\-sional hypersurface and consider a (generic) projection of $X$ onto a $k$--dimensional smooth space $S$. Then the fibres $X_s, s\in S$, of this projection are $k$--codimensional hypersurfaces, which we know to resolve by induction for each individual $s\in S$. If the family $\{X_s, s\in S\}$ is ''sufficiently equisingular'' (in the sense that, for a given point $x_0\in X_{s}$, the singularities of the nearby fibres are of the ''same type'' as the singularity of $(X_{s},x_0)$), then the procedure to resolve the singularities of a single fibre $X_{s}$ should resolve the nearby fibres simultaneously and hence the singularities of $X$.
\medskip

In the above cited papers Zariski considered the case that $\dim S= \dim X-1$ and that the projection $f:X\to S$ has a section $\sigma:S\to X$ such that $W=\sigma(S)$ is the (smooth) singular locus of $X$. In this case the fibres $X_s$ are plane curves and Zariski gave a precise definition of the intuitive idea of equisingularity. He defined {\em equisingularity of $X$ along $W$}, i.e. for families of plane curve singularities $(X_s, \sigma(s))$, and proved several equivalent characterizations. One of these characterizations is that $f:X\to S$ admits a simultaneous resolution by blowing up $X$ at $W$, such that the blown up family admits again a simultaneous resolution by blowing up sections etc., until the special fibre (and then all nearby fibres) are resolved. Hence, if the hypersurface $X$ has a smooth singular locus and if an equisingular projection $f: X\to S$ exists, then the singularities of $X$ can be resolved. Unfortunately, such an equisingular projection does not always exist and resolution of singularities is in general not possible by this method.
\medskip

It was Heisuke Hironaka who proved the resolution of singularities in general for arbitrary algebraic varieties over an algebraically closed field of characteristic $0$ in 1964 by blowing up smooth centres and by an ingenious induction, which reduces certain invariants of the singularities. This is a deep result for which Hironaka received the fields medal in 1970. Nowadays there exist simpler algorithmic proofs and even implementations, e.g. in {\sc Singular} (see \cite{FP13} in \cite{DGPS15}). For a historical account of resolution in characteristic $0$ see \cite{Ha00}. The resolution of singularities over fields of positive characteristic is an important and still open problem for varieties of dimension bigger than 3, see \cite{Ha10} for an overview. 
\medskip

Independant of resolution of singularities, Zariski's Studies in Equisingularity had great influence on the development of singularity theory. We give a short account of equisingularity for families of plane curves which is of interest to us and which motivated the study of simultaneous normalization, to be dicussed in the subsequent sections. In the following we mainly work in the category of complex analytic space and germs.
Zariski defines the equivalence of two reduced plane curve singularities in three different ways (which he proves to be equivalent), by induction on the number of successive blowing ups needed to resolve the singularities.
\medskip

We recall his definition of \textbf{(a)--equivalence}. Let $(C,0)$ and $(D,0)$ be two reduced curve singularities in $(\C^2, 0)$. If $(C,0)$ is smooth resp. an ordinary node (i.e. analytically isomorphic to $x\cdot y=0$) then $(D,0)$ is equivalent to $(C,0)$ if $(D,0)$ is smooth resp. an ordinary node. Denote by $\sigma^\ast (C,0)$ the minimal number of blowing ups needed to obtain a good embedded resolution of $(C,0)$. This means that the reduced total transform has only ordinary nodes as singularities. $\sigma^\ast(C,0)=0$ means that $(C,0)$ is either smooth or has an ordinary node. Now let $\sigma^\ast(C,0)>0$ and assume that for curve singularities $(\Gamma, 0)$  with $\sigma^\ast(\Gamma,0)<\sigma^\ast(C,0)$ the notion of equivalence has already been defined.
\medskip

Denote by $C'$ the strict transform of $(C,0)$ obtained by blowing up $0\in \C^2$ and let $t$ denote the number of connected components of $C'$, corresponding to the different tangents of $(C,0)$. Let $\{0_1, \ldots, 0_t\}$ be the intersection points of $C'$ with the exceptional divisor $E$ and denote the germs of $C'$  at $0_\nu$ by $(C', 0_\nu)$. Then $\sigma^\ast(C', 0_\nu)$ is smaller than $\sigma^\ast(C,0)$. Now Zariski's definition  (in \cite{ZaI65}) is as follows:
\medskip

 An $(a)$--equivalence between $(C,0)$ and $(D,0)$ is a bijection $\pi$ between the set of branches $\{(C_i, 0)\}$ of $(C,0)$ and $\{(D_i, 0)\}$ of  $(D, 0), i=1, \dots, r,$ satisfying 
\begin{enumerate}
\item [(i)] $\pi$ is tangentially stable,  i.e. $\pi(C_i, 0)$ and $\pi(C_j, 0)$ have the same tangent iff this holds for $(C_i, 0)$ and $(C_j, 0)$ for any $i, j$,
\item [(ii)] if $(D_i, 0)=\pi(C_i, 0)$, then $\text{mt}(D_i, 0)=\text{mt}(C_i, 0)$,  $i=1, \ldots, r$,
\item [(iii)] the induced bijection $\pi'$ between the branches of $(C', 0_\nu)$ and of $(D', 0_\nu)$ is an $(a)$--equivalence.
\end{enumerate}

Zariski introduces also $(b)$--equivalence and formal equivalence and proves in \cite{ZaI65} that all three definitions are equivalent. We call any of these equivalent conditions \textbf{Zariski equivalence}.

The following theorem gives other characterizations which are more intuitive.

\begin{Theorem} {\bf(Zariski)}\\
Let $(C,0)$ and $(D,0)$ be two reduced plane curve singularities. The following are equivalent.

\begin{enumerate}
\item [(1)] $(C,0)$ and $(D,0)$ are Zariski--equivalent.
\item [(2)] There exists a bijection $\pi$ between the set of branches $\{(C_i, 0)\}$ of $(C,0)$ and $\{(D_i, 0)\}$ of $(D,0)$, $i=1, \ldots, r$, such that 
\begin{enumerate}
\item [(i)] $(C_i, 0)$ and $\pi(C_i, 0)$ have the same Puiseux pairs for $i=1, \ldots, r$,
\item [(ii)] the intersection multiplicity of $(C_i, 0)$ and $(C_j, 0)$ coincides with the intersection multiplicity of $\pi(C_i, 0)$ and $\pi(C_j, 0)$ for all $i\neq j$,
\end{enumerate}
\item [(3)] $(C,0)$ and $(D,0)$ are embedded topologically equivalent.
\end{enumerate}
\end{Theorem}


We say that two embedded complex germs  $(X,x) \subset (\C^n, x)$ and $(Y,y) \subset(\C^n, y)$ with representatives $X$ resp. $Y$ in $\C^n$ are \textbf{embedded topologically equivalent} if there exists a sufficiently small $\varepsilon_0>0$ such that for each $0<\varepsilon<\varepsilon_0$ there exists a homeomorphism $h:B_\varepsilon(x)\xrightarrow{\approx} B_\varepsilon(y)$ with $h(X)=Y$ and $h(x)=y$, where $B_\varepsilon(x)=\{z\in \C^n | \|z-x\|<\varepsilon\}$ is called a {\bf Milnor ball} for sufficiently small $\varepsilon$. Since topological maps disregard the non--reduced structure, embedded topological equivalence is a condition about the respective reductions.

The equivalence of the three characterizations is basically due to Zariski (\cite{ZaIII68} and \cite{Za71}, using the theory of saturation, see also \cite{Te75}, \S 2, Appendix).
\medskip

Recall that for a convergent power series $f\in \C\{x_1, \ldots, x_n\}$ or for a holomorphic function germ $f:(\C^n, 0)\to (\C,0)$ and analytic coordinates $x_1, \ldots, x_n$, the \textbf{Milnor number} $\mu(f)$ is defined as
\[
\mu(f)=\dim_\C \C\{x_1, \ldots, x_n\}/\langle \partial f/\partial x_1, \ldots, \partial f/\partial x_n\rangle.
\]

If the power series $f\in \C\{ x,y \}$ defines a reduced plane curve singularity $(C,0)$ = $(f^{-1}(0), 0)$ then $\mu$ and $\delta$ are related by the very useful formula of John Milnor (cf. [Mi68]),
\[
\mu(C,0)=2\delta(C,0)-r(C,0)+1,
\]
where $r(C,0)$ is the number of branches of $(C,0)$ and 
$$ \delta(C,0) = \dim_\C (n_{*} \ko_{\overline{C}})_0 / \ko_{C,0},$$
with $n: \overline C \to C$ the normalization of $C$, the {\bf delta-invariant} of $(C,0)$ (see Definiton \ref{Def5.2} for a generalization).

It is a non--trivial fact, proved by L\^{e} D\~{u}ng Tr$\acute{\text{a}}$ng in \cite{Le74}, that $\mu(f)$ depends only on the topological type of $(C,0)\subset (\C^2, 0)$ and therefore we write $\mu(C,0)$. 
\medskip

We consider now a $1$--parameter deformation $\Phi:(\kc,0)\to (\C,0)$ of a reduced plane curve singularity $(C, 0)\subset (\C^2, 0)$, together with a section $\sigma:(\C, 0)\to (\kc, 0)$. 

\begin{Theorem}\label{Theorem2.2} {\bf (Zariski, Teissier,  L\^{e}, Lejeune, Ramanujam, Timourian)}\\
Let $\Phi:(\kc, 0)\to (\C, 0)$ a deformation of the plane curve singularity $(C,0)$ with section $\sigma$ and  $\Phi :\kc\to S$ a good representative. Then the following conditions are equivalent:
\begin{enumerate}
\item [(1)] the germs $(\kc_s, \sigma(s))$ are pairwise Zariski--equivalent for $s\in S$,
\item [(2)] the Milnor number $\mu(\kc_s, \sigma (s))$ is constant for $s\in S$,
\item [(3)] the delta--invariant $\delta(\kc_s, \sigma (s))$ and the number of branches $r(\kc_s, \sigma(s))$ are constant for $s\in S$,
\item [(4)] the multiplicity $\text{mt}(\kc_s, \sigma(s))$ is constant and $\Phi:\kc\to S$ admits a simultaneous (embedded) resolution by blowing up $\sigma(S)$, which results in a family of multigerms of plane curve singularities over $S$ along finitely many sections with constant multiplicities. These sections can be inductively blown up until the special fibre $(\kc_0, \sigma(0))\cong (C, 0)$ (and then all fibres $(\kc_s, \sigma (s))\subset (B,0)$) are resolved, with normal crossings of the reduced total transform,
\item [(5)] the germs $(\kc_s, \sigma (s))$ are pairwise embedded topologically equivalent for $s\in S$,

\item [(6)] $\Phi:\kc\to S$ is induced by a topologically trivial family of functions $F_s: (\C^2,0)\to (\C,0), s\in S$. \end{enumerate}
\end{Theorem}

The proof Theorem \ref{Theorem2.2} is due to several people: (1) $\Leftrightarrow$ (5) and (4) $\Rightarrow$ (1) are due to Zariski, see Theorem 2.1, while (1) $\Rightarrow$ (4) is due to Teissier (see \cite{Te75}, \S 2). (5) $\Rightarrow$ (6) was proved by Timourian \cite{Ti77} (for a precise statement see below), while (6) $\Rightarrow$ (5) is trivial. (5) $\Rightarrow$ (2) was proved by L\^{e} D\~{u}ng Tr$\acute{\text{a}}$ng \cite{Le74} and (2) $\Rightarrow$ (5) by L\^{e} and Ramamanujam [LR76]. (2) $\Rightarrow$ (3) is due to Teissier \cite{Te76} and (3) $\Rightarrow$ (2) follows from Milnors formula relating $\mu$ and $\delta$ above. The implication (2)  $\Leftrightarrow$ (1) for irreducible curves was first proved in a short note by Monique Lejeune, L\^e D\~ung Tr\'ang and Bernard Teissier in \cite{LLT70}.
For a proof of several further characterizations and for some historical remarks see \cite{Te76}, Theorem 5.3.1. and \cite{Te77}, 3.7. The result of Timourian has been generalized by Bobadilla, who introduces in \cite{Bo13} the notion of ''cuts'' and a general technique for topological trivialization.
The notion of ''good representative'' is explained in the Appendix.

\begin{Remark}
\rm{Let $(X,0)\subset (\C^n, 0)$ be an isolated hypersurface singularity and $\Phi:\kx \to S\subset \C$ a good representative of a $1$--parameter deformation of $(X,0)$ {\em without} a given section, but with constant Milnor number $\mu(X,0)=\mu(\kx_s) := \sum\limits_{x\in \kx_s}\mu(\kx_s,x)$ for all $s\in S$. Then there exists a section $\sigma: S\to \kx$ of $\Phi$ such that $\sigma (s)$ is the only singular point of $\kx_s$ for all $s\in S$. The uniqueness of the singular point in $\kx_s$ was proved by Fulvio Lazzeri \cite{La74} and Gabrielov \cite{Ga74}, the existence of a section by Teissier \cite{Te77}.}
\end{Remark}

We explain, what topologically trivial family of functions means. For the good representative $\Phi:\kc\to S$ there exists a holomorphic function $F:B\times S\to \C$ such that $\kc$ is the hypersurface $F^{-1}(0)\subset B\times S$, and for $s\in S$ we have $\kc_s=F_s^{-1}(0)$ with $F_s=F(-, s):B\to \C$. We call $\{F_s, s\in S\}$ an \textbf{embedded topologically trivial family of functions} if there exists a neighbourhood $U$ of $(0,0)$ in $B\times S$ and a homeomorphism $h: U\xrightarrow{\approx}B\times S$, $h(0,s)=(0,s)$, such that the following diagram commutes
\[
\xymatrix{
U\ar[rr]^h_\approx\ar[dr]_G && B\times S\ar[dl]^{F_0 \times  \id}\\
 & \ \C\times S
}
\]
with $G(z,s)=(F(z, s), s)$.

The restriction of $h$ induces a homeomorphism $\kc = G^{-1}(\{0\} \times S) \xrightarrow{\approx} (F_0 \times$id$)^{-1}(0) = \kc_0 \times S$, making the $\Phi:\kc\to S $ an\textbf{ embedded topologically trivial} family of hypersurface singularities (respecting the section $\sigma(s) = (0,s)$). This implies several weaker notions of topological triviality.

In particular, that $\Phi:\kc\to S$ is (non--embedded) topologically trivial, where a family $\Phi: \kc\to S$ (resp. with section $\sigma$) is called \textbf{topologically trivial} if there exists a homeomorphism $h:\kc \to \kc_0\times S$ over $S$ (resp. s.t. $h(\sigma(s))=(0,s))$.

Moreover, setting $U_s=h^{-1}(\kb\times \{s\})$ and $h_s:U_s\xrightarrow{\approx}B$ the induced homeomorphism, then $F_s=F_0\circ h_s$. That is, $F_0$ and $F_s$ are \textbf{topologically right equivalent} for $s\in S$. It follows from \cite{LR76} that for $s,t\in S$ the germ $(\kc_s, 0)$ and $(\kc_t, 0)$ in $(\C^2,0)$ are  \textbf{embedded topological equivalent} (with the same Milnor ball in $\C^2$). \\
\smallskip

\section{$\delta$--constant and $\mu$--constant families of reduced curves}

In this section we consider families of (not necessarily plane) reduced curves. We give an overview of numerical characterizations of simultaneous normalization due to Bernard Teissier and Hung-Jen Chiang Hsieh and Joseph Lipman and of topologial triviality, due to Ragnar Buchweitz and the author of this article.

In the late 1970ies Bernard Teisssier, reconsidering the work of Zariski and Hironaka, studied again simultaneous resolutions of singularities and he initiated the study of \textbf{simultaneous normalization}. Note that for a reduced curve singularity $(C,0)$ the normalization
\[
n:(\overline{C}, \overline{0})\to (C,0),
\]
where $\overline{0}$ denotes the finite set $n^{-1}(0)$, provides a resolution of $(C,0)$ since $(\overline{C}, \overline{0})$ is smooth. However, in contrast to the previous section, where we considered embedded resolutions of plane curve singularities, the normalization gives only an abstract, non--embedded resolution, even if $(C,0)$ is plane.

Teissier introduces in [Te76] several notions of simultaneous resolution, which we recall.

\begin{Definition}\label{Def2.1}
Let $f:(X,x)\to (S,0)$ be a flat morphism of complex germs with reduced fibres and $(S,0)$ reduced.
\begin{enumerate}
\item [(1)] $f$ admits a \textbf{very weak simultaneous resolution} if, for sufficiently small representatives, there exists a morphism $\pi:\widetilde{X}\to X$ such that 
\begin{enumerate}
\item [(i)]  $\pi$ is proper and bimeromorphic\footnote{A morphism $\pi$ is bimeromorphic if there exists a nowhere dense analytic subset $A\subset X$ such that $\pi^{-1}(A)$ is nowhere dense in $\widetilde{X}$ and $\pi:\widetilde{X}\smallsetminus \pi^{-1}(A)\to X\smallsetminus A$ is an isomorphism. A subset $A$ of a complex space $X$ is called {\em nowhere dense} if its closure has no interior points. By Ritt's lemma (\cite{GR84} Chapter 5, \S3.1) an analytic set $A$ is nowhere dense if and only if $\dim(A,x)<\dim(X,x)$ for all $x\in A$.},
\item [(ii)] $\widetilde{f}=f\circ\pi:\widetilde{X}\to S$ is flat with non--singular fibres,
\item [(iii)] the induced map $\pi_s:\widetilde{X}_s\to X_s$ on the fibres is a resolution of singularities of $X_s$ for $s\in S$.
\end{enumerate}
\item [(2)] $f$ admits a \textbf{weak simultaneous resolution} along the image $(W,0)=\sigma(S,0)$ of a section $\sigma:(S,0)\to (X,x)$ of $f$ if, in addition to the properties (i)-(iii) of (1), the following holds with $\widetilde{W}:=\pi^{-1}(W)$:
\begin{enumerate}
\item [(iv)] the induced morphism $\widetilde{g}=\widetilde{f}|\widetilde{W}:\widetilde{W}\to S$ satisfies: for each $\widetilde{x}\in\widetilde{W}$ there exists an isomorphism $h:(\widetilde{W}^\red, \widetilde{x})\overset{\cong}{\rightarrow}(\widetilde{g}^{-1}(0)^\red, \widetilde{x})\times (S,0)$ such that $\widetilde{g}|(\widetilde{W}^\red,\widetilde{x})=pr\circ h$, with $pr$ the projection onto $(S,0)$.

\end{enumerate}
\item [(3)] $f$ admits a \textbf{strong simultaneous resolution} if, with the notions of (2), there exists an isomorphism $h:(\widetilde{W}, \widetilde{x})\overset{\cong}{\rightarrow}(\widetilde{g}^{-1} (0),\widetilde{x})\times (S,0)$ such that $\widetilde{g}=pr\circ h$.
\item [(4)] If (1) (i) holds and if 

\begin{enumerate}
\item [(ii)] $\widetilde{f}$ is flat with normal fibres,
\item [(iii)] $\pi_s:\widetilde{X}_s\to X_s$ is the normalization of $X_s$ for $s\in S$, 
\end{enumerate}
then $\pi: \widetilde{X}\to X$ is called a \textbf{simultaneous normalization} of $f$.
\end{enumerate}
\end{Definition}


Note that a simultaneous resolution of a family may be considered as a deformation of the resolution of the special fibre that blows down to a deformation of the special fibre. This problem is quite subtle and was studied for surface singularities by Mike Artin, Egbert Brieskorn, Jonathan Wahl and others, including Henry Laufer in \cite{La73}. For families of reduced {\em plane} curve singularities we have the following numerical characterization of strong, weak and very weak simultaneous resolution, proved by Teissier in \cite{Te75}, Theorem II.5.3.1. resp. in [Te76, Theorem I.1.3.2.] and Zariski (the equivalence between (iii) and (iv) below). 

\begin{Theorem} {\bf (Teissier, Zariski)} \\
Let $f:(X,x)\to (S,0)$ be a deformation with section $\sigma:(S,0)\to (X,x)$ of the reduced plane curve singularity $(X_0,x)$ with $(S,0)$ smooth. Then the following are equivalent (for good representatives):
\begin{enumerate}
\item [(i)] The Milnor number $\mu(X_s, \sigma(s))$ is constant for $s\in S$,
\item [(ii)] the delta--invariant $\delta(X_s, \sigma(s))$ and the number of branches $r(X_s, \sigma(s))$ are constant for $s\in S$,
\item [(iii)] $f$ admits a weak simultaneous resolution,
\item[(iv)] $f$ admits a strong simultaneous resolution.
\end{enumerate}
\end{Theorem}

For simultaneous normalization we have the following result by Teissier:

\begin{Theorem}\label{Theorem2.3} {\bf (Teissier)} \\
Let $f:(X,x)\to (\C,0)$ be a deformation of the (not necessarily plane) reduced curve singularity  $(X_0,x)$. Then the following are equivalent (for a good representative $f:X\to S$):
\begin{enumerate}
\item [(i)]  The delta--invariant $\delta(X_s) = \sum\limits_{y\in X_s}\delta(X_s,y)$ is constant for $s\in S$,
\item [(ii)] $f$ admits a very weak simultaneous resolution $\pi:\widetilde{X}\to X$, 
\item [(iii)] $f$ admits a simultaneous normalization $\pi:\widetilde{X}\to X$.
\end{enumerate}
Moreover, in this case $\widetilde{X}$ is necessarily the normalization of $X$. 

\end{Theorem}

Theorem \ref{Theorem2.3} says that a flat family $f:X\to S\subset \C$ of reduced curves (without section) admits a simultaneous normalization iff the function $s\mapsto \delta(X_s)$ is constant on $S$. Teissier considers in \cite{Te76}, I.1.3.2., also families over higher dimensional base spaces $S$ and states the same result for $S$ normal, attributing the proof to Michel Raynaud. It was observed much later by Chiang-Hsieh and Lipman in \cite{CL06} that the proof contained a gap, which they closed. 

\begin{Theorem} {\bf (Teissier, Raynaud, Chiang-Hsieh, Lipman)}\label{Theorem2.4} \\
	Let $f:(X,x)\to (S,0)$ be a flat morphism of complex germs with $(S,0)$ normal, $(X,x)$ equidimensional and with fibre $(X_0,x)$ a reduced curve singularity. Then $f$ admits a simultaneous normalization $\pi:\widetilde{X}\to X$ iff $\delta(X_s)$ is constant for $s\in S$ (for a good representative $f:X\to S$). Moreover, $\pi$ is the normalization of $X$.
\end{Theorem}

\begin{Remark}
\rm{\begin{enumerate}
\item [(1)] The proof for higher dimensional base spaces is surprisingly much more difficult then for $1$--dimensional base spaces. Chiang--Hsieh and Lipman had to use a global argument, namely the existences of the quot--scheme for coherent sheaves.
\item [(2)] The proof in \cite{CL06} is given for families of schemes over a perfect field and then the authors derive the complex analytic case from that. A simplified direct proof along the same lines in the complex analytic setting is given in \cite{GLS07}, II.2.6..
\end{enumerate}}
\end{Remark}

\begin{Example}
\rm{Let $X=\{(x,y,s)\in \C^3|x^2-y^3-sy^2=0\}$ and $f:(X,0)\to (\C,0)$ the projection on the $s$--axis. Then $X_0=\{(x,y)\in \C^2|x^2-y^3=0\}$ is the ordinary cusp and $X_s=\{(x,y)\in\C^2|x^2-y^2(y-s)=0\}$ has 1 singular point at $0$ which is an ordinary node if $s\neq 0$. Hence $f$ does not admit a strong simultaneous resolution, according to Theorem \ref{Theorem2.2}. Note that $r (X_0, 0)=1$ while $r(X_s, 0)=2$ for $s\neq 0$.

However, the family is $\delta$--constant as $\delta(X_s, 0)=1$ for all $s$. By Theorem \ref{Theorem2.3} $f$ admits a simultaneous normalization, with total space the normalization $\overline{X}$ of $X$. We can compute $\overline{X}$ as
\[
\overline{X}=\{(y,x,s,T)\in \C^4|x-Ty=y^2+y_s-T_x=y-T-s^2=x^2-y^3-sy^2=0\}
\]
with normalization map $n:\overline{X}\to X$ the projection to the first three coordinates (computed with  {\tt normal.lib} {\cite {GLP15}, from \sc Singular} \cite{DGPS15}). $\overline{X}$ is nonsingular (in fact isomorphic to $\{y-T^2=0\}$ by eliminating $x$), hence normal, and the singular point $(0,0,s)\in X_s$ has two preimage points $(0,0,s,\pm x)$ in $\overline{X}_s$ for $s\neq 0$.}
\medskip

We consider now \textbf{topological triviality} of families of reduced but not necessarily plane curves. A key role plays the Milnor number for a reduced curve singularity $(C,0)\subset(\C^n, 0)$, introduced by Buchweitz and the author in \cite{BuG79} and \cite{BuG80}. It is defined as
\[
\mu(C,0)=\dim_\C(\omega_{C,0}/d\ko_{C,0}),
\]
where $\omega_C$ is the dualizing sheaf of $C$. Moreover, it was shown in \cite{BuG80}, Proposition 1.2.1, that 
\[
\mu(C,0)=2\delta(C,0)-r(C,0)+1
\]
and that $\mu(C,0)=0$ iff $(C,0)$ is smooth.

Now we consider a situation as in Theorem 3.3, but with $\mu$ being constant. The key result is the following (Theorem 4.2.2 in \cite{BuG80}).
\end{Example}

\begin{Proposition}\label{Prop3.7}
Let $f:X\to S$ be a good representative of a deformation $f:(X,x)\to (\C,0)$ of a reduced curve singularity $(X_0, x)\subset (\C^n,0)$. The following holds for $s\in S$:
\begin{enumerate}
\item [(1)] $X_s$ is connected,
\item [(2)] $\mu(X_0, x)-\mu(X_s)=\dim_\C H^1(X_s, \C) = 1 - \chi (X_s)$,
\item [(3)] $\mu(X_0, x)-\mu(X_s)\geq \delta(X_0,x)-\delta(X_s)\geq 0$.
\end{enumerate}
\end{Proposition}

Here $H^1(X_s, \C)$ is the first singular cohomology group of the topological space $X_s$ and $\chi$ the topological Euler characteristic. An easy consequence is the following (Theorem 4.2.4 in \cite{BuG80}).

\begin{Theorem}\label{Theo3.8} {\bf (Buchweitz, Greuel)} \\
	Let $f: X\to S$ be as in Proposition \ref{Prop3.7}. Then the following are equivalent: 
\begin{enumerate}
\item [(1)] $\mu(X_s)$ is constant for $s\in S$.
\item [(2)] $\delta(X_s)$ and $r'(X_s) = \sum\limits_{y\in X_s} r'(X_s,y)$ are constant for $s\in S$.
\item [(3)] $H^1(X_s)=0$ for $s\in S$.
\item [(4)] $X_s$ is contractible for $s\in S$.
\end{enumerate}
\end{Theorem}

Here $r'(X_s, y)=r(X_s, y)-1$, which is 0 at smooth points $y$ of $X_s$. Note that if $(X_0, x)$ is not planar, then $\mu(X_s)= $ constant does in general not imply that $X_s$ has only 1 singular point or that $X_s$ is homeomorphic to $X_0$ (for examples see \cite{BuG80}).  Therefore one has to assume the existence of a section $\sigma:S\to X$ such $X_s\smallsetminus\sigma(s)$ is non--singular, which follows from $\mu(X_s, \sigma(s))$ being constant.

\begin{Theorem}\label{Theo3.9} {\bf (Buchweitz, Greuel)} \\
	Let $f:X\to S$ be as in Proposition \ref{Prop3.7} and assume that there is a section $\sigma: S\to X$ of $f$ such that $X_s\smallsetminus\sigma (s)$ is smooth for $s\in S$. The following conditions are equivalent:
\begin{enumerate}
\item [(1)] $\mu(X_s, \sigma(s))$ is constant for $s\in S$.
\item [(2)] $\delta(X_s, \sigma(s))$ and $r(X_s, \sigma(s))$ are constant for $s\in S$,
\item [(3)] $\dim_\C H^1(X_s, \C)=0$ for $s\in S$,
\item [(4)] $f$ admits a weak simultaneous resolution,
\item [(5)] $(X_0, x)$ and $(X_s, \sigma (s))$ are embedded topologically equivalent in $(\C^n, 0)$ for $s\in S$.
\item [(6)] The family $f:X\to S$ with section is topologically trivial. 
\end{enumerate}
\end{Theorem}

The proof is given in \cite{BuG80}, Theorem 5.3.1. For the definition of ''embedded topologically equivalent'' and ''topologically trivial'' see the previous section.\\

The following Corollary is well--known. It is an immediate consequence of Theorem \ref{Theo3.9}, although it was not explicitely mentioned in \cite{BuG80}.

\begin{Corollary} \label{Cor3.10}
$f$ admits a strong simultaneous resolution iff $\mu(X_s, \sigma (s))$ and $\mt(X_s, \sigma(s))$ are constant. 
\end{Corollary}

\noindent\textbf{Proof}:
\rm{By Theorem \ref{Theo3.9} we have to show that a weak simultaneous resolution is strong iff $\mt(X_s, \sigma (s))$ is constant. This is proved for generically reduced curve singularities in  Lemma \ref{Lem7.16}.
}\hfill\qed
\bigskip

By a result of Teissier \cite{Te76} strong simultaneous resolution is also equivalent to the validity of the Whitney conditions (a) and (b)  for the pair $(X\smallsetminus W, W)$. Moreover, the Whitney conditions to hold is equivalent to $\mu(X_s, \sigma(s))$ and $\mt(X_s, \sigma(s))$ being constant. This was proved by Brian\c{c}on, Galligo and Granger in \cite{BGG80}, based on \cite{BuG80}, where the authors consider further equisingularity conditions, see Section 10.1 for an overview.\\

The results of this section are generalized to flat families of generically reduced curves in Section 7 (for simultaneous normalization) and in Section 9 (for topological triviality).

\section{Deformation of the normalization\\and the $\delta$--constant stratum}

In this section we study equinormalizable deformations of a reduced plane curve singularity $(C,0)\subset(\C^2, 0)$ and, in connection with this, the $\delta$--constant stratum of a deformation. This is done by considering the deformation of the normalization resp. of  the parametrization of $(C,0)$, which is easier than the corresponding deformation of the equation of $(C,0)$. 

The main result is a theorem, published first in \cite{GLS07}, saying that the restriction of the semiuniversal deformation of $(C,0)$ to the $\delta$--constant stratum $\Delta^\delta$ is, after pull back to the normalization of $\Delta^\delta$, isomorphic to the semiuniversal deformation of the normalization map $n:(\overline{C}, \overline{0})\to (C,0)$. This gives not only a nice desription of $\delta$--constant deformations, it shows also that $\Delta^\delta$ is irreducible with smooth normalization.\\

Let $\Phi:(\kc, 0)\to (S,0)$ be a deformation of the plane curve singularity $(C, 0)$ over an arbitrary complex germ $(S,0)$, and $\Phi:\kc\to S$ a good representative. We denote by $\Sing(\Phi)$ the points in $\kc$ where $\Phi$ is not regular and call
\[
\Delta_\Phi:=\Phi(\Sing(\Phi))\subset S\ ,
\]
the \textbf{discriminant} of $\Phi$ and for $k\geq 0$
\[
\begin{array}{lcl}
\Delta^\delta_\Phi(k): & = & \{s\in S|\delta(\kc_s)\geq k\},\\
\Delta^\mu_\Phi(k): & = & \{s\in S|\mu (\kc_s)\geq k\}\ .
\end{array}
\]
Since the fibres $\kc_s$ of $\Phi$ have only finitely many isolated singularities, $\delta(\kc_s)$ and $\mu(\kc_s)$ are finite numbers and the restriction of $\Phi$ to $\Sing(\Phi)$ is a finite morphism. The discriminant $\Delta_\Phi$ of $\Phi$ is therefore a closed analytic subset of $S$ which we endow with the Fitting structure (cf. \cite{GLS07}, Def.I.1.4.5.). By \cite{GLS07}, Prop.II.2.57,  the sets $\Delta_\Phi^\delta(k)$ resp. $\Delta^\mu_\Phi(k)$ are closed analytic subsets of $S$, which we equip here with the reduced structure. The smallest of these analytic sets are those with $k=\delta:=\delta(C,0)$ resp. $k=\mu:=\mu(C,0)$, called the $\mathbf{\delta}$\textbf{--constant stratum} resp. the $\mathbf{\mu}$\textbf{--constant stratum} of $\Phi$ and are denoted by $\Delta^\delta_\Phi$ resp. $\Delta^\mu_\Phi$. The strata $\Delta^\delta_\Phi(k)$ are also called the Severi strata of $\Phi$.

Although we do not need it here, we like to mention that $\Delta^\mu_\Phi$ carries a natural, not necessarily reduced scheme structure since it can be interpreted as the stratum in $S$ parametrizing all equisingular deformations of $(C,0)$ that can be induced from $\Phi$ (cf. \cite{GLS07}, Section II.2.1.). Moreover, it was recently shown in \cite{CMS16}, that all Severi strata are defined by equations coming from a natural symplectic form on $S$.\\

We come now to the description of $\Delta^\delta_\Phi$ by using deformations of the normalization map $n:(\overline{C}, \overline{0})\to (C,0)$ (cf. Appendix).  We have the forgetful functor from deformations of $n:(\overline{C}, \overline{0})\to (C,0)$ to deformations of $(C,0)$ , or in other words, from \textbf{deformations of the normalization} of $(C,0)$ to \textbf{deformations of the equation} of $(C,0)$.

The following theorem explains the relation between deformations of the normalization of $(C,0)$ and $\delta$--constant deformations of the equation of $(C,0)\subset(\C^2, 0)$. For this we recall the seminuniversal deformation of (the equation of) $(C,0)$. Let $f\in \ko_{\C^2,0}$ define $(C,0)$ and $g_1:=1, g_2, \ldots, g_\tau\in \ko_{\C^2,0}$ be a basis of the Tjurina algebra $ T_f = \C\{x,y\}/ \langle f, \partial f/\partial x, \partial f/\partial y\rangle $ with $\tau=\tau(f) = \dim_\C T_f $ the Tjurina number of $f$.
Define with $s=(s_1, \ldots, s_\tau)$,
\[
(\kd,0):=\{(z,s)\in (\C^2\times \C^\tau, 0)\ |\ f(z)-\sum\limits^\tau_{i=1} g_i(z)s_i=0\},
\]
and let
\[
\Phi:(\kd, 0)\to (B_C, 0):=(\C^\tau,0)
\]
be the projection on the second factor. Then $\Phi$ is the seminuniversal deformation of $(C,0)$ with base space $(B_C, 0)$ and we call (for a good representative) the analytic subsets
\[
\Delta^\delta:=\Delta^\delta_\Phi \text{ resp. } \Delta^\mu:=\\\Delta^\mu_\Phi
\]
of $B_C$ the $\mathbf{\delta}$\textbf{--constant stratum of} $(C,0$) resp. the $\mathbf{\mu}$\textbf{--constant stratum of} $(C,0)$. It consists of those points $s\in B_C$ such that the fibre $\kd_s$ satisfies $\delta(\kd_s)=\delta(C,0)$ resp. $\mu(\kd_s)=\mu(C,0)$.\\

The following theorem from \cite{GLS07}, Theorem II.2.59, is the main result of this section. A short summary of the deformation theoretic background is given in the Appendix.

\begin{Theorem} \label{Theo4.1}  {\bf (Greuel, Lossen, Shustin)} \\
Let $n:(\overline{C},0)\to (C,0)$ be the normalization of the reduced plane curve singularity $(C,0)$.
\begin{enumerate}
\item [(1)] If $(\overline{\kc}, \overline{0})\to (\kc,0)\to (B_{\overline{C}\to C}, 0)$ denotes the semiuniversal deformation of $n$, then its base space $(B_{\overline{C}\to C}, 0)$ is smooth of dimension $\tau -\delta$.
\item [(2)] Let $(\kd, 0)\to (B_C, 0)$  be the semiuniversal deformation of $(C,0)$. Then the forgetful transformation $\underline{\Def}_{\overline{C}\to C}\to \underline{\Def}_C$ and the versality of $(\kd, 0)\to (B_C, 0)$ imply the existence of a morphism
\[
\alpha:(B_{\overline{C}\to C}, 0)\to (B_C, 0)
\]
satisfying:
\begin{enumerate}
\item [(i)] the image of $\alpha$ is the $\delta$--constant stratum $(\Delta^\delta, 0)$  of $(C,0)$ and $\alpha: (B_{\overline{C}\to C}, 0)$ $\to (\Delta^\delta,0)$ is the normalization of $(\Delta^\delta, 0)$.
\item [(ii)] The pull back of the semiuniversal deformation $(\kd, 0)\to (B_C, 0)$ of $(C,0)$ via $\alpha$ is isomorphic to $(\kc,0)\to (B_{\overline{C}\to C}, 0)$ and hence lifts to the semiuniversal deformation $(\overline{\kc}, \overline{0}) \to (\kc,0)\to (\B_{\overline{C}\to C}, 0)$ of the normalization $n$ of $(C,0)$.
\end{enumerate}
\item [(3)] 
\begin{enumerate}
\item [(i)] The $\delta$--constant stratum $(\Delta^\delta, 0)\subset (B_C, 0)$ is irreducible of dimension $\tau-\delta$ and has a smooth normalization.
\item [(ii)] $s\in \Delta^\delta$ is a smooth point of $\Delta^\delta$ iff each singularity of the fibre $\kd_s$ has only smooth branches.
\item [(iii)] There exists an open dense set $U\subset\Delta^\delta$ such that each fibre $\kd_s, s\in U$, has only ordinary nodes as singularities.
\end{enumerate}
\end{enumerate}
\end{Theorem}

\begin{Remark}
\begin{enumerate}
\rm{\item [(1)] The existence of  $\alpha$ follows from deformation theoretic properties. In general it is unique only up to first order. However, in our situation it is unique since it is the normalization of $(\Delta^\delta,0)$.
\item [(2)] Statement (3) (i) and (ii) was proved before with global methods by Stephan Diaz and Joseph Harris in \cite{DH88} in 1988 using an earlier result by Arbarello and Cornalba. Theorem \ref{Theo4.1} above was first published in \cite{GLS07}, although the result was already known to the author of this article since 1988 (he scetched a proof since then on several conferences).}
\end{enumerate}
\end{Remark}

\begin{Remark}\label{Rem4.3}
\rm{For the full proof of Theorem \ref{Theo4.1} we refer to \cite{GLS07}, but we like to make a few comments on the proof.
\begin{enumerate}
\item [(1)] We use at several places that $(C,0)$ is a plane curve singularity. We need this assumption to prove that $B_{\overline{C}\to C}$ is smooth. The proof for this is indirect, by showing first that the semiuniversal deformation of the normalization $n:(\overline{C}, \overline{0})\to (C,0)$ is isomorphic to the semiuniversal deformation of the parametrization $\varphi=j\circ n:(\overline{C}, \overline{0})\rightarrow(\C^2,0)$, with $j:(C,0)\hookrightarrow (\C^2,0)$ an embedding of $(C,0)$. In fact, any deformation of the the parametrization $\varphi$ can be lifted to a deformation of the pair $(n,j)$, $(\overline{C}, \overline{0})\xrightarrow{n} (C,0)\overset{j}{\hookrightarrow} (\C^2, 0)$, see \cite{GLS07}, Prop. II.2.23. This is in general not true for non--plane curve singularities. Since the parametrization $\varphi$ is a morphism between smooth (multi--)germs, its semiuniversal deformation is given by varying the coefficients of $\varphi$ without any flatness condition, and hence its base space is smooth (cf. \cite{GLS07}, Theorem II.2.38.). This is described in the next section.

\item [(2)] Another place where we use that $(C,0)$ is planar, is the proof that $\alpha:(B_{\overline{C}\to C},0)\to (\Delta^\delta, 0)$ is the normalization. It is well known that for a generic point $s\in \Delta^\delta$ the fibre $\kd_s$ has only $\delta$ ordinary nodes as singularities (statement (3) (iii)). A local computation for a node $(D,0)$ shows that its $\delta$--constant stratum is smooth (a reduced point) and isomorphic to $(B_{\overline{D}\to D}, 0)$. By openness of versality for deformations of $(C,0)$ and of $(\overline{C}, \overline{0})\to (C,0)$ it follows that $(\Delta^\delta, s)$ is smooth of codimension $\delta$.
\end{enumerate}}
\end{Remark}
\bigskip

\section{Deformation of the parametrization\\and the $\mu$--constant stratum}

We come back to equisingular deformations 
(or {\bf es--deformation} for short) 
of a reduced plane curve singularity $(C,0)\subset(\C^2,0)$. Here we describe it not by deformation of the equation as in Section 1 but by deformation of the parametrization. This approach is very natural and much easier, since deformations of the parametrization are almost effortless to describe. Both concepts can be defined for deformations over non--reduced base spaces, giving  a deformation functor in the sense of Michel Schlessinger \cite{Sch68} (see also \cite{GLS07}, App. C) and a canonical (not necessarily reduced) scheme structure on the base space of its semiuniversal deformation. When we compare es--deformations of the parametrization with es-deformations of the equation, the main problems arise in connection with deformations over non--reduced base spaces.
\bigskip

Without going through the formalities for which we refer to \cite{GLS07}, Section 2.3, let us describe deformations of the parametrization in concrete terms. Consider the commutative diagram of complex (multi--)germs
\[
\xymatrix{
(\overline{C}, \overline{0})\ar[d]_n\ar[dr]^\varphi &\\
(C,0)\ar@{^{(}->}[r]^j & (\C^2. 0)
}
\]
with $j$ the given embedding, $n$ the normalization and $\varphi$ the parametrization of $(C,0)$. If $(C,0)=(C_1, 0)\cup\cdots\cup(C_r, 0)$ is the decomposition of $(C,0)$ into irreducible components, then $(\overline{C}, \overline{0})=(\overline{C}_1, \overline{0}_1)\coprod\cdots\coprod(\overline{C}_r, \overline{0}_r)$ is a multigerm with $(\overline{C}_i, \overline{0}_i)\cong (\C,0)$ mapped via $n$ onto $(C_i, 0)$, inducing the normalization of the branch $(C_i, 0)$. We fix local coordinates $x,y$ of $(\C^2,0)$ and $t_i$ of $(\overline{C}_i, \overline{0}_i), i=1, \ldots, r$, identifying it with $(\C,0)$. Then the parametrization $\varphi=\{\varphi_i\}_{i=1, \ldots, r}$ is given by $r$ holomorphic map germs
\[
\varphi_i:(\C,0)\to (\C^2, 0)\ ,\ t_i\mapsto (x_i(t_i), y_i(t_i)),
\]
parametrizing the branch $(C_i, 0)$. Since any deformation of the smooth germs $(\overline{C}_i, \overline{0}_i)$ and $(\C^2, 0)$ is trivial, a \textbf{deformation of the parametrization  $\mathbf{\varphi}$} over a complex germ $(S,0)$ {\bf with sections} $\sigma, \overline{\sigma}$ is given by a Cartesian diagram

\[
\xymatrix{
\coprod\limits^r_{i=1}(\overline{C}_i, \overline{0}_i) = \hspace{-4ex}& (\overline{C}, \overline{0}\ar@{^{(}->}[r]\ar[d]_\varphi & (\overline{\kc}, \overline{0})\ar[d]^\Phi\ar[r]^-\cong & \coprod\limits^r_{i=1}(\C\times S, 0)\ar[d]\\
(M,0)=\hspace{-7.4ex}& (\C^2,0)\ar[d]\ar@{^{(}->}[r] & (\km, 0)\ar[d]\ar[r]^\cong & (\C^2\times S, 0)\ar[d]^{pr}\\
&\{0\}\ar@{^{(}->}[r] & (S,0)\ar@/^1pc/[u]^\sigma\ar@/_2pc/[uu]_/.9em/{\overline{\sigma}}\ar[r]^{=} & (S,0)
}
\]

\noindent with $pr$ the projection, compatible sections $\sigma$ and $\overline{\sigma}$ and a morphism $\Phi=\{\Phi_i\}_{i=1, \cdots, r}:(\overline{\kc}, \overline{0})=\coprod\limits^r_{i=1} (\overline{\kc}_i, \overline{0}_i)\to(\km, 0)\cong(\C^2\times S,0)$, given by holomorphic map germs
\[
\Phi_i:(\overline{\kc}_i, \overline{0}_i)\cong(\C\times S,0)\to(\C^2\times S,0),\]
\[
(t_i, s)\mapsto(X_i(t_i,s), Y_i(t_i,s),s),
\]
with $X_i(t_i,0)=x_i(t_i), Y_i(t_i, 0)=y_i(t_i)$. We may assume (by \cite{GLS07}, Prop. 2.2) that $\sigma$ and $\overline{\sigma}=\{\overline{\sigma}_i\}_ {i=1, \ldots, r}$ are the trivial sections, i.e. $\sigma(s)=(0,s)$ and $\overline{\sigma}_i(s)=(0,s)$ (where $0$ always denotes the origin of the corresponding germ).
\bigskip

The deformation $\Phi$ is thus given by two holomorphic map germs
\[
\begin{array}{lclcl}
X_i(t_i,s) & = & x_i(t_i)+a_i(t_i, s),\ a_i(t_i, 0)&=&0,\\
Y_i (t_i, s) & = & y_i(t_i)+b_i(t_i, s),\ b_i(t_i, 0)&=&0\ .
\end{array}
\]
Setting $\ord\varphi_i:=\min\{\ord(x_i(t_i)), \ord (y_i (t_i))\}$ we call the deformation $\Phi$ \textbf{equimultiple} (along the trivial sections $\sigma$ and $\overline{\sigma}_i)$ if for $i=1, \ldots, r$,
\[
\ord\varphi_i=\min\{\ord_{t_i}X_i(t_i, s), \ord_{t_i}Y_i(t_i, s)\}\ .
\]

Here $\ord_t X(t,s)$ denotes the \textbf{$\mathbf t$--order} of $X$, i.e. the lowest degree of non--zero terms of $X(t,s)$ considered as power series in $t$ with coefficients in $\C\{s\}$.

The definition makes sense also for non--reduced base spaces. For a reduced base space $S$, equimultiple means that the multiplicity of the branch of the plane curve singularity parametrized by $t_i\mapsto (X_i(t_i, s), Y_i(t_i, s))$ is independent of $s\in S$. 
\bigskip

We describe now equisingular deformations of the parametrization $\varphi$ of $(C,0)$.

We denote by $M\subset \C^2$ a small neighbourhood of $0$ and consider an infinitely near point $p\in \widetilde{M}$ (including $0\in M$) of $(C,0)$. That is, $\widetilde{M}$ is obtained from $M$ by a finite squence of blowing up points with  $\widetilde{C}$ the strict transform of $C$ and $p$ an intersection point of $\widetilde{C}$ with the exceptional divisor of the blowing up. Then the germ $(\widetilde{C}, p)$ is the strict transform of a collection of branches of $(C,0)$, denoted by $(C_p, 0)$. Let $(\overline{C}, \overline{p})=\varphi^{-1}(C_p, 0)\subset (\overline{C},\overline{0})$ and $\varphi_p: (\overline{C}, \overline{p})\to (\C^2, 0)$ be the restriction of $\varphi$. Since $(C_p, 0)$ and $(\widetilde{C}, p)$ have the same normalization $(\overline{C}, \overline{p}), \varphi_p$ factors through $(\widetilde{C}, p)$ and the induced map
\[
\widetilde{\varphi}_p: (\overline{C}, \overline{p})\to (\widetilde{C}, p) \hookrightarrow (\widetilde{M}, p)\cong (\C^2, 0)
\]
is a parametrization of $(\widetilde{C}, p)$.

\begin{Definition}
A deformation with section $(\Phi, \sigma, \overline{\sigma})$ of $\varphi$ over a complex germ $(S,0)$ is called \textbf{equisingular} or an \textbf{es--deformation}  of $\varphi$ if the following holds:
\begin{enumerate}
\item [(i)] $\Phi$ is equimultiple along $\sigma$.
\item [(ii)] For each infinitely near point $p$ of $(C,0)$, $\Phi$ induces an equimultiple deformation $\Phi_p:(\overline{C}, \overline{p})\times (S,0)\to (\widetilde{M}, p)\times (S,0)$ of the parametrization $\widetilde{\varphi}_p$ along some section $\sigma_p$, compatible with $\overline{\sigma}_p=\{\overline{\sigma}_i|\overline{\sigma}_i (0)\in\overline{p}\}$.
\item [(iii)] The collection $\{(\Phi_p, \sigma_p, \overline{\sigma}_p)|\ p \text{ infinitely near to } (C,0)\}$ is compatible with blowing ups relating two infinitely near points.
\item [(iv)] If $(\widetilde{M}', p')$ is obtained by a single blowing up of $(\widetilde{M}, p)$ then the source of $\Phi_{p'}$ is the blow up of the source of $\Phi_p$ along the section $\sigma_p$.
\end{enumerate}
\end{Definition}

\noindent We denote the category of es--deformations of $\varphi$ over $(S,0)$ by 
\[
\Def^{es}_{(\overline{C}, \overline{0})\to(\C^2,0)}(S,0), \text{ or by } \Def^{es}_\varphi (S,0),
\]
where the morphisms are morphisms of deformations with section, and by $\underline{\Def}^{es}_\varphi (S,0)$ the functor of isomorphism classes. For details see \cite{GLS07}, Def. 2.36, and the subsequent discussion. Since equimultiple deformations are defined for arbitrary $(S,0)$, this is also the case for es--deformations. In order to describe the semiuniversal es--deformation of $\varphi$ we have to introduce the equisingularity module of $\varphi$. We have
\[
\ko_{(\overline{C}, \overline{0})}=\overset{r}{\underset{i=1}{\oplus}} \ko_{(\overline{C}_i, \overline{0}_i)}\cong \overset{r}{\underset{i=1}{\oplus}}\C\{t_i\}
\]
with $\overline{\mathfrak{m}}\cong \overset{r}{\underset{i=1}{\oplus}} t_i\C\{t_i\}$ the Jacobson radical of $\ko_{(\overline{C},\overline{0})}$. Writing elementes of $\ko_{(\overline{C}, \overline{0})}$ as column vector we define the \textbf{equisingularity module}   {\bf $I^{es}_\varphi$ of $\varphi$} as the set of elements
\[
\left(\begin{smallmatrix}a_1\\\vdots\\a_r\end{smallmatrix}\right)\frac{\partial}{\partial x} + \left(\begin{smallmatrix}b_1\\\vdots\\b_r\end{smallmatrix}\right)\frac{\partial}{\partial y}\in \overline{\mathfrak{m}}\frac{\partial}{\partial x}\oplus\overline{\mathfrak{m}} \frac{\partial}{\partial y}, 
\]
\[\text{ such that } \{x_i(t_i)+\varepsilon a_i(t_i), y_i(t_i)+\varepsilon b_i(t_i)\} \text{ is an es--deformation of  } \varphi \text{ over } T_\varepsilon ,
\]
(for $T_\varepsilon$ see Appendix (7)).
It can be shown that $I^{es}_\varphi$ is an $\ko_{\C^2, 0}$--submodule of $\overline{\mathfrak{m}}\frac{\partial}{\partial x}\oplus \overline{\mathfrak{m}}\frac{\partial}{\partial y}$ (\cite{GLS07}, Prop. II.2.40.).
\bigskip

We identify now the submodule of $I^{es}_\varphi$ defining trivial deformations of $\varphi$. The algebra homomorphismus $\varphi^\# : \ko_{\C^2, 0}=\C\{x,y\}\to \ko_{\overline{C}, \overline{0}}=\overset{r}{\underset{i=1}{\oplus}}\C\{t_i\}$ induced by $\varphi$ maps $x$ to $\left(\begin{smallmatrix}x_1(t_1)\\\vdots\\x_r(t_r)\end{smallmatrix}\right)$ and $y$ to $\left(\begin{smallmatrix}y_1(t_1)\\\vdots\\y_r(t_r)\end{smallmatrix}\right)$ mapping the maximal ideal $\frak{m}=\langle x,y\rangle$ of $\ko_{\C^2,0}$ to the $\ko_{\C^2,0}$--submodule $\varphi^\#(\frak{m})$ of $\ko_{\overline{C}, \overline{0}}$. Then
\[
\varphi^\#(\frak{m})\frac{\partial}{\partial x}+\varphi^\#(\frak{m})\frac{\partial}{\partial y}
\]
is a $\ko_{\C^2,0}$--submodule of $\overline{\frak{m}}\frac{\partial}{\partial x}\oplus \overline{\frak{m}}\frac{\partial}{\partial y}$. Setting $\dot{x_i}(t_i)=\frac{\partial x_i}{\partial t_i}(t_i)$, $\dot{y_i}(t_i)=\frac{\partial y_i}{\partial t_i}(t_i)$ and 
\[
\dot{\varphi}:= \left(\begin{smallmatrix}\dot{x_1}(t_1)\\\vdots\\\dot{x}_r(t_r)\end{smallmatrix}\right)\frac{\partial}{\partial x}+\left(\begin{smallmatrix}\dot{y_1}(t_1)\\\vdots\\\dot{y}_r(t_r)\end{smallmatrix}\right)\frac{\partial}{\partial y}\ ,
\]

we get the following $\ko_{\overline{C}, \overline{0}}$--submodule of $\overline{\mathfrak{m}}\frac{\partial}{\partial x}\oplus\overline{\mathfrak{m}}\frac{\partial}{\partial y}$,
\[
\dot{\varphi}\cdot \overline{\mathfrak{m}}=\left\{\left(\begin{smallmatrix}c_1\dot{x_1}\\\vdots\\c_r\dot{x}_r\end{smallmatrix}\right)\frac{\partial}{\partial x}+\left(\begin{smallmatrix}c_1\dot{y_1}\\\vdots\\c_r\dot{y}_r\end{smallmatrix}\right)\frac{\partial}{\partial y}\ \big|\ c_i\in t_i\C\{t_i\}\right\}\ .
\]

Then $\{x_i(t_i)+\varepsilon  a_i(t_i), y_i(t_i)+\varepsilon b_i(t_i\}_{i=1, \ldots, r}$ is a trivial deformation (with trivial section) of $\varphi$ over $T_\varepsilon$ iff
\[
\left(\begin{smallmatrix}a_1\\\vdots\\a_r\end{smallmatrix}\right)\frac{\partial}{\partial x} + \left(\begin{smallmatrix}b_1\\\vdots\\b_r\end{smallmatrix}\right)\frac{\partial}{\partial y}\in \dot{\varphi}\cdot \overline{\mathfrak{m}}+\varphi^\#(\mathfrak{m})\frac{\partial}{\partial x}\oplus \varphi^\# (\mathfrak{m})\frac{\partial}{\partial y} \ .
\]

(cf. [GLS07], proof of Prop. II.27). We define the $\ko_{\C^2. 0}$--module
\[
T^{1, es}_\varphi:= I^{es}_\varphi / (\dot{\varphi}\overline{\mathfrak{m}}+\varphi^\#  (\mathfrak{m})\frac{\partial}{\partial x}\oplus \varphi^\# / (\mathfrak{m})\frac{\partial}{\partial y})
\]
and it follows that $T^{1,es}_\varphi$ is isomorphic as a $\C$--vector space to $\underline{\Def}^{es}_\varphi (T_\varepsilon)$, the space of isomorphism classes of {\bf infinitesimal es--deformations of $\mathbf{\varphi}$}  (\cite{GLS07}, Lemma II.2.7.).
\bigskip

We are now in the position to formulate the main theorem about es--deformations of the parametrization (\cite{GLS07}, Theorem II.238).

\begin{Theorem} \label{Theo4.2} {\bf (Greuel, Lossen, Shustin)} \\
Let $\varphi:(\overline{C}, \overline{0})=\coprod\limits^r_{i=1}(\C,0)\to (\C^2,0), t_i\mapsto(x_i(t_i), y_i(t_i))$ be a parametrization of the reduced plane curve singularity $(C,0)$ with $r$ branches. If
\[
\left(\begin{smallmatrix}a_1^j\\\vdots\\a_r^j\end{smallmatrix}\right)\frac{\partial}{\partial x}+ \left(\begin{smallmatrix}b_1^j\\\vdots\\b_r^j\end{smallmatrix}\right)\frac{\partial}{\partial y}\in I^{es}_\varphi\ ,\ j=1, \ldots, k\ ,
\]
represent a basis (resp. a system of generators) of $T^{1, es}_\varphi$, then $\Phi: (\overline{\kc},\overline{0})=\coprod\limits^r_{i=1}(\C\times \C^k, 0)\to (\C^2\times \C^k, 0)\ ,\ (t_i, s)\mapsto (X_i(t_i, s), Y_i(t_i, s), s),$
with 
\[
\begin{array}{lcl}
X_i(t_i, s) & = & x_i(t_i)+\sum\limits^k_{j=1} a_i^j(t_i) s_j,\\
Y_i(t_i, s) & = & y_i(t_i)+\sum\limits^k_{j=1}b^j_i(t_i)s_j,
\end{array}
\]
$s=(s_1, \ldots, s_k)\in (\C^k, 0)$, is a semiuniversal (resp. versal) es--deformation of $\varphi$ with trivial sections over $(\C^k, 0)$.
\end{Theorem}

\begin{Corollary}
The semiuniversal es--deformation of a parametrization of a reduced plane curve singularity has a smooth base space of dimension $\dim_\C T^{1, es}_\varphi$.
\end{Corollary}

\begin{Remark}\label{Remark4.4}
\rm{\begin{enumerate}
\item [(1)] The base space of the semiuniversal es--deformation of $\varphi$ is not only smooth, but a linear subspace of the vector space
\[
T^{1, sec}_\varphi:=\overline{\mathfrak{m}}\frac{\partial}{\partial x}\oplus\overline{\mathfrak{m}}\frac{\partial}{\partial y}/(\dot{\varphi}\overline{\mathfrak{m}}+\varphi^\# (\mathfrak{m})\frac{\partial}{\partial x}\oplus\varphi^\#(\mathfrak{m})\frac{\partial}{\partial y}),
\]
which can be identified with the base space of the semiuniversal deformation with sections of $\varphi$ (cf. \cite{GLS07}, Prop II.2.27).
\item [(2)] In the previous section we considered deformations of the normalization $n:(\overline {C}, \overline{0})$ $\rightarrow (C,0)$ {\em without section} and used them to describe the $\delta$--constant stratum of $(C,0)$ (Theorem \ref{Theo4.1}). Let  $\Def_n$ (resp. $\Def_n^{\sec}$) be the deformations (resp. deformations with section) of  $n$ and $\Def_\varphi$ (resp. $\Def^{sec}_\varphi$) the deformations of $\varphi$. We prove in \cite{GLS07}, Prop. II.2.23, that we can revover $(C,0)$ from the parametrization $\varphi:(\overline{C}, \overline{0})\to (\C^2,0)$ and that this holds for any deformation of $\varphi$. In this way we get a surjection $\Def_n\to \Def_\varphi$, inducing an isomorphism $\underline{\Def}_n\cong\underline{\Def}_\varphi$ of isomorphism classes between deformations of $n$ and that of $\varphi$ (without section), see also Remark \ref{Rem4.3}. The same holds for deformations with section. In particular, for infinitesimal deformations of $n$ and of $\varphi$ we have

\[
\begin{array}{l}\underline{Def}_n(T_\varepsilon)\cong\underline{Def}_\varphi(T_\varepsilon)\cong T^1_\varphi: =\\[2.0ex]
\ko_{\overline{C}, \overline{0}}\frac{\partial}{\partial x}\oplus\ko_{\overline{C}, \overline{0}}\frac{\partial}{\partial y}/(\dot{\varphi}\ko_{\overline{C}, \overline{0}}+\varphi^\#(\ko_{\overline{C}, \overline{0}})\frac{\partial}{\partial x}\oplus\varphi^\#(\ko_{\overline{C}, \overline{0}})\frac{\partial}{\partial Y}),
\end{array}
\]

which are vector spaces of dimension $\tau(C,0)-\delta(C,0)$. 

For deformation with sections we have
\[
\underline{\Def}^{sec}_n(T_\varepsilon)\cong T_\varphi^{1, sec}\ ,
\]
which has dimension $\tau(C,0)-\delta(C, 0)-r(C,0)+2$ by \cite{GLS07}, Prop. II.2.34 (if $(C,0)$ is not smooth).

\item [(3)] The theory of es--deformations of the parametrization was first developed in \cite{CGL07a}, even for algebroid plane curve singularities defined over an algebraically closed field of arbitrary characteristic. It is essentially the same as over the complex number. An algorithm to compute $T^{1,es}_\varphi$ and hence the seminuniversal es--deformations of $\varphi$ is given in \cite{CGL07b}, Algorithm 1 and Remark 3.9.
\end{enumerate}}
\end{Remark}

The nice thing with es--deformations of the parametrization is that they are a ''linear subfunctor'' of arbitrary deformations of the parametrization and hence easy to describe. We used them  in \cite{GLS07} to derive properties of es--deformation of the equation, in particular to show that the $\mu$--constant stratum of a plane curve singularity is smooth, as we are going to review now.
\bigskip

In Section 1 of this article  we considered 1--parameter (i.e. over $(\C, 0)$) es--deformations of the equation of a reduced plane curve singularity $(C,0)\subset (\C^2, 0)$ given by $f\in \ko_{\C^2, 0}=\C\{x,y\}$. Now we define es--deformations of $(C,0)$ over a (not necessarily reduced) complex germ $(S,0)$ of arbitrary dimension. 
\medskip

Consider a deformation $\Phi:(\kc, 0)\to (S,0)$ of $(C,0)$  with a section $\sigma:(S,0)\to (\kc, 0)$. It is isomorphic to an embedded deformation of $(C,0)$ with trivial section by Section 2.(8). That is, if $(S,0)$ is a subgerm of some $(\C^k, 0)$, then there is an isomorphism $(\kc, 0)\cong (F^{-1}(0), 0)\subset (\C^2\times S, (0,0))$ for some holomorphic map $F:(\C^2\times S,0)\to (\C,0)$ with 
\[
F(x,y,s)=f(x,y)+\sum\limits^k_{i=1} s_ig_i(x,y,s),\ \ s=(s_1, \ldots, s_k)
\]
such that $\Phi$ is the composition $(\kc,0)\xrightarrow{\cong}(F^{-1}(0), 0)\xrightarrow{pr} (S,0)$ with $pr$ the second projection and $\sigma(s)=(0,0,s)$ under this isomorphism.

\begin{Definition}
The embedded deformation $\Phi: (\kc, 0)\to (S,0)$ of $(C,0)$  with trivial section $\sigma$ is \textbf{equisingular along} $\sigma$ or an \textbf{es--deformation with section} of $(C,0)$ if
\begin{enumerate}
\item [(i)] it is equimultiple along $\sigma$, i.e. $\ord_{(x,y)}(F)=\mt(C,0)$,
\item [(ii)] after blowing up $\sigma$, there exist sections through the infinitely near points in the first neighbourhood of $(C,0)$ such that the respective reduced total transform of $(\kc, 0)$ is equisingular along these sections,
\item [(iii)] for a nodal curve $(C,0)$ (defined by $f(x,y)=xy$) an es--deformation is an equimultiple deformation.
\end{enumerate}
\end{Definition}

\noindent For an es--deformation of $(C,0)$ with section the section is unique (cf. \cite{GLS07} Prop.II.2.8). We call a deformation $\Phi$ (without section) of $(C,0)$ \textbf{equisingular} or an \textbf{es--deformation} if there exists a section $\sigma$ such that $\Phi$ is equisingular along $\sigma$. We denote the category of $es$--deformations of $(C,0)$ by $\Def_{(C,0)}^{es}$, which is a full subcategory of $\Def_{(C,0)}$. If $\varphi$ is a parametrization of $(C,0)$ we denote es--deformations (with section) of $\varphi$ by $\Def^{es}_\varphi$. The base spaces of the corresponding semiuniversal deformations are denoted by
\[
(B_C, 0)\ ,\ (B_C^{es}, 0) \text{ and } (B_\varphi^{es}, 0).
\]

With these notations we have the main result of this section about es-deformations of isolated plane curve singularities.

\begin{Theorem}  {\bf (Campillo, Greuel, Lossen)}\label{Theo4.6}
\begin{enumerate}
\item [(1)] Every es--deformation of $\varphi$ induces a unique es--deformation of $(C,0)$, providing  a functor $\Def_\varphi^{es}\to \Def_C^{es}$.
\item [(2)] Every es--deformation of $(C,0)$ comes from an es--deformation of $\varphi$, i.e. $\Def_\varphi^{es}\to D_C^{es}$ is surjective.
\item [(3)] Any two es--deformations of $\varphi$ that induces isomorphic es--deformation of $(C,0)$ are isomorphic, i.e. the induced transformation of isomorphism classes $\underline{\Def}_\varphi^{es}\to \underline{\Def}_C^{es}$ is an isomorphism of functors.
\item [(4)] There exists a morphism of base spaces $\alpha:(B_\varphi^{es}, 0)\to (B_C, 0)$ mapping $(B_\varphi^{es}, 0)$ ismorphically onto $(B_C^{es}, 0)$.
\item [(5)] $(B_C^{es}, 0)$ coincides with the $\mu$--constant stratum $(\Delta^\mu, 0)\subset (B_C, 0)$.
\end{enumerate}
\end{Theorem}
\medskip

\begin{Corollary}
The $\mu$--constant stratum $(\Delta^\mu, 0)$ of the semiuniversl deformation of a reduced plane curve singularity is smooth of dimension $T^ {1, es}_\varphi$.
\end{Corollary}

For the proof see \cite{GLS07}, Theorem II.2.64 and Theorem II.2.61. The proof is by induction on the number of blowing ups and uses the notion of equi-intersectional deformations over arbitrary (not necessarily reduced) base spaces, introduced in \cite{GLS07}.

\begin{Remark}
\rm{\begin{enumerate}
\item [(1)] The first proof of the smoothness of the $\mu$--constant stratum for a reduced plane curve singularity was given by Jonathan Wahl in \cite{Wa74}, using deformation theory of global objects, namely of divisors supported on the exceptional divisor of the resolution of $(C,0)$.
Our proof, first published in \cite{CGL07a} and \cite{GLS07}, using es--deformations of the parametrization, appears to be simpler than Wahl's.
\item [(2)] An algorithm to compute the $\mu$--constant stratum is given in \cite{CGL07b}, Algorithm 2 and Remark 3.12. It is implemented in the {\sc Singular} library {\tt equising.lib} \cite{LM07}.
\end{enumerate}}
\end{Remark}
\bigskip

\section{A $\delta$-- and $\mu$--invariant for isolated non--normal singularities}

We introduce in this section a $\delta$--invariant and a $\mu$--invariant for isolated non--normal singularities. Moreover, we fix some notations to be used in the rest of this article. Most definitions and results of this section have been given in \cite{BrG90} for generically reduced curves. The generalization to isolated non-normal singularities given here is straightforward.
\medskip

For a complex space $X$ let $X^{\red}$ denote its reduction and $i:X^{\red}\hookrightarrow X$ the closed embedding. $X$ and $X^\red$ coincide as topological spaces but for the structure sheaves we have $\ko_{X^\red}=\ko_X/Nil(\ko_X)$, where $Nil(\ko_X)$ is the sheaf of nilpotent elements of $\ko_X$. 
If $n:\overline{X}\to X^\red$ denotes the normalization of $X^\red$, we call the composition
\[
\nu:\overline{X}\xrightarrow{n} X^\red\overset{i}{\hookrightarrow} X
\]
the \textbf{normalization} of $X$. $X$ is \textbf{normal} if $\nu$ is an isomorphism. 

A normal space is reduced and \textbf {locally irreducible}, i.e. the germs $(X,x)$ are irreducible for all $x \in X$. A germ is called locally irreducible if this holds for some representative.
\medskip

Moreover, we denote by $w:\widehat{X}\to X^\red$ the weak normalization of $X^\red$ and by 
\[
\omega: \widehat{X}\xrightarrow{w} X^\red\hookrightarrow X
\]
the \textbf{weak normalization} of $X$. $X$ is weakly normal if $\omega$ is an isomorphism. 

We call a morphism 
\[
\pi: \widetilde {X}\to X  
\]
a \textbf{partial normalization} of $X$ if $\pi$ factors as 
$\pi: \widetilde{X} \to X^\red\hookrightarrow X$, and if $\nu$ factors as 
$\nu:\overline{X}\to \widetilde{X}\xrightarrow{\pi} X$. 
Examples are the normalization, the reduction and the weak normalization of $X$.

For a morphism $f:X\to T$ of complex spaces we denote by
\[
\begin{array}{rcl}
 f^\red  :  X^\red\to T^\red, \\ 
\bar{f}=f\circ \nu: \overline{X} \to T, \\
\widehat{f}=f\circ \omega: \widehat{X} \to T, \\
\widetilde{f}:=f\circ\pi:\widetilde{X}\to T
\end{array}
\]
the induced morphisms and by 
\[
X_t=f^{-1}(t),\ X^\red_t=(f^\red)^{-1}(t),\  \overline{X}_t=\bar{f}^{-1}(t),\ \widehat{X}_t=\widehat{f}^{-1}(t),\
\widetilde{X}_t=\widetilde{f}^{-1}(t) 
\]
the corresponding fibres over $t \in T$.\\

Recall that $X$ is called \textbf{weakly normal} (or \textbf{maximal}) if every continuous function on $X$, which is holomorphic on the regular locus of $X$, is holomorphic on $X$.

$\widehat{X}$ is reduced and obtained from the normalization $\overline{X}$ by identifying points with the same image in $X$. $w$ is a homeomorphism (cf. \cite{Fi76}, Section 2.29). A germ $(X,x)$ is called normal resp. weakly normal if this holds for some representative. 
\medskip

A curve singularity of embedding dimension $r$ is weakly normal iff it consists of $r$ smooth branches with independent tangent directions, i.e. it is isomorphic to the coordinate axes in $\C^r$. The Whitney umbrella $\{y^2-x^2z=0\}$ is an example of a weakly normal surface (cf. \cite{Vi11}); it has a smooth normalization.
\medskip

By an \textbf{irreducible component} of a complex space or a germ we always mean an irreducible component of its {\em reduction}. $(X,x)$ is called \textbf{irreducible} if $(X^\red, x)$ is irreducible (note that an irreducible germ need not be reduced).  If $(X,x)$ is a curve singularity (i.e. a $1$--dimensional germ), an irreducible component of $(X,x)$ is also called a \textbf{branch} of $(X,x)$.
\medskip

If a complex germ $(X,x)$ has $r$ irreducible components $(X_i, x), i=1, \ldots, r$, its normalization is given by
\[
\nu:(\overline{X}, \overline{x})=\coprod\limits_{i=1}^{r} (\overline{X}_, \overline{x}_i)\xrightarrow{n}(X^\red, x)\overset{i}{\hookrightarrow}(X,x)\ ,
\]
with $(\overline{X}, \overline{x})$  a multigerm, $\overline{x}=n^{-1}(x)$, and $(\overline{X}_i, \overline{x}_i)\to (X_i, x)$ the normalization of $(X_i, x)$.
$(X,x)$ is normal iff $\nu:(\overline{X}, \overline{x})\to (X,x)$ is an isomorphism.

It is convenient to introduce
\[
\begin{array}{rcl}
\mathbf{r'(X,x)} & = & r(X,x)-1,\\
\end{array}
\]
if $r(X,x)$ denotes the number of irreducible components of $(X,x)$.
\medskip

If $X$ is a complex space and $\gamma(X,x)\in \Z\cup\{\infty\}$ an invariant defined for the germs $(X,x)$ such that $\gamma(X,x)\neq 0$ for only finitely many $x\in X$, we write
\[
\gamma(X):=\sum\limits_{x\in X}\gamma(X,x).\
\]

\begin{Definition} \rm \label{def-ri} 
\begin{enumerate}
\item[(1)] For a complex space $X$, resp. a germ $(X,x)$, we denote by $X^i$, resp. by $(X^i, x)$ the union the $i$--dimensional irreducible components of $X^\red$, resp. of $(X^\red, x)$. We set
\[
X^{>j}:=\underset{i>j}{\cup} X^i\ ,\ (X^{>j}, x):=\underset{i>j}{\cup}(X^i, x), \text{and}
\]
\[
\begin{array}{c}
r_i(X):=\sharp\{\text{ irreducible componentes of } X^i\},\\
r_i(X,x):=\sharp \{\mbox{irreducible components of } (X^i,x)\}.
\end{array}
\]
The latter is the same as the number of $i$--dimensional minimal primes of $\ko_{X,x}$.

\item[(2)] We say that $x$ is an isolated non--normal point of $X$ or $(X,x)$ is an \textbf{isolated non--normal singularity} (\textbf{INNS} for short) if there exists a neighbourhood $U$ of $x$ such that $U\smallsetminus \{x\}$ is normal.
\end{enumerate}
E.g. any isolated singularity is an INNS.
\end{Definition}

Let $\nu:\overline{X}\xrightarrow{n} X^\red\overset{i}{\hookrightarrow} X$ be the normalization of  the complex space $X$ and 
$$\nu^{\sharp}: \ko_X \sr \nu_*\ko_{\gt{X}} $$ 
the induced map of the structure sheaves. We have $\Ker(\nu^\sharp)=Nil(\ko_X)$ and $\Coker(\nu^\sharp)=\nu_*\ko_{\gt{X}}/\ko_{X^{red}}$. Since $\nu$ is finite, $\Ker(\nu^\sharp)$ resp. $\Coker(\nu^\sharp)$ are coherent sheaves of $\ko_X$-modules with support the non-reduced locus resp. the non-normal locus of $X$. 
\medskip

If $x\in X$ is an INNS, then $\Ker(\nu^\sharp)_x$ and $\Coker(\nu^\sharp)_x$ are finite dimensional vector spaces and, if $\dim(X,x)\geq 1$,  we have $\Ker(\nu^\sharp)_x=\Nil(\ko_{X,x})=H^0_{\{x\}}(\ko_X)$,  where $H^0_{\{x\}}$ denotes local cohomology. If $\dim(X,x)=0$ we have  $\dim_\C H^0_x(\ko_X)=\dim_\C\Nil(\ko_{X,x})+1$.
\bigskip

For the definition of $\delta$ the isolated points of $X$ play a special  role. Let 
$$ X^{>0}=X^{red}\tru \{\mbox{isolated points of } X^{red}\},  $$
be the positive dimensional part of $X^\red$, and denote by $\nu: \gt{X^{>0}} \mtn X^{>0}$ the normalization, 
(with $\ko_{\gt{{X}^{>0}}}=\ko_{X^>0}=0$ if $\dim X=0$).

\begin{Definition} \label{Def5.2} \rm 
Let $(X,x)$ be an INNS. We define
 $$\epsilon(X,x):=\dim_\C H^0_{\{x\}}(\ko_{X}),$$
   the \textbf{epsilon-invariant}, 
$$ \delta(X,x):= \dim_\C \big((\nu_*\ko_{\gt{X^{>0}}})_x/\ko_{X^{>0},x}\big) - \epsilon(X,x)$$
the \textbf{delta-invariant}, and
$$ \mu(X,x)=2\delta(X,x)-r'(X,x)$$
the $\mathbf{\mu}$\textbf{--invariant} of $(X,x)$.
\end{Definition}

\begin{Remark}\label{Rem5.3}\rm 
\begin{enumerate}
\item[(1)] The above definition coincides with the classical $\delta$--inva\-riant resp. Milnor number $\mu(X,x)$ if $(X,x)$ is a reduced curve singularity. It generalizes $\delta$ and $\mu$ defined in \cite{BrG90} for $\dim (X,x)\leq 1$.  We have 
$$ \delta(X,x)=\delta(X^{>0},x)-\epsilon(X,x), $$
with $\delta(X^{>0},x)\geq 0$  and $\delta(X^{>0},x)=0$ if and only if $X^{>0}$ is normal at $x$ or $x\not \in X^{>0}$.
\item[(2)] If $\dim (X,x)=0$ then $\delta(X,x)=-\dim_\C \ko_{X,x}=-\epsilon(X,x)$. In particular, $\delta(X,x)=-1$ for $x$  an isolated and reduced (hence normal) point of $X$.
\item[(3)] Let $\dim (X,x) >0$. Then $\delta(X,x)=\delta(X^\red, x)-\varepsilon(X,x)$.  If $(X,x)$ is normal, then $\delta(X,x) = \epsilon(X,x)= \mu(X,x) =0$.
We have $ \mu(X^\red, x)\geq \delta(X^{red},x) \geq 0$, and $\delta(X^{red},x) = 0 \Leftrightarrow \mu(X^\red, x)=0 \Leftrightarrow (X^{red},x)$ is normal.  However  we may have
 $\delta < 0$ and $\mu<0$ if $(X,x)$ is not reduced and $(X^{red},x)$ is not normal (e.g. $\delta(X,0) = -k$ for  the irreducible curve singularity $(X,0)$ defined by the ideal 
$\langle{x^2-y^3}\rangle \cap \langle{x^{3+k},y}\rangle , k\geq 0$).

\item[(4)] We observe that for any INNS we have that $-\delta(X,x)$  is the Euler characteristic $\chi ({\ko}^\bullet_{X,x})=\dim_\C H^0({\ko}^\bullet_{X,x})-\dim_\C H^1({\ko}^\bullet_{X,x})$ of the two term complex 
$${\ko}^\bullet_{X,x}: 0 \mtn \ko^{0}_{X,x}:= \ko_{X,x} \mtn \ko^{1}_{X,x}:=(\nu_*\ko_{\gt{X^{>0}}})_x \mtn 0, $$
which turns out to be useful in proofs.

\item [(5)] If $\dim(X,x)\geq 2$, we prefer the notion $\mu$--invariant instead of Milnor number, since its properties are then rather different from the classical Milnor number for isolated hypersurface or complete intersection singularities. 
\end{enumerate}
\end{Remark}
\medskip

For topological considerations the weak normalization turns out to be useful. An example of a weakly normal germ $(X,x)$ is the $1$--point union (or bouquet) of normal singularities in independent directions. That is, if $(X_i,0)\subset (\C^{n_i}, 0), i=1, \ldots, r$ are normal germs, then $(X,0)\subset(\C^{n_1}\times \cdots\times\C^{n_r}, 0)$ is the union of the $(X_i,0)$, embedded in the $i$--th component $(0, \ldots, 0, \C^{n_i}, 0\ldots, 0)$. Any singularity isomorphic to $(X,0)$ is called an \textbf{ordinary INNS}.

\begin{Lemma}\label{lemma8.1} Let $(X,x)$ be an INNS, with irreducible components $(X_i, x), i=1, \ldots, r$, \ and $(\widehat{X}, \widehat{x})$ the weak normalization of $(X,x)$. The following holds for good representatives.
\begin{enumerate}
\item [(1)] $X_i\cap X_j=\{x\}$ for $i\neq j$.
\item [(2)] The weak normalization $\widehat{X}\to X$ is a homeomorphism.
\item [(3)] The irreducible components $\widehat{X}_i$ of $\widehat{X}$ are normal and $\widehat{X}_i\to X_i$ is the normalization of $X_i$ for $i=1, \ldots, r$.
\item [(4)] $(\widehat{X}, \widehat{x})$ is an ordinary INSS.
\end{enumerate}
\end{Lemma} 

\noindent\textbf{Proof:} (1) follows since $x$ is an isolated non--normal point in $X$ and $X_i\cap X_j$ is non--normal in any of its points, while (2) is a general property of the weak normalization (cf. \cite{Fi76}, section 2.29). Since $X\smallsetminus\{x\}$ is normal, the germ $(\widehat{X}, \widehat{x})$ is constructed from the normalization $(\overline{X}, \overline{x})=\coprod\limits_{i=1}^r(\overline{X}_i, \overline{x}_i)$ by identifying the $(\overline{X}_i, \overline{x}_i)$ in one point $\widehat{x}$. This follows from the proof of the maximalization theorem in \cite{Fi76}, Section 2.29 (see also the proof of the next lemma) and proves (3) and (4).\hfill\qed \\

\begin{Lemma}\label{lemma8.2} Let $(X,x)$  be a reduced germ of dimension $\geq$ 1.  Then $\delta(X,x) \geq r'(X,x)$ and the following are equivalent:
\begin{enumerate}
	\item [(i)] $(X,x)$ is a weakly normal INNS,
	\item [(ii)] $(X,x)$ is an ordinary INNS,
	\item [(iii)] $\delta(X,x)=r'(X,x)$.
\end{enumerate}
\end{Lemma}

\noindent\textbf{Proof:} 
For $X$ weakly normal and $U\subset X$ open we have $\ko_X(U)\subset n_\ast\ko_{\overline{X}}(U)$ with
\[
\begin{array}{c}
\ko_X(U)=\{f:U\to\C \text{ continuous}, f|U\smallsetminus\Sing(X) \text{ holomorphic}\},\\ n_\ast\ko_{\overline{X}}(U)=\{f:U\smallsetminus\Sing(X)\to \C \text{ holomorphic and locally bounded on }U\}. 
\end{array}
\]

To prove the equivalences of the lemma, let $(X,x)$ be an INNS. Since $X\smallsetminus \{x\}$ is normal, we have $(n_\ast\ko_{\overline{X}})(U\smallsetminus\{x\})=\ko_X(U\smallsetminus\{x\})$ 
and $U\smallsetminus \{x\}=\coprod\limits_{i=1}^rU_i\smallsetminus \{\overline{x}_i\}$, 
with $X_i$ the irreducible components of $X$ and $U_i\subset \overline{X_i}$ an open neighbourhood of $\overline{x}_i$. 
Hence, for $U$ sufficiently small,
\[
n_\ast\ko_{\overline{X}}(U)=\{f:U\smallsetminus\{x\}\to\C\text{ holomorphic and bounded in }x\}.
\]
Therefore the Jacobson radical $J$ of $(n_\ast\ko_{\overline{X}})_x$, represented by $f \in n_\ast\ko_{\overline{X}}(U)$ with $f(\overline{x}_i)=0$ for all $i$, is contained in $\ko_{X,x}$ and we get, since $(X,x)$ is reduced,
\[
\delta(X,x) = \dim_\C (n_\ast\ko_{\overline{X}})_x/\ko_{X,x}=r-1
\]
(the constants of $\ko_{X,x}$ are diagonally embedded in $(n_\ast\ko_{\overline{X}})_x$).\\
The argument shows also that $(X,x)$ is an ordinary INNS and that $\delta(X,x) \geq r-1$ with equalitiy iff $J$ is contained in $\ko_{X,x}$, i.e. iff $(X,x)$ is an ordinary INNS. This proves the lemma.\hfill\qed \\

The following definition is quite useful in our context. It should however not be confused with the intersection multiplicity defined by Serre in \cite{Se65}.
\medskip

For two complex germs $(X', 0), (X'', 0)\subset(\C^n, 0)$ defined by ideals $I', I''\subset \ko_{\C^n, 0}$ we define the \textbf{intersection number} as
\[
(X', X'')_0:=\dim_\C\ko_{\C^n, 0}/(I'+I''),
\]
which is finite iff $(X', 0)\cap(X'', 0)=\{0\}$. \\

\begin{Proposition}\label{Prop5.4}
Let $(X_i, 0)\subset (\C^n, 0), i=1, \ldots, p$, be isolated non--normal singularities, defined by ideals $I_i\subset \ko_{\C^n, 0}$ such that, for $i\neq j$, $(X_i, 0)\cap (X_j, 0)=\{0\}$.  We set $(X,0):=\overset{p}{\underset{i=1}{\cup}}(X_i, 0)$ defined by $I:=\overset{p}{\underset{i=1}{\cap}}I_i.$

Then $(X,0)$ is an INNS and we have:
\begin{enumerate}
\item [(i)]$\delta(X,0)=\overset{p}{\underset{i=1}{\sum}}\delta(X_i, 0)+\overset{p}{\underset{i=1}{\sum}}(X_i, X_{i+1}\cup\cdots\cup X_p)_0$,

\item [(ii)] 
$\varepsilon(X,0)=\overset{p}{\underset{i=1}{\sum}}\varepsilon(X_i, 0)  - \overset{p}{\underset{i=1}{\sum}}(X_i,X_{i+1}\cup\cdots\cup X_p)_0\\
\hspace*{3.45cm}       + \overset{p}{\underset{i=1}{\sum}}(X_i^{>0},X_{i+1}^{>0}\cup \cdots\cup X_p^{>0}))_0
$

\end{enumerate}

In  particular, $\varepsilon(X,0)=\overset{p}{\underset{i=1}{\sum}}\varepsilon(X_i, 0)$ if $\dim(X_i,0) > 0$ for ${i, \ldots,p}$.

Here the union of germs is defined by the intersection of the corresponding ideals.
\end{Proposition}

\noindent\textbf{Proof}: (c.f. \cite{BrG90} for isolated curve singularities). The proof is by induction on $p$.

For $p=2$ we have to show
\[
\begin{array}{lcl}
\delta(X) & = & \delta(X_1)+\delta(X_2)+(X_1, X_2)_0,\\
\varepsilon (X) & = & \varepsilon(X_1)+\varepsilon(X_2)-(X_1, X_2)_0+(X_1^{>0}, X_2^{>0})_0\ .
\end{array}
\]
Consider, with $\ko :=\ko_{\C^n, 0}$, the commutative diagram 
\[
\xymatrix{0\ar[r] & \ko/I\ar[d]\ar[r] & \ko/I_1\oplus\ko/I_2\ar[d]\ar[r] & \ko/(I_1+I_2)\ar[d]^\gamma\ar[r] & 0\\
0\ar[r] & \ko/I^{>0}\ar[r] & \ko/I_1^{>0}\oplus \ko/I_2^{>0}\ar[r] & \ko/(I_1^{>0}+I_2^{>0})\ar[r] & 0\ ,
}
\]
with exact rows and surjective vertical arrows. Note that $I^{>0}=\sqrt{I}$ if $\dim(X,0) >0$ and $I^{>0}=\ko$ if $\dim(X,0)=0$. In any case we get the exact kernel sequence
\[
0\to H^0_x(\ko_X)\to H^0_x(\ko_{X_1})\oplus H^0_x(\ko_{X_2})\to \Ker(\gamma)\to 0\ ,
\]
with $\dim_\C\Ker(\gamma)=(X_1, X_2)_0-(X_1^{>0}, X_2^{>0})_0$. This proves (ii). 

Similarly we have a diagram with exact rows
\[
\xymatrix{0\ar[r] & \ko/I\ar[d]^\alpha\ar[r] & \ko/I_1\oplus\ko/I_2\ar[r]\ar[d]^\beta & \ko/I_1+I_2\ar[d]\ar[r] & 0\\
0 \ar[r] & \overline{\ko/I^{>0}}\ar[r] & \overline{\ko/I_1^{>0}}\oplus \overline{\ko/I_2^{>0}}\ar[r] & 0 ,
}
\]
where $\overline{\phantom{xx}}$ denotes integeral closure. Note that $\alpha$ is the morphism of the complex $\ko^{\bullet}_{X,0}$ and $\beta$ that of $\ko^{\bullet}_{X_1, 0}\oplus\ko^\bullet_{X_2, 0}$. Taking Euler characterstics, we get
\[
\chi(\ko^{\bullet}_{X_1, 0})+\chi(\ko^\bullet_{X_2,0})=\chi(\ko^\bullet_{X,0})+(X_1, X_2)_0
\]
which implies (i) by Remark \ref{Rem5.3} (4).

The case $p>2$ follows by induction.\hfill\qed 
\\

In Proposition \ref{Prop5.4} it is not assumed that the $X_i$ are irreducible. It may even happen that $I_i\subset I_j$ for some $i\neq j$, i.e. $X_j\subset X_i$, which implies that $I_j$ is an $\frak{m}$--primary ideal, with $\frak{m}\subset \ko_{\C^n,0}$ the maximal ideal.\\

Consider now the primary decomposition of an ideal $I\subset \ko_ {\C^n,0}$ defining the INNS $(X,0)\subset(\C^n, 0)$, \[
I=\overset{s}{\underset{i=1}{\cap}} Q_i, 
\]
with $Q_i$ being $P_i$--primary, $P_i=\sqrt{Q_i}$, and $\{P_1, \ldots, P_s\}$ the set of associated prime ideals of $I$.

Let $\{P_1, \ldots, P_r\}, r\leq s$, be the set of minimal primes and $(X_i, 0)$ the irrducible components of $(X,0)$ defined by the minimal prime $P_i$, $i=1,\ldots,r$. We have $(X_i, 0)\cap(X_j, 0)=\{0\}$ for $i\neq j$, since any point of $X_i \cap X_j$ is a non--normal point of $X$. Moreover, since $X\smallsetminus \{x\}$ is reduced, there is only one embedded prime, the maximal ideal $\frak{m}$. Since $(X^\red, x)$ is defined by $\sqrt{I}=\overset{r}{\underset{i=1}{\cap}} P_i$, we get the primary decomposition of $I$ as in the following Corollary.

\begin{Corollary}\label{Cor5.5}
Let $(X,0)\subset(\C^n, 0)$ be an INNS, $\dim(X,0)\geq 0$, defined by the ideal  $I\subset \ko_{\C^n, 0}$. Then $I$ has an irredundant primary decomposition of the form
\[
I=P_1\cap\cdots\cap P_r\cap Q,
\]
with prime ideals $P_i$, defining the irreducible components $(X_i, 0), i=1, \ldots, r$ of $(X,0)$, and a (not unique) $\frak{m}$--primary ideal $Q$, defining the embedded component of $(X,0)$. We can use (i) resp. (ii) of Proposition \ref{Prop5.4} to compute $\delta(X,0)$ resp. $\varepsilon(X,0)$. In particular:
\begin{enumerate}
\item [(i)] $\delta(X^{>0}) =\overset{r}{\underset{i=1}{\sum}}\delta(X_i, 0)+\overset{r}{\underset{i=1}{\sum}}(X_i, X_{i+1}\cup\cdots\cup X_r)_0$,
	
\item [(ii)] $\delta(X,0)=\delta(X^{>0}, 0)-\varepsilon(X,0)$,

\item [(iii)] $\varepsilon(X,0)=\dim_\C\ko_{\C^n,0}/Q-\dim_\C\ko_{\C^n, 0}/(I^{>0}+Q)$.
\end{enumerate}
Here $I^{>0}$ defines $(X^{>0}, 0)$, i.e. $I^{>0}=\sqrt{I}=P_1\cap\cdots\cap P_r$ if $\dim(X,0)>0$ and $I^{>0}=\ko_{\C^n, 0}$ if $\dim(X,0)=0$.
\end{Corollary}
\medskip

\begin{Example}\label{Ex5.6} \rm 
Let $(X,0)\subset\C^3, 0)$ be the union of an $A_3$--singularity $X_1$ in the $(x,y)$--plane with the $z$--axis $X_2$ and an embedded component $X_3$, defined by the respective ideals $I_1=\langle z, x^2-y^4\rangle$, $I_2=\langle x,y\rangle$ and $Q=\langle z, x^3, y^5\rangle$ in $\ko=\ko_{C^3, 0}$.

We have $I^{>0}=I_1\cap I_2=\langle xz, yz, x^3-y^4\rangle, I^{>0}+Q=\langle z, x^2-y^4, xy^4, y^5\rangle$ and $I_1+I_2=\langle x,y,z\rangle$. Then
\[
\begin{array}{lcl}
\delta(X^{>0}) & = & \delta(X_1)+\delta(X_2)+\dim_\C\ko/(I_1+I_2)=2+0+1=3,\\
\varepsilon(X) & = & \dim_\C\ko/Q-\dim_\C\ko/(I^{>0}+Q)=15-9=6,
\end{array}
\]
and hence $\delta(X)=-3$.

Note the $I$ has also a different primary decomposition $I_1\cap I_2\cap Q_1$ with $Q_1=\langle z, x^4, x^3y, y^5\rangle$. Then $\dim_\C\ko/Q_1=16$ and $\dim \ko/(I^{>0}+Q_1)=10$, giving of course the same $\varepsilon(X)$.
\medskip

Proposition \ref{Prop5.4} is very useful for practical computations. We provide the {\sc Singular} code for the above examples. We compute $\varepsilon(X)$ in two different ways: by Corollary \ref{Cor5.5} (ii) and as $\dim_\C(\sqrt{I}/I)$. Since the computation of the dimension for a finite dimensional quotient $J/I$ is not implemented in {\sc Singular}, we provide a small procedure for computing $\dim_\C(J/I)$.

\begin{verbatim}
LIB "all.lib";
ring R = 0,(x,y,z),ds;  
ideal I1 = z,x2-y4;                   
ideal I2 = x,y; 
ideal II = intersect(I1,I2);          // defines X^red
ideal Q  = z,x3,y5;                   // defines emb. component
ideal I  = std(intersect(II,Q));      // defines X

list nor = normal(I,"withDelta"); 
nor[size(nor)];                       // displays delta(X^red) 
// ->  3 = delta(X^red)
vdim(std(Q));                         // dim_C(R/Q) = 15
vdim(std(II+Q));                      // (II,Q)_0   = 9 
// epsilon(X) =  15 - 9 = 6 
// delta(X) = 3 - epsilon(X) = -3

primdecGTZ(I);                        // primary decomp. of I
// shows that I has also Q1 = <z,y5,x3y,x4> as emb. component
ideal Q1 = z,y5,x3y,x4;
vdim(std(Q1));                        // dim_C(R/Q1) = 16
vdim(std(II+Q1));                     // (II,Q1)_0   = 10
// epsilon(X) = 16 - 10 = 6 

// If dim(X)>0, epsilon is the C-dimension of radical(I)/I 
// See procedure Vdim below to compute radical(I)/I 
ideal rI = radical(I);
Vdim(rI,I);                           // epsilon(X)  = 6

proc Vdim (ideal i, ideal j) 
//computes vectorspace dimension of i/j, requires poly.lib
{
   j = lead(std(j));
   i = lead(std(i));
   int d = maxdeg1(j);
   j = j,maxideal(d+1);
   i = i,maxideal(d+1);
   int ki = vdim(std(i));
   int kj = vdim(std(j));
   return (kj-ki);
}
\end{verbatim}

{\em Warning}: If $I\subset \Q[X_1, \ldots, X_n]$ is an ideal, then the procedure {\tt normal (I,"withDelta")} computes the total $\delta$--invariant $\delta(X)$ of the affine variety $X=V(I)\subset\C^n$, even if we compute w.r.t. a local ordering. It computes $\delta(X,0)$ if $0$ is the only non--normal point of $X$ (as in our case).
\end{Example}
\bigskip

\section{Simultaneous normalization}

In this section we introduce the general concepts of simultaneous normalizations and equinormalizability and prove some consequences for maps that admit a simultaneous normalization. Moreover, we recall the main results from Chiang--Hsieh and Lipman in \cite{CL06} and
Koll$\acute{\text{a}}$r in \cite{Ko11} about simultaneous normalization of projective morphisms.
\bigskip

Let $f: X \mtn S$ be a morphism of complex spaces. We call $f$ {\bf regular} (resp. {\bf normal}, resp. {\bf weakly normal}, resp. {\bf reduced}) at $x\in X$ if $f$ is flat at $x$ and if the fibre $X_{f(x)}$ is regular (resp. normal, resp. weakly normal, resp. reduced) at $x$. The morphism $f$ is called regular (resp.  normal, resp. weakly normal, resp. reduced) if this holds at every $x\in X$. Recall that a complex space $X$ is said to be weakly normal  (or maximal) if every continuous function on $X$ which is holomorphic on the regular locus of $X$ is holomorphic on $X$.\\
The following sets
\begin{align}
\mbox{Flat}(f) & = \{x\in X | f \mbox{ is flat in  } x\},\nonumber\\
\mbox{Reg}(f) &= \{x\in X | f \mbox{ is regular in  } x\},\nonumber\\
\mbox{Nor}(f) &= \{x\in X | f \mbox{ is normal in  } x\},\nonumber\\
\mbox{Wor}(f) &= \{x\in X | f \mbox{ is weakly normal in  } x\},\nonumber\\
 \mbox{Red}(f) &= \{x\in X | f \mbox{ is reduced in  } x\}\nonumber
\end{align}
are analytically open subsets of $X$. For Flat cf. \cite{Fr67}, Th\'{e}or\`{e}me (IV,9), and for Reg, Nor,  Wor, Red see \cite{BF93}, Theorem 1.1 (2), Theorem 2.1 (2).\\
$\NFlat(f), \Sing(f), \NNor(f), \NWor(f), \NRed(f)$ denote the corresponding ana\-ly\-tic sets of  non-flat, non-regular or singular, non-normal, non-weakly normal, non-reduced  points of $f$. \\
 We denote by  $\mbox{Nor}(X)$, resp. $\mbox{Wor}(X)$,  resp. $\mbox{Red}(X)$ the set of  normal, resp. weakly normal,  resp. reduced points of $X$ and with an N in front the corresponding complements.  The set of regular points of $X$ is denoted by $\mbox{Reg}(X)$, and its complement by $\Sing(X)$. Let us recall the following well-known result:

\begin{Proposition} \label{Prop6.1} Let $f: X \mtn S$ be a morphism of complex spaces.
\begin{enumerate}
\item[(1)] If $f$ is reduced then $X$ is reduced if and only if $f(X)\subseteq \mbox{Red}(S)$. 
\item [(2)] 
\begin{enumerate} 
 \item [(i)] If $f$ is regular then $X$ is regular if and only if $f(X)\subseteq \mbox{Reg}(S)$.
 \item [(ii)] If $X$ is reduced and $f$ flat, then Reg$(f)$ is dense in $X$.
\end{enumerate} 
\item [(3)] If $f$ is normal then $X$ is normal (resp. weakly normal) if and only if $f(X)\subseteq \mbox{Nor}(S)$ (resp. $f(X)\subseteq \mbox{Wor}(S)$).
\end{enumerate}
\end{Proposition}

For the proof of (1) and (3) we refer to \cite{Mat86}, Corollary to Theorem 23.9, \cite{Ni81}, (2.4), and \cite{Man80}, (III.3). For (2) see \cite{GLS07} (Definition II.1.112 and Theorem II.1.115 for (i), and Corollary II.1.116 for (ii)).\\

We need also the following theorem.

\begin{Theorem} \label{thr-neglectible} Let $f: X \mtn S$ be a morphism of complex spaces.
\begin{enumerate}
\item[(1)] If $S$ is generically\footnotemark\footnotetext{We say that a property $\kp$ holds \emph{generically } on $X$, if the set of points where $\kp$ holds, contains an analytically open dense subset of $X$, i.e. the complement is analytic and nowhere dense in $X$. } reduced,  then $f(\NFlat(f))$ is neglectible\footnotemark \footnotetext{A subset $N$ of a complex space $X$ is called \textit{neglectible}, if $N$ is contained in a countable union of nowhere dense locally analytic subsets of $X$. Then $X\setminus N$ is dense in $X$ by Baire's theorem.} in $S$.
\item[(2)] Assume there is an open dense subset $V\subseteq S$ with $V\subseteq \mbox{Reg}(S)$.
\it
\item[(i)] If $f^{-1}(V)\subseteq \mbox{Reg}(X)$,  then $f(\Sing(f))$ is neglectible in $S$.
\item[(ii)]  If $f^{-1}(V)\subseteq \mbox{Nor}(X)$,  then $f(\NNor(f))$ is neglectible in $S$.
\item[(iii)]  If $f^{-1}(V)\subseteq \mbox{Wor}(X)$,  then $f(\NWor(f))$ is neglectible in $S$.
\item[(iv)]  If $f^{-1}(V)\subseteq \mbox{Red}(X)$,  then $f(\NRed(f))$ is neglectible in $S$.
\end{enumerate}
If the restriction of $f$ to $\Sing(f)$ (resp. $\NNor(f)$, resp. $\NWor(f)$, resp. $\NRed(f)$) is proper, then in statements (i) to (iv), ''neglectible'' can be replaced by ''a nowhere dense closed analytic subset''.
\end{Theorem}
Statement (1) is due to Frisch (\cite{Fr67}), Prop. IV. 14). Statement (2)(i) is Sard's theorem, while the rest is proved in \cite{BF93}, Theorem 2.1(3).
\bigskip

Following \cite{CL06} and \cite{Ko11} we make the following

\begin{Definition}   \label{defEquinormal}  Let $f: X\sr S$ be a morphism of complex spaces. 
\begin{enumerate}
\item [(1)] A \textbf{simultaneous normalization of $f$} is a morphism
$n: \nga{X} \sr X$ such that
 \begin{enumerate}
\item[(i)]  $n$ is finite,
\item[(ii)] $\nga{f}:=f\circ n: \nga{X}\mtn S$ is normal,
\item[(iii)] the induced map $n_s: \nga{X}_s:=\nga{f}^{-1}(s) \sr X_s$ is bimeromorphic for each $s\in f(X)$.
\end{enumerate}

\item [(2)] $f$ is called \textbf{equinormalizable} if  the normalization $\nu: \gt{X}\mtn X$ is a simultaneous normalization of $f$.
$f$ is called equinormalizable at $x\in X$ if the restriction of $f$ to some neighborhood of $x$ is equinormalizable.

\item [(3)] We say that $n: \nga{X} \sr X$ is a  \textbf{simultaneous weak normalization of $f$} if it satisfies (i) and (iii) of (1) and if $\nga{f}: \nga{X}\mtn S$ is weakly normal.

\end{enumerate}
\rm{If $f: (X,x) \sr (S,s)$ is a morphism  of germs, then a simultaneous normalization of $f$ is a morphism $n$ from a multi-germ $(\nga{X}, n^{-1}(x))$ to $(X,x)$ such that some representative of $n$ is a simultaneous normalization of a
representative
of $f$. The germ $f$ is called equinormalizable if some representative of $f$ is equinormalizable.}
 \end{Definition}

\begin{Remark} \label{rem1.5}\rm

\begin{enumerate}
\item [(1)] $n_s:\widetilde{X}_s\to X_s$ is finite, hence proper. A proper bimeromorphic morphism is surjective since the image is closed and contains an open dense subset. It follows that $\widetilde{X}_s$ and $X_s$ are generically isomorphic.
\item[(2)] Simultaneous normalizations need not be unique (cf. \cite{Ko11}, Example 1.8). They are unique if $S$ is normal, since then $n$ is the normalization of $X$ by Proposition \ref{pro1.2} (5).
\end{enumerate}
 \end{Remark}

\begin{Lemma}\label{lemma5.5}
If $f:X\to S$ admits a simultaneous normalization, then $f^\red: X^\red\to S^\red$ admits also a simultaneous normalization. The converse holds if the non--empty fibres $X_s$ of $f$ are generically reduced.
\end{Lemma}

\noindent\textbf{Proof}:
\rm{Let $n:\widetilde{X}\to X$ be a simultaneous normalization of $f$ and let $\widetilde{f}^\red:\widetilde{X}'\to S^\red$ be the map obtained from $\widetilde{f}:\widetilde{X}\to S$ by base change $S^\red\hookrightarrow S$.
Then $\widetilde{f}^\red$ is flat with normal, hence reduced fibres and therefore $\widetilde{X}'$ is reduced by Proposition \ref{Prop6.1} (1). Since $\widetilde{X}'$ and $\widetilde{X}$ are topologically the same spaces it follows that $\widetilde{X}'$ is the reduction $\widetilde{X}^\red$ of $\widetilde{X}$. Then $n$ induces a map 
\[
n^\red:\widetilde{X}^\red\to X^\red
\]
which is easily seen to be a simultaneous normalization of $f^\red$.} For the converse we note that $X_s^\red\hookrightarrow X_s$ is bimeromorphic if $X_s$ is generically reduced.\hfill\qed
\bigskip

Although we do not require $f$ to be flat in Definition \ref{defEquinormal}, the following proposition shows that $f$ has to be generically flat in order  to admit a simultaneous normalization and, at non-flat points, $f$ must be ''almost flat'':

\begin{Proposition} \label{pro1.2}
Assume that $f: X\mtn S$ admits a simultaneous normalization $n: \nga{X} \mtn X$.
\begin{enumerate}
\item[(1)] $n_s: \nga{X}_s \mtn X_s$ is the normalization of $X_s$ and $X_s$ is generically normal, hence generically reduced for each $s \in f(X)$.
\item[(2)] $f$ is open and for each  $x\in X$ the dimension formula holds,
$$\dim (X,x) = \dim (X_{f(x)},x) + \dim (S,f(x)). $$
\item[(3)] If $X_{f(x)}$ is reduced at $x$ then $f$ is flat at $x$. In particular, $f$ is generically flat on $X$ and $(X,x)$ is reduced if $(X_{f(x)},x)$ and $(S, f(x))$ are reduced..
\item[(4)] $N:=\NNor(f)$ is nowhere dense in $X$ and  $n^{-1}(N)$ is nowhere dense in $\nga{X}$.  $n: \nga{X}\mtn X$ is surjective and $n: \nga{X}\tru n^{-1}(N) \mtn X\tru N$ is an isomorphism. Moreover,
    $$\NNor(f)\cap X_s = \NNor(X_s) \mbox{ for } s \in f(X).$$
\item[(5)] $\nga{X}$ is normal if and only if $n$ is the normalization of $X$, and this is equivalent to  $f(X) \subseteq \mbox{Nor}(S)$.
In particular, if $S$ is normal then $n$ is the normalization of $X$ and if $f(X)\nsubseteq \mbox{Nor}(S)$ then $f$ is not equinormalizable.
\item[(6)] If $S'\mtn S$ is any morphism, then the pullback $X_{S'}\mtn S'$ of $f$ admits a simultaneous normalization. In particular, if $\gt{S} \mtn S$ is the normalization of $S$ then $X_{\gt{S}} \mtn \gt{S}$ is equinormalizable.
\item[(7)] $\nga{X}$ is weakly normal (resp. reduced) if  and only if $f(X)\subseteq \mbox{Wor}(S)$ (resp. $f(X)\subseteq \mbox{Red}(S)$).
\item[(8)] $\nga{X}$ and $X$ are generically reduced if $S$ is generically reduced.
\end{enumerate}
\end{Proposition}

\noindent \textbf{\large Proof:} \rm{
(1)  Since $\nga{X}_s$ is normal, $n_s$ factors through $(X_s)^{red}$. Hence $n_s: \nga{X}_s \mtn (X_s)^{red}$ is finite and bimeromorphic and therefore the normalization map (cf. \cite{GR84}, Ch. 8, \S 3.3 and \S 4.2). Moreover, since $\widetilde{X}_s$ and $X_s$ are generically isomorphic by Remark \ref{rem1.5}, $X_s$ is generically normal.
\medskip

(2) Since $\nga{f}$ is flat, it is open and therefore $f$ is open. 
To see the dimension formula, note that $\dim (\nga{X},\tilde{x}) = \dim (\nga{X}_s,\tilde{x}) + \dim (S,s)$ for $\tilde{x}\in n^{-1}(x)$ and $s=\tilde{f}(\tilde{x})=f(x)$, since $\tilde{f}$ is flat. As $\nga{X}_s$ is the normalization of $X_s$ we can choose $\tilde{x}\in n^{-1}(x)$ such that
 $\dim (\nga{X}_s,\tilde{x}) = \dim(X_s,x)$.  Then
$$ \dim (X,x) \geq \dim (\nga{X},\tilde{x}) = \dim(X_s,x)+\dim(S,s). $$
Since the other inequality is always true, the dimension formula follows. 
\medskip

(3) If $(X_{s},x)$ is reduced for $s=f(x)$  then $\ko_{X_{s},x}\mtn ({n_s}_*\ko_{\nga{X}_s})_x$ is injective. By \cite {GLS07}, Prop. B.5.3, it follows that $({n}_*\ko_{\nga{X}}/\ko_{X})_x$ is $\ko_{S,s}$-flat and since $(n_*\ko_{\nga{X}})_x$ is $\ko_{S,s}$-flat, we have  that $\ko_{X,x}$ is $\ko_{S,s}$-flat by \cite{GLS07}, Cor. B.5.2. Since the non--empty fibres are generically reduced by (1), $f$ is generically flat. For the rest use Proposition \ref{Prop6.1} (1).
\medskip

(4) $n: \nga{X}\mtn X$ is surjective since $n_s: \nga{X}_s \mtn X_s$ is surjective  for $s\in f(X)$ by (1). 
For $y \in n^{-1}(\mbox{Nor}(X_s)), n: (\nga{X},y)\mtn (X,n(y))$ induces an isomorphism of the fibres by (1).
Since $\tilde{f}$ is flat, $n$ is itself an isomorphism at $y$ (cf. \cite{GLS07}, Lemma I.1.86). Hence $f$ is normal at $n(y)$ and we conclude $ \mbox{NNor}(f)\cap X_s\subseteq \mbox{NNor}(X_s)$. Since the other inclusion is always true, we get $N\cap X_s = \NNor(X_s)$. Moreover, $n$ is bijective over $X\tru N$ and hence $n: \nga{X}\tru n^{-1}(N) \mtn X\tru N$ is an isomorphism. 
Since $\nga{X}_s \mtn X_s$ is bimeromorphic, $\mbox{Nor}(X_s)$ is open and dense in $X_s$ and hence $N\cap X_s=\NNor(X_s)$ has everywhere smaller dimension than $X_s$.  Hence for $x \in N$ we get
\begin{align*}
 \dim (N,x) &\leq \dim(N\cap X_s,x) + \dim (S,s)\\
  & < \dim (X_s,x) + \dim (S,s) = \dim (X,x).
\end{align*}

Therefore $N$ is nowhere dense in $X$ by Ritt's lemma and since $n$ is finite, a dimension argument shows that $n^{-1}(N)$ is nowhere dense in $\nga{X}$. 
\medskip

(5) Since $\nga{f}=f\circ n$ is normal we get that $\nga{X}$ is normal if and only if $\nga{f}(\nga{X}) = f(X) \subseteq \mbox{Nor}(S)$ by Proposition \ref{Prop6.1} If $\nga{X}$ is normal, $n$ is the normalization by (4). If $f(X)\nsubseteq \mbox{Nor}(S)$ then $\widetilde{X}$ is not normal and hence the normalization $\nu: \gt{X} \mtn X$ cannot be a simultaneous normalization of $f$. 
\medskip

(6) The preservation of simultaneous normalization by base change follows from that of flatness. The rest is a consequence of (5).
\medskip

(7) This is an immediate consequence of Proposition \ref{Prop6.1} applied to $\nga{f}$. 
\medskip

(8) Since $\nga{f}$ is reduced, the restriction induces a map $\NRed(\nga{X}) \mtn \NRed(S)$ by Proposition \ref{Prop6.1} This implies for $s=\nga{f}(\nga{x})$,
$$\dim (\NRed(\nga{X}), \nga{x}) \leq \dim (\NRed(S),s) + \dim (\nga{X}_s,\nga{x}) < \dim (\nga{X}, \nga{x}) $$
since $\dim (\NRed(S),s) < \dim (S,s)$ and since $\nga{f}$ is flat. Hence $\NRed(\nga{X})$ is nowhere dense in $\nga{X}$. Since $\nga{X}$ and $X$ are generically isomorphic by (4) we get also that $\NRed(X)$ is nowhere dense in $X$.\hfill\qed

\begin{Remark} \rm \label{remark1.1}
\begin{enumerate}
\item [(1)] Statement (5) was proved in \cite{CL06}, Theorem 2.3 for $f$ a reduced map.
\item [(2)] Simultaneous normalizations are only interesting for morphisms with positive fibre dimension. Namely, if  a finite morphism $f : X \mtn S$ admits a simultaneous normalization, then $f: (X,x)\mtn (S,f(x))$ is an isomorphism for each $x \in X$. \\
In fact, since $f$ is finite, $x$ is an isolated point of $X_s, s = f(x).$ Since $n_s: \nga{X}_s \mtn X_s$ is bimeromorphic, $(\nga{X}_s, n^{-1}(x)) \cong (X_s,x)$ and hence $x$ is a smooth  point of $X_s$.  Therefore $f$ is flat at $x$ (Proposition \ref{pro1.2}(3))  with smooth $0$-dimensional fibre and hence an isomorphism at $x$.
\end{enumerate}
\end{Remark}

Note that $(X, x)$ does not need to be pure dimensional if $f: (X,x)\to (S,s)$ admits a simultaneous normalization, even if $(S,s)$ is normal (cf. Theorem \ref{theorem1}). By the following lemma it is, however, not allowed that the nearby fibres of $f$ have isolated points if the special fibre has positive dimension. Moreover, there is  an intimate  relation between the irreducible components of the total space and of the fibres of $f$.

\begin{Lemma} \label{lemma1.1} Let $f: (X,x) \mtn (S,s)$ have a simultaneous normalization and assume that $(S,s)$ is locally irreducible (e.g. if $(S,s)$ is normal).
\begin{enumerate}
\item[(1)] The germs $(X,x)$ and $(X_s,x)$ have the same number of irreducible components.
\item[(2)] For each irreducible component $(X_s',x)$ of $(X_s,x)$ there exists a unique irreducible component $(X',x)$ of $(X,x)$ containing $(X_s',x)$. The corresponding components satisfy
    $$ \dim (X',x) = \dim (X_s',x) + \dim (S,s).$$
    In particular, if $\dim (X_s,x) > 0 $ then each irreducible component of  $(X,x)$ has dimension $> \dim (S,s)$.
\end{enumerate}
\end{Lemma}

\noindent\textbf{Proof}: \rm{Let $\gt{S} \mtn S$ be the normalization of $S$ and $X_{\gt{S}} \mtn \gt{S}$ the pullback of $f$. Since $(S,s')$ is irreducible for each $s'$ in a neighbourhood of $s$, $(S,s)$ is homeomorphic to $(\gt{S},\bar{s})$ for a unique $\bar{s}\in \gt{S}$. Hence $(X,x)$ is homeomorphic to $(X_{\gt{S}},\bar{x})$ for the unique preimage $\bar{x}$ of $x$ and we get a bijection between the  irreducible components of $(X,x)$ and $(X_{\gt{S}},\bar{x})$ and between those of $(X_s,x)$ and $\big((X_{\gt{S}})_{\bar{s}}, \bar{x}\big)$. By Proposition \ref{pro1.2} (6) we may therefore assume that $(S,s)$ is normal.

Now let  $n : \nga{X} \mtn X$ be a simultaneous normalization of $f$. Then
$n_s: (\nga{X}_s,n^{-1}(x)) \mtn (X_s,x)$ is the normalization of the fibre $(X_s, x)$ of $f$ and $n: (\nga{X}, n^{-1}(x)) \mtn (X,x)$ the normalization of the total space (Proposition \ref{pro1.2}). (1) follows since in both cases the irreducible components are in 1-1 correspondance to the points of $n^{-1}(x)$. Moreover, for each $\nga{x}\in n^{-1}(x), $ $(\nga{X}_s, \nga{x})$ resp. $(\nga{X}, \nga{x})$ is the normalization of the corresponding component $(X_s',x)$  resp. $(X',x)$, proving the first part of (2). For the unique $\tilde{x}\in \nga{X}$ belonging to $(X', x)$ we get
\begin{align}
\dim (X',x) & = \dim (\nga{X},\nga{x}) = \dim (\nga{X}_s,\nga{x}) + \dim (S,s)\nonumber \\
& = \dim (X_s',x) + \dim (S,s).\nonumber
\end{align}

Finally, if $\dim (X_s,x) >0$, then $(X_s,x)$ is not a point and all irreducible components of $(X_s,x)$ have dimension $>0$ and the last statement follows.}\hfill\qed
\bigskip

Let $f:X\to S$ be a morphism of complex spaces such that the fibres $X_s, s\in S$, have only finitely many non--normal points. We call $f$ (locally) \textbf {delta--constant} if the function $s\mapsto \delta(X_s)$ is (locally) constant on $S$. A morphism of germs is $\delta$--constant, if it has a $\delta$--constant representative. 
\medskip

In \cite{CL06} Chiang--Hsieh and Lipman reconsider Teissier's $\delta$--constant criterion and, in addition, consider projective morphisms $f:X\to S$ with fibres of arbitrary dimension, replacing the $\delta$--invariant by the Hilbert polynomial. They prove:

\begin{Theorem} \label{Theo10} {\bf (Chiang--Hsieh, Lipman)}\\
Let $f:X\hookrightarrow S\times \P^r\to S$ be a flat, projective morphism of complex spaces with $S$ irreducible and normal, $X$ locally equidimensional and all fibres $X_s$ reduced. For $s\in S$ let $\nu_s:\overline{X}_s\to X_s$ be the normalization of $X_s$ and with  $\overline{\ko}_{X_s}:=(\nu_{s})_\ast \ko_{\overline{X_s}}$ let
\[
\kh(\overline{\ko}_{X_s})(n)=\chi(\overline{\ko}_{X_s} \otimes \ko_{\P_r}(n))
\]
be the Hilbert polynomial of the coherent sheaf $\overline{\ko}_{X_s}$ on $X_s\subset \P^r$. Then $f$ is equinormalizable iff the Hilbert polynomial $\kh(\overline{\ko}_{X_s})$ is constant on $S$.
\end{Theorem}

\begin{Remark}
\rm{\begin{enumerate}
\item [(1)] Locally equidimensional means that for each $x\in X$ the irreducible components of the germ $(X,x)$ have the same dimension. Since $f$ is flat and $S$ locally irreducible, this implies that $X_s$ is locally equidimensional. Moreover, it follows that $X$ is reduced.
\item [(2)] Theorem \ref{Theo10} is in \cite{CL06} deduced from an analogous theorem for algebraic schemes over a perfect field.
\item [(3)] Since $f$ is flat, the Hilbert polynomial $\kh(\ko_{X_s})$ is constant on $S$. It follows from the exact sequence
\[
0\to \ko_{X_s} \to\overline{\ko}_{X_s}\to \overline{\ko}_{X_s}/\ko_{X_s}\to 0
\]
that $\kh(\overline{\ko}_{X_s})$ is constant iff $\kh(\overline{\ko}_{X_s}/\ko_{X_s})$ is constant. If the fibres are reduced curves, then $\kh(\overline{\ko}_{X_s}/\ko_{X_s})=\delta(X_s)$ and the theorem follows from Theorem 4.4 since the projective curves $X_s$ have only finitely many singularities.
\end{enumerate}}
\end{Remark}

The assumptions in Theorem \ref{Theo10} can be weakened as has been shown by Koll$\acute{\text{a}}$r in \cite{Ko11}.

\begin{Theorem}\label{Theo-ProjMorph} {\bf (Kollar)}\\
Let $f:X\to S$ be a projective morphism of algebraic schemes over a perfect field with $S$ semi--normal, such that all fibres $X_s$ are locally equidimensional of the same dimension and generically reduced. Then $f$ admits a simultaneous normalization iff the Hilbert polynomial $\kh(\overline{\ko}_{X_s})$ is constant on $S$.
\end{Theorem}

For the notion of semi--normal we refer to the reference given in \cite{Ko11}. We have {\em normal} $\Rightarrow$ {\em weakly normal} $\Rightarrow$ {\em semi--normal} $\Rightarrow$ {\em reduced}. In characteristic $0$ weakly normal coincides with semi--normal. A survey on weak-- and semi--normality can be found in \cite{Vi11}.

\begin{Remark}
\rm{\begin{enumerate}
\item [(1)] In Theorem \ref{Theo-ProjMorph} it is not assumed that $f$ is flat, but it follows from the assumptions that $f$ is generically flat on $X$ and that $X$ is generically reduced and equidimensional if $S$ is connected, by Propositon \ref{pro1.2} (3).
\item [(2)] Passing from reduced fibres to generically reduced fibres is a non--trivial generalization, at least from the point of view in \cite{CL06}, where the authors make essential use of reducedness of the fibres.  It follows from Proposition \ref{pro1.2} (1) that generically reduced is a necessary condition for simultaneous normalization.
\end{enumerate}}
\end{Remark}

The approach of  Koll$\acute{\text{a}}$r is different from that of Chiang--Hsieh and Lipman. For a morphism $f:X\to S$ of algebraic schemes over a field, Koll$\acute{\text{a}}$r considers the \textbf{functor of simultaneous normalizations} of $f$,  $\mathbf \SimNor(f)$, which associates to a morphism $T\to S$ the set of simultaneous normalizations of the pull back $X_T\to T$.  Koll$\acute{\text{a}}$r proves in \cite{Ko11}, Theorem 1.

\begin{Theorem} \label{Theo6.13} {\bf (Kollar)} \\
Let $f:X\to S$ be a proper morphism such that all fibres $X_s$ are generically geometrically reduced. Then there exists a morphism of schemes $\pi: S^n\to S$ which represent the functor $\SimNor(f)$.

In particular, for any morphism $\varphi: T\to S$, $X_T\to T$ has a simultaneous normalization iff $\varphi$ factors through $\pi:S^n\to S$.
\end{Theorem}

The construction of $S^n$ and hence of a fine moduli space for simultaneous normalizations is a new idea of  Koll$\acute{\text{a}}$r. It works only for proper and hence for projective morphisms. A similar construction in the local or affine setting is not known.\\
\smallskip

\section{Families of isolated non--normal singularities}

In this section we present our results about $1$--parametric families $f:X\to T$, $T$ a 1--dimensional complex manifold, of isolated non--normal singularities of arbitrary dimension. We use the notations of Section 5. 

In particular, $\nu:\overline{X}\to X$ is the normalization,  $\omega: \widehat{X}\to X$ the weak normalization and
we have $f^\red=f|X^\red, \ \bar{f}=f\circ \nu,\ \widehat{f}=f\circ \omega,\  f^{>1}=f|X^{>1}$, with fibres
\[
 X^\red_t=(f^\red)^{-1}(t),\  \overline{X}_t=\bar{f}^{-1}(t),\ \widehat{X}_t=\widehat{f}^{-1}(t),\ X^{>1}_t= (f^{>1})^{-1}(t). 
\]

Note that $X^\red_t,\ \overline{X}_t$ and $ \widehat{X}_t$ need not be reduced, normal and weakly normal respectively. The corresponding reduction, normalization and weak normalization are denoted correspondingly by $(X_t)^\red, (\overline{X_t})$ and $(\widehat{X_t})$. 
$\overline{X}^{>1}_t$ denotes the fibre of $ \bar{f}^{>1}  = f^{>1}\circ \nu^{>1}  : \overline{X}^{>1} \to T$. We use similar notations for germs.

\bigskip

For topological considerations the flatness assumption is often too strong. For $1$--parametric families we may consider the weaker notion of ''activeness'' which we introduce now.

\begin{Definition}
Let $f:X\to T$ be morphism of complex spaces with $T$ a $1$--dimension complex manifold. We call $f$ \textbf{active} at $x$ if $f^\red:X^\red\to T$ is flat at $x$, and $f$ is called active if it is active at every $x\in X$.
\end{Definition}

\begin{Remark}\label{rem1} \rm
\begin{itemize}
\item [(1)] The notion of acticity was introduced by Grauert and Remmert in \cite{GR84}, \S2.1, who call an element $f\in \ko_{X,x}$ active if the image $f^\red\in \ko_{X^\red,x}$ is a non--zero divisor (n.z.d.) of $\ko_{X^\red,x}$. Since the morphism of germs $f^\red:(X^\red,x)\to (\C,0)$ is flat, iff $f^\red$ is a non--zero divisor of $\ko_{X^\red,x}$, both notions coincide for $f\in \mathfrak{m}_{X,x}$, the maximal ideal of $\ko_{X,x}$. 

\item [(2)] The \textbf{active set of $f$}
\[
\Act(f)=\{x\in X|f\text{ is active at }x\}
\]
coincides with Flat$(f^\red)$ and hence is analytically open in $X$ by the theorem of Frisch. A simple argument is given in  \cite{GR84}. Since $X^\red$ is reduced and $f^\red$ is flat, it follows from \cite{GLS07}, Cor. II, 1.116, that $\Reg(f^{\text{\red}})$ is dense in $X$ for $f$ active. 

\item [(3)]  Since the set of zerodivisors of $\ko_{X,x}$ (resp. $\ko_{X^{red},x}$)  is the union of all associated (resp. minimal)  prime ideals  of $\ko_{X,x}$, $f$ is flat (resp. active) iff $f$ is not contained in any associated (resp. minimal) prime ideal of $\ko_{X,x}$. 
In particular, we get for a morphism $f: (X,x) \to (\C,0)$:\\
\[
f \text{ is active } \Leftrightarrow \dim \ko_{X_i,x}/ f = \dim \ko_{X_i,x} - 1  
\]
for each irreducible component $(X_i,x)$ of $(X,x)$.

\item [(4)] Activeness is a purely topological condition. We have for $f: X\to T$
\[
f \text{ is active } \Leftrightarrow f \text{ is open}.
\]
The openness follows since flat maps are open (by \cite{Fi76}, 3.10, Lemma 1) and since $f$ and $f^\red$ are topologically the same maps.

Conversely, let $f$ be open and $x\in X$. If $(X,x)$ is irreducible, then $f(X) \neq 0$ and we are done. Otherwise let $(X', x)$ be an irreducible component of $(X,x)$ and $(X'',x)$ the union of the remaining irreducible components of $(X,x)$. Then there exist $y\in X'\smallsetminus X''$ arbitrary close to $x$ and we can choose an open neigbourhood $U$ in $X$ with $U\cap X'' = \emptyset$. Since $f(U)$ is open by assumption, $f|X'$ is not constant. This shows that $f$ is not contained in any minimal prime of $\ko_{X,x}$ and hence $f$ is active.

\end{itemize}
\end{Remark}

The next lemma and Lemma \ref{Lem7.4} show that activeness is a natural assumption.

\begin{Lemma} \label{lemma7.3} Let $f: (X,x) \sr (\C,0)$  be a morphism of complex germs.
\begin{enumerate}
\item [(1)] If $f$ is flat, then $f$ is active.
\item [(2)] $f$ is active if and only if $\bar{f}:(\overline{X}, \bar{x})\to (\C,0)$ is flat.
\end{enumerate}
\end{Lemma}

\noindent\textbf{Proof:} \rm{ (1) follows from Remark \ref{rem1} (3).

(2) If $\bar{f}$ is flat at $\bar x$, then $\bar{f}$ is a n.z.d. of $\ko_{\gt{X},\nu^{-1}(x)}$. Since $\ko_{X^{red},x}$ is a subring of $(\nu_*\ko_{\gt{X}})_x$, $f^{red}$ is a n.z.d. of $\ko_{X^{red},x}$, i.e. $f$ is active at $x$.

For the converse consider the conductor
\[\kc:=\Ann_{\ko_{X^{red}}}\big( n_*\ko_{\gt{X}}/\ko_{X^{red}}\big),
\]
which is  coherent  because the sheaf $n_*\ko_{\gt{X}}/\ko_{X^{red}}$ is coherent. Denote by $\NNor(X^{red})$   the non-normal locus of
$X^{red}$ which is nowhere dense in $X^{red}$. Then $\NNor(X^{red}) $ is the vanishing locus of $\kc$.  It follows from the prime avoidance theorem  that there exists
an element $h \in \kc_x$ such that $h$ is not contained in any associated prime of $\ko_{X^{red},x}$. Thus $h$ is a n.z.d. of
$\ko_{X^{red},x}$ and $h(n_*\ko_{\gt{X}})_x\subseteq\ko_{X^{red},x}$. Therefore we have $(n_*\ko_{\gt{X}})_x\subseteq
h^{-1}\ko_{X^{red},x}$.
Since $f$ is active,  $f^{red}$ is a n.z.d. of $\ko_{X^{red},x}\cong h^{-1}\ko_{X^{red},x}$, hence of
$(n_*\ko_{\gt{X}})_x$. This implies that $\bar{f}$ is a n.z.d. of $\ko_{\gt{X},\nu^{-1}(x)}$.\hfill\qed
\bigskip

\noindent More generally, we have the following result:

\begin{Lemma} \label{lemma7.4} 
\begin{enumerate}
\item [(1)] If $f$ is flat and $(X',x)$ is a subgerm of $(X,x)$ defined by the intersection of some primary or prime ideals of $\ko_{X,x}$,  then $f\mid_{(X',x)}$ is also flat.
\item [(2)] If $f$ is active and $\pi:(\widetilde{X},\widetilde{x})\to (X,x)$  a \emph{partial normalization} of $(X,x)$, then $\widetilde{f}=f\circ \pi: \widetilde{X}\to \C$ is flat at every point of \ $\widetilde{x}=\pi^{-1}(x)$.
\end{enumerate}
\end{Lemma}

\noindent \textbf{\large Proof:} 
\rm{ (1) follows again from Remark \ref{rem1} (3).

(2) Since $\bar{f}$ is a n.z.d. of  $(\nu_*\ko_{\gt{X}})_x$ by Lemma \ref{lemma7.3}, $\widetilde{f}$ is a n.z.d. of $(\pi_\ast\ko_{\widetilde{X}})_x\subset(\nu_\ast\ko_{\overline{X}})_x$, since $\pi$ factors through the reduction of $X$.
\hfill\qed
\bigskip

\noindent Note that the statement of Lemma \ref{lemma7.3} (1) is false for a morphism $f: (X,x)\mtn (S,s)$ with $(S,s)$ not smooth or of dimension bigger than one. See \cite{Fi76}, 3.13, Example,  for an example with $(S,s)$ reduced of dimension $1$ and \cite{CN76} with $f: (X,x) \mtn (\C^2,0)$, $f$  flat but $f^{red}$ not flat.\hfill\qed
\bigskip

The following lemma shows that activeness is a necessary condition for equinormalizability (while flatness is not: we may add an extra fat point to $(X,x)$ of arbitrary length, destroying flatness but neither activeness nor simultaneous normalization, cf. Lemma \ref {lemma5.5}).
\medskip

\begin{Lemma}\label{lemma0} If $f: (X,x) \mtn (\C,0)$  admits a simultaneous normalization, then $f$ is active.
\end{Lemma}

\noindent \textbf{\large Proof:}  This follows from Lemma \ref{lemma7.3} (2) .\hfill\qed
\medskip

\begin{Lemma} \label{lemma1} Let $f: (X,x) \sr (\C,0)$ be active and $(X_0,x)$  reduced. Then $(X,x)$ is reduced and $f$ is flat. 

Moreover, $(X,x)$ is either 1-dimensional and smooth or $\depth(X,x)\geq 2$ and $\dim(X_i,x)\geq 2$ for each irreducible component $(X_i,x)$ of $(X,x)$.

\end{Lemma}

\noindent \textbf{\large Proof:} Consider the diagram
\[
\xymatrix{0 \ar[r] & \mathcal{N}\text{il}(\ko_X)\ar[d]^f \ar[r] & \ko_X \ar[r]\ar[d]^f & \ko_{X^\red} \ar[r]\ar[d]^{f^\red} & 0\\
0 \ar[r] & \mathcal{N}\text{il}(\ko_X) \ar[r] & \ko_X \ar[r] & \ko_{X^\red} \ar[r] & 0.
}
\]
Since $f$ is active, $f^\red$ is injective and since $\ko_{X_0}$ is reduced, the sanke lemma implies that $\mathcal{N}\text{il}(\ko_X)/f \mathcal{N}\text{il} (\ko_X)=0$. By Nakayama's lemma $\mathcal{N}\text{il}(\ko_X)=0$, $X$ is reduced and $f$ is flat. If $\dim (X,x) = 1$ then $(X_0,x)$ is a reduced point and $f$ is locally an isomorphism. If $\dim (X,x) \geq 2$ then $\depth(X_0,x) \geq 1$ since it is reduced of positive dimension, hence $\depth(X,x) \geq 2$. \hfill\qed\\

\begin{Proposition} \label{prop1} Let $f:(X,x)\to (\C, 0)$ be active, $\nu:(\overline{X}, \nu^{-1}(x))\to (X,x)$ the normalization and $\omega:  (\widehat{X}, \widehat{x})\to (X,x)$ the weak normalization. For good representatives $f: X\to T, \bar{f}: \overline{X}\to T$ and $\widehat{f}:\widehat{X}\to T$ the following holds.
\begin{enumerate}
\item [(1)] Let $X_0\smallsetminus \{x\}$ be reduced.
\begin{enumerate}
\item [(i)] $f$ is flat on $X\smallsetminus \{x\}$, $\bar{f}$ is flat on $\overline{X}$ and $\widehat{f}$ is flat on $\widehat{X}$.
\item [(ii)] $\NRed(f)$ is finite over $T$ and hence all fibres $X_t, t\in T$, are reduced outside finitely many points.
\item [(iii)] $X$ is generically reduced and reduced at  $y\in X_0\smallsetminus\{x\}$. We have $X\smallsetminus \{x\}$ reduced $\Leftrightarrow f$ reduced on $X\smallsetminus \{x\}\Leftrightarrow X_t$ reduced for $t\in T\smallsetminus\{0\}$.
\item [(iv)] If $f$ is flat, then $X$ is reduced iff $X_t$ is reduced for $t\in T\smallsetminus\{0\}$.
\end{enumerate}

\item [(2)] Let $X_0\smallsetminus\{x\}$ be normal.
\begin{enumerate}
\item [(i)] $\NNor(f)$ is finite over $T$ and all fibres $X_t,\  t\in T$, are normal outside finitely many points.
\item [(ii)] $\bar{f}$ is normal outside $\nu^{-1}(x)$ and $\overline{X}_t\to X_t$ is the normalization for $t\in T\smallsetminus\{0\}$.
\item [(iii)] $X$ is generically normal and we have $X\smallsetminus\{x\}$ normal $\Leftrightarrow f$ normal on $X\smallsetminus\{x\}\Leftrightarrow X_t$ normal for $t\in T\smallsetminus\{0\}$. 
\item [(iv)] $\overline{X}_t$ is reduced for $t\in T\smallsetminus\{0\}$. 
\end{enumerate}

\item [(3)] Let $X_0\smallsetminus\{x\}$ be weakly normal.
\begin{enumerate}
\item [(i)] $\NWor(f)$ is finite over $T$ and all fibres $X_t, \ t\in T$, are weakly normal outside finitely many points. 
\item [(ii)] $\widehat{f}$ is weakly normal outside $\widehat{x}=\pi^{-1}(x)$ and $\widehat{X}_t\to X_t$ is the weak normalization for $t\in T\smallsetminus\{0\}$. 
\item [(iii)] $X$ is generically weakly normal and we have $X\smallsetminus \{x\}$ weakly normal $\Leftrightarrow f$ weakly normal on $X\smallsetminus \{x\}\Leftrightarrow X_t$ weakly normal for $t\in T\smallsetminus \{0\}$ .
\item [(iv)] $\widehat{X}_t$ is reduced for $t\in T\smallsetminus\{0\}$. 
\end{enumerate}
\end{enumerate}
\end{Proposition}

\noindent\textbf{\large Proof:} 
It follows from Lemma \ref{lemma1} that in any of the three cases $f$ is flat at points of $ X_0\smallsetminus\{x\}$.
\medskip

(1)(i) Since $f$ is flat at points $y\in X_0\smallsetminus\{x\}$, the analytic set $\NFlat(f)$ is finite over $T$. Then $f(\NFlat(f))$ is a closed analytic set in $T$, which is nowhere dense by Theorem \ref{thr-neglectible} (1) and therefore equal to $\{0\}$. It follows that $f$ is flat on $X\smallsetminus\{x\}$.

The flatness of $\bar{f}$ and $\widehat{f}$ follows now from Lemma \ref{lemma7.3} and Lemma \ref{lemma7.4}.
\medskip

(ii) Since $f$ is flat on $X\smallsetminus \{x\}$ by (i), we have
\[
\NRed(f)\cap (X_0\smallsetminus\{x\})=\NRed(X_0\smallsetminus\{x\})=\emptyset
\]
and hence $\NRed (f)$ is finite over $T$.
\medskip

(iii) $X$ is generically reduced by (ii) and Proposition \ref{Prop6.1} (1) and $X$ is reduced at $y\in X_0\smallsetminus\{x\}$.\\
If $X\smallsetminus\{x\}$ is reduced then we can apply Theorem \ref{thr-neglectible} (2) (iv) to get that $f(\NRed(f))=\{0\}$  and that $f$ is reduced on $X\smallsetminus\{x\}$. Hence $X_t, t\in T\smallsetminus\{0\}$, is reduced. If $X_t , t\in T\smallsetminus\{0\}$, is reduced, then $f$ is $X$ is reduced at $y\in X_t$ by (1), and hence $X \smallsetminus \{x\}$ reduced since it is reduced at $y\in X_0\smallsetminus\{x\}$.
\medskip

(iv) Let $f$ be flat and $X_t$ reduced for $t\neq 0$. Then $X\smallsetminus \{x\}$ is reduced by (iii). Hence $\Nil(\ko_{X,x})$ has support $\{x\}$ and is therefore killed by some power of $f$. Since $f$ is flat, every power of $f$ is a n.z.d of $\ko_{X,x}$, implying $\Nil(\ko_{X,x})=0$.
\bigskip

(2) (i) follows as in (1)(ii).
\medskip

(ii) $\NNor (\bar{f})$ is finite over $T$. Since $\bar{f}$ is flat by (1) (i), it follows from Theorem \ref{thr-neglectible} (2) (ii) applied to $\bar{f}$ that $\bar{f}(\NNor(\bar{f})=\{0\}$ and hence that $\bar{f}$ is normal outside $\nu^{-1}(x)$. It follows that $\overline{X}_t$ is normal for $t\neq 0$ and that $\overline{X}_t\to X_t$ is the normalization.
\medskip

(iii) Applying Theorem \ref{thr-neglectible} (2) (ii) to $\NNor(f)$, this follows as in (1) (iii).
\smallskip

(iv) Since $\bar f$ is flat and $\bar X$ reduced, this follows from (1) (iv).
\bigskip

(3) The proof is the same as for (2), using Theorem \ref{thr-neglectible} (2) (iii).\hfill\qed
\bigskip

\begin{Remark} \label{rem2.9}
\rm \begin{enumerate}
\item [(1)] If $f:(X,x)\to (\C,0)$ is active and $f:X\to T$ a good representative, then we get for the number of 1-dimensional components of $(X,x)$: $r_1(X,x)=r_0(X_t)=\sharp \{\text{ isolated points of } X_t\}$ for $t\in T\smallsetminus\{0\}$.

\item [(2)] If $r_1(X,x)\not = 0$, let $X^1$ be the union of the $1$--dimensional irreducible components of $X$. Since $X^1$ is reduced, the fibre $X^1_t$ of $f^1=f | X^1$ is reduced for $t\neq 0$ (Proposition \ref{prop1}(1)) with $\sharp X^1_t= \sharp\{\text{isolated points of }X_t\}$ (note that $X$ and $X_t$ need not be reduced at points $y\in X_t^1$). $f^1$ is flat and hence (by \cite{GLS07}, Theorem 1.81 (2))
$$ \varepsilon(X_0^1) = \dim_\C\ko_{X^1,x}/\langle f \rangle = \sharp\{\text{isolated points of }X_t, t\neq 0\} .$$
\end{enumerate}
\end{Remark}
\bigskip

\begin{Proposition}\label{lemma5}
Let $f: (X,x)\mtn (\C,0)$ be active with $(X_0,x)$ an INNS. 
\begin{enumerate}
\item [(i)] If $r_1(X,x)=0$, then $(\overline{X}_0, \bar{x})=(\bar{f}^{-1}(0), \bar{x})$ is reduced.

\item [(ii)] If $r_1(X,x)\neq 0$, let $(X',x)$ be an irreducible component of $(X^1,x)$ with  $x'\in\nu^{-1}(x)$ the corresponding point in  $\overline{X}$. Then
then $(\overline{X}_0, x')$ is reduced iff  $(X', x)$ is smooth and $ \varepsilon(X'_0)=1$. 

\item [(iii)] The fibres $\widehat{X}_t$ of $\widehat{f}$ and $X_t$ of $f$ have the same normalization for all $t \in T$. If $\dim(X,x) \geq 2$, then $(\overline{X}_0,\overline{x})$ and $(X_0,x)$ have the same normalization iff $r_1(X,x)=0$.

\item [(iv)] Let $f$ be flat. Then $\depth (X,x) \geq 2$ iff $r_1(X,x)=0$ and $(X_0,x)$ is reduced. 

If $\dim(X,x)\geq 2$, then $(X_0,x)$ reduced implies $r_1(X,x)=0$ and hence $(X_0,x)$ reduced is equivalent to $\depth (X,x) \geq 2$.
\end{enumerate}
\end{Proposition}

\noindent \textbf{\large Proof:}
\rm{(i) Follows from (iv) applied to $\bar f$.
\medskip

(ii) If $(\overline{X}_0, x')$ is a reduced point, then $\bar f : (\overline{X},x') \to (\C,0)$ is an isomorphism and hence $f : (X',x) \to (\C,0)$ is an isomorphism.
\medskip

(iii) Since $\widehat{X}_t\to X_t$ is bimeromorphic, the normalization of $\widehat{X}_t$ is also the normalization of $X_t$. If $\dim (X,x)\geq 2$, then $\dim (X_0, x)\geq 1$ and hence its normalization $(\overline{X_0})$ satisfies $\dim (\overline{X_0}, x')\geq 1$ for each $x'\in \overline{x}$. Since $r_1(X,x)$ is the number of isolated points of the fibre $\overline{X}_0$, $r_1(X,x)=0$}.
\medskip

(iv) If $\depth (X,x) \geq 2$,  all components of $(X,x)$ have dimension $\geq 2$ and $\depth (X_0,x) \geq 1$, implying $(X_0,x)$ reduced since it is an INNS. Conversely, let $\dim(X,x)\geq 2$ and $(X_0,x)$ reduced. Then $\dim(X_0,x)\geq 1$ and $\depth (X_0,x) \geq 1$, which implies $\depth (X,x) \geq 2$ and therefore $r_1(X,x)=0$, proving the claim.
\hfill\qed
\bigskip

For a deeper investigation of the behavior of $\delta$ and other invariants in active families $f:X\to T$, we need to consider partial normalizations of $f$ and the restriction of $f$ to subspaces of $X$. In order to treat these situations simultaneously, we introduce the notion of ''moderation''.

As isolated points play a special role for the definition of $\delta$, so does the $1$--dimensional part of $X$ when we study 1-parametric $\delta$--constant families of isolated non--normal singularities. This $1$--dimensional part can be changed by a moderation which we define now.

\begin{Definition} Let $f:X\to T$ be active with $T$ a pure $1$--dimensional complex manifold and with $\NNor(f)$ finite over $T$. A morphism $\pi:\widetilde{X}\to X$ is called a \textbf{moderation} of $f$ if
\begin{enumerate}
\item [(i)] $\pi$ is finite.
\item [(ii)] There exists a closed analytic subset $N\subset X$, finite over $T$, such that, with $\widetilde{N}=\pi^{-1}(N), \pi:\widetilde{X}\smallsetminus \widetilde{N}\to X\smallsetminus N$ is an isomorphism.
\item [(iii)] $\widetilde{f}=f\circ\pi:\widetilde{X}\to T$ is active.
\end{enumerate}
\end{Definition}

Examples of a moderation are the reduction $X^\red\hookrightarrow X$, the inclusion $X^{>1}\hookrightarrow X$, the normalization  $\overline{X}\to X$ and the weak normalization $\widehat{X}\to X$.\\

\begin{Lemma}\label{lemma6.13}
Let $\pi:\widetilde{X}\to X$ be a moderation of $f:X\to T$. 
\begin{enumerate}
\item [(1)] If $\pi^\sharp:\ko_X\to\pi_\ast\ko_{\widetilde{X}}$ is the induced map of structure sheaves, then $\Ker(\pi^\sharp)$ and $\Coker(\pi^\sharp)$ are finite over $T$.
\item [(2)] $\NNor(\widetilde{f})$ is finite over $T$.
\item [(3)] $\pi:\widetilde{X}^{>1}\to X^{>1}$ is bimeromorphic.
\item [(4)] For $t\in T$ the fibre $\widetilde{X}^{>1}_t$ of $\widetilde{f}|\widetilde{X}^{>1}$ has the same normalization as the fibre $X^{>1}_t$ of $f|X^{>1}$.
\end{enumerate}
\end{Lemma}

\noindent\textbf{Proof:} (1) $\Ker(\pi^{\sharp})$ and $\Coker(\pi^{\sharp})$ are contained in $N$, hence finite over $T$.
\medskip

(2) $\NNor(\widetilde{f})$ is finite over $T$ since it is contained in $\widetilde{N}$.
\medskip

(3) Since $\dim \widetilde{N}=\dim N\leq 1$, $\widetilde{N}$ and $N$ are nowhere dense in $\widetilde{X}^{>1}$ and $X^{>1}$ respectively.
\medskip

(4) Since $N$ and $\widetilde{N}$ are finite over $T$, $N\cap X^{>1}_t$ and $\widetilde{N}\cap \widetilde{X}^{>1}_t$ consist of finitely many points and $\pi_t:\widetilde{X}^{>1}_t\to X^{>1}_t$ is bimeromorphic. If $\nu_t:(\overline{X_t^{>1}})\to \widetilde{X}^{>1}_t$ is the normalization of $\widetilde{X}^{>1}_t$, then $\pi_t\circ\nu_t$ is the normalization of $X_t^{>1}$.\hfill\qed
\bigskip

The following lemma is used when we study the fibres of a moderation.
\medskip

\begin{Lemma}\label{lemma6.11} Let $\pi:(\widetilde{Z},\widetilde{z})\to (Z,z) , \ \widetilde{z}=\pi^{-1}(z),$ be finite such that for small representatives  $Z\smallsetminus\{z\}$ is normal and $\pi:\widetilde{Z}\smallsetminus\widetilde{z}\to Z\smallsetminus\{z\}$, is an isomorphism. Let
\[
\begin{array}{rcl}
\varphi & : & \ko_{Z,z}\to \pi_\ast(\ko_{\widetilde{Z}})_z\text{ and}\\
\varphi^{>0}& : &  \ko_{Z^{>0}, z}\to \pi_\ast(\ko_{\widetilde{Z}^{>0}})_z
\end{array}
\]
be the maps induced by $\pi$ and $\pi^{>0}:\widetilde{Z}^{>0}\to Z^{>0}$.

Then $(Z^{>0}, z)$ and $(\widetilde{Z}^{>0}, \widetilde{z}'), \ \widetilde{z}'=(\pi^{>0})^{-1} (z)$, have the same normalization and $\varphi$ has a finite dimensional kernel and cokernel, satisfying
\[
\begin{array}{lcl}
\dim_\C\Coker (\varphi)-\dim_\C\Ker (\varphi) & = & \delta(Z,z)-\delta(\widetilde {Z},\widetilde{z})\\
& = & \varepsilon(\widetilde{Z}, \widetilde{z})-\varepsilon(Z,z)+\dim_\C\Coker(\varphi^{>0})
\end{array}
\] 
\end{Lemma}

\noindent\textbf{Proof:}  Since $\pi:\widetilde{Z}\smallsetminus\widetilde{z}\to Z\smallsetminus\{z\}$ is an isomorphism, $\varphi$ as well as $\varphi^{>0}$ have finite dimensional kernel and cokernel. If $\dim(Z,z)=0$, then $\widetilde{z}$ consists of finitely many isolated points of $\widetilde{Z}$ and we have $\dim_\C\ko_{Z,z}=\varepsilon(Z,z)=-\delta(Z,z)$ and $\dim_\C(\pi_\ast, \ko_{\widetilde{Z}})_z=\varepsilon(\widetilde{Z}, \widetilde{z})=-\delta(\widetilde{Z}, \widetilde{z})$. The claim follows with $\Coker(\varphi^{>0})=0$.

Let $\dim(Z,z)\geq 1$. Then $(Z^\red, z)=(Z^{>0}, z)$ and $(\widetilde{Z}^\red, \widetilde{z})=(\widetilde{Z}^{>0}, \widetilde{z}')\cup \{\text{finitely many isolated points of }\widetilde{Z}\}$. Let $n:(\overline{Z^{>0}}, \bar{z}) \to (\widetilde{Z}^{>0}, \widetilde{z}')$ be the normalization map. Since $\widetilde{Z}^{>0}\smallsetminus \widetilde{z}'$ is dense in $\widetilde{Z}^{>0}$ and $\pi^{>0}:\widetilde{Z}^{>0}\smallsetminus \widetilde{z}'\to Z^{>0}\smallsetminus\{z\}$ is an isomorphism, $\nu=\pi\circ n:(\overline{Z}^{>0}, \bar{z})\to (Z^{>0}, z)$ is the normalization.

 Now consider the $2$--term complexes
\[
\begin{array}{l}
\ko^\bullet : 0\to\ko_{Z,z}\to (\nu_\ast\ko_{\overline{Z^{>0}}})_z\to 0\\
\widetilde{\ko}^\bullet : 0\to \pi_\ast(\ko_{\widetilde{Z}})_z\to \pi_\ast(n_\ast\ko_{\overline{Z^{>0}}})_z\to 0
\end{array}\raisebox{2ex}{\text{ and }}
\]
and the morphism $\varphi^\bullet: \ko^\bullet\to \widetilde{\ko}^\bullet$, with $\varphi^0=\varphi$ and  the identity in degree $1$.

Let $K^\bullet$ resp. $C^\bullet$ be the $1$--term complexes $\Ker(\varphi^\bullet)$ resp. $\Coker (\varphi^\bullet)$, concentrated in degree $0$. Then we have the exact sequence of complexes
\[
0\to K^\bullet\to \ko^\bullet\to \widetilde{\ko}^\bullet\to C^\bullet\to 0\ .
\]
Taking Euler characteristics, we get $\dim_\C\Coker(\varphi)-\dim_\C\Ker(\varphi)=\chi(\widetilde{\ko}^\bullet)-\chi(\ko^\bullet)=\delta(Z)-\delta(\widetilde{Z})$, showing the first equality. 

Since $\delta(Z)=\delta(Z^{>0})-\varepsilon(Z)$ we get
\[
\delta(Z)-\delta(\widetilde{Z})=\varepsilon(\widetilde{Z})-\varepsilon(Z)+\delta(Z^{>0})-\delta(\widetilde{Z}^{>0}).
\]
From the inclusions $\ko_{Z^{>0}, z}\hookrightarrow(\pi_\ast\ko_{\widetilde{Z}^{>0}})_z\hookrightarrow(\nu_\ast\ko_{\overline{Z^{>0}}})_z$ the second equality follows.\hfill\qed
\medskip

\begin{Proposition}\label{proposition6.14}
Let $f:(X,x)\to(\C, 0)$ be flat with $(X_0, x)$ an INNS and let $\pi:(\widetilde{X}, \widetilde{x})\to (X,x)$ be a moderation of $f$ such that $\widetilde{f}=f\circ\pi$ is flat. Then we have for good representatives $f:X\to T$ and $\widetilde{f}:\widetilde{X}\to T$ with fibres $X_t=f^{-1}(t)$ and $\widetilde{X}_t=\widetilde{f}^{-1}(t), t\in T$:
\begin{enumerate}
\item [(i)] $\delta(X_0)-\delta(X_t)=\delta(\widetilde{X}_0)-\delta(\widetilde{X}_t)$,
\item[(ii)] $\varepsilon(X_0)-\varepsilon(X_t)=\varepsilon(\widetilde{X}_0)-\varepsilon(\widetilde{X}_t)_x+\dim_\C\Coker(\varphi_0)^{>0}-\dim_\C\Coker(\varphi_t)^{>0}$.
\end{enumerate}
\end{Proposition}
Here $(\varphi_t)^{>0}:\ko_{(X_t)^{>0}}\to \pi_\ast \ko_{(\widetilde{X}_t)^{>0}}$ is the map induced by $\pi$, and $(X_t)^{>0}$ is the positive dimensional part of $(X_t)^\red$.\\

\noindent\textbf{Proof:} $\pi: (\widetilde{X}_0, \widetilde{x})\to (X_0, x)$ satisfies the assumptions of Lemma \ref{lemma6.11} and hence we get
\[
\dim_\C\Coker(\varphi_0)-\dim_\C\Ker(\varphi_0)=\delta(X_0)-\delta(\widetilde{X}_0)\ 
\]
for 
\[
\varphi_0:\ko_{X_0, x}\to (\pi_\ast\ko_{\widetilde{X}_0})_x,
\]
the map induced by $\ko_{X,x}\to (\pi_\ast\ko_{\widetilde{X}})_x$.

Now consider the commutative diagram for fixed $t\in T$
\[
\xymatrix{0 \ar[r] & \ko_X\ar[d]\ar[r]^{f-t} & \ko_X\ar[r]\ar[d] & \ko_{X_t}\ar[r]\ar[d]^{\varphi_t}\ar[r] & 0\\
\Theta\ar[r] & \pi_\ast\ko_{\widetilde{X}}\ar[r]^{\widetilde{f}-t} & \pi_\ast\ko_{\widetilde{X}}\ar[r] & \pi_\ast\ko_{\widetilde{X}_t}\ar[r] & 0\ ,
}
\]
with exact rows, since $f$ and $\widetilde{f}$ are flat and $\pi_\ast$ is an exact functor $(\pi$ is finite). We set
\[
\begin{array}{lcl}
\kn: & = & \Ker(\ko_X\to \pi_\ast\ko_{\widetilde{X}})\ ,\\
\km: & = & \Coker(\ko_X\to \pi_\ast\ko_{\widetilde{X}})\ .
\end{array}
\]
Since $\kn$ and $\km$ are finite over $T$, $f_\ast\kn$ and $f_\ast\km$ are $\ko_T$--coherent. Since $\ko_{T, 0}$ is a principal ideal domain we have, for sufficiently small $T$, a decomposition
\[
f_\ast\km=\kf\oplus\kt
\]
with $\kf$ a free $\ko_T$--module and $\kt$ an $\ko_T$--torsion sheaf concentrated on $\{0\}$. $f_\ast\kn$ is $\ko_T$--free (since $f$ is flat and $\kn\subset \ko_X$) of rank
\[
\dim_\C(f_\ast\kn) (0)=\dim_\C(f_\ast\kn)(t),
\]
where $\kg(t):=\kg_t\otimes_{\ko_{T,t}}\!\C$ for an $\ko_T$--module $\kg$. We have also
\[
\dim_\C\kf(0)=\dim_\C \kf(t)\ .
\]
The snake lemma, applied to the diagram above, gives the exact sequences
\[
\begin{array}{c}
0 \to f_\ast\kn\xrightarrow{f-t} f_\ast\kn\to\Ker(\varphi_t)\to f_\ast\km\xrightarrow{\widetilde{f}-t} f_\ast\km\to \Coker (\varphi_t)\to 0,\\[0.5ex]
0\to f_\ast\kn(t)\to \Ker(\varphi_t)\to \Ker(\widetilde{f}-t)\to 0,\\[0.5ex]
0\to \Ker(\widetilde{f}-t)\to \kf_t\oplus\kt_t\xrightarrow{\widetilde{f}-t}\kf_t\oplus\kt_t\to \Coker(\varphi_t)\to 0\ .
\end{array}
\]
We have $\kt_t=0$ for $t\neq 0$ and $\dim_\C\Ker(\widetilde{f}|{\kt_0})=\dim_\C\Coker(\widetilde{f}|{\kt_0})$ since $\kt_0$ is finite dimensional. Therefore
\[
\begin{array}{c}
\dim_\C\Coker(\varphi_0)-\dim_\C\Ker(\varphi_0)=\dim_\C\kf(0)-\dim_\C(f_\ast \kn)(0)\\=\dim_\C\kf(t)-\dim_\C(f_\ast\kn)(t)=\dim_\C\Coker(\varphi_t)-\dim_\C\Ker(\varphi_t).
\end{array}
\]
Applying Lemma \ref{lemma6.11} to $\varphi_t$ we get
\[
\dim_\C\Coker(\varphi_t)-\dim_\C \Ker(\varphi_t)=\delta(X_t)-\delta(\widetilde{X}_t)\ ,
\]
which proves (i).
\medskip

To see (ii) we apply Lemma \ref{lemma6.11} and get for $t\in T$
\[ \varepsilon(\widetilde{X}_t)-\varepsilon(X_t)=\dim_\C\Coker(\varphi_t)-\dim_\C\Ker(\varphi_t)-\dim_\C\Coker(\varphi_t)^{>0}.
\]
Using again
\[
\dim_\C\Coker(\varphi_0)-\dim_\C\Ker(\varphi_0)=\dim_\C\Coker(\varphi_t)-\dim_\C\Ker(\varphi_t)
\]
we get the result. \hfill\qed
\bigskip

Applying Proposition \ref{proposition6.14} to some special moderations we get the first main result.
\medskip

\begin{Theorem}\label{theorem6.15}
	Let $f:(X,x)\to (\C, 0)$ be flat with fibre $(X_0, x)$ an isolated non--normal singularity. For a good representative $f:X \to T$ the following holds.
	\begin{enumerate}
		\item [(1)] If $t\neq 0$ then $\delta(X_0)-\delta(X_t)$ is equal to
		\begin{enumerate}
			\item [(i)] $\delta(X^\red_0)-\delta(X^\red_t)=\delta(X^\red_0)-\delta((X_t)^\red)$,
			\item [(ii)] $\delta(\overline{X}_0)-\delta(\overline{X}_t)=\delta(\overline{X}_0)+r_1(X,x)$,
			\item [(iii)] $\delta(X^{>1}_0)-\delta(X^{>1}_t)=\delta(\overline{X}^{>1}_0)\geq 0$,
			\item[(iv)]
			$\delta(\widehat{X}_0)-\delta(\widehat{X}_t)
			 =  r'(X_0)-r'(X_t)+r_1(X,x)+\dim_\C\Coker(w^\sharp)-\varepsilon(\widehat{X}_0)$.		
		\end{enumerate}
		\item [(2)] $\varepsilon (X_0)-\varepsilon(X_t)=\varepsilon(X^{>1}_0)\geq 0$ for $t\neq 0$.
	\end{enumerate}
\end{Theorem}

Here $w:(\widehat{X_0})\to\widehat{X}_0$ is the weak normalization of $\widehat{X}_0=\widehat{f}^{-1}(0)$ and $w^\sharp : \ko_{\widehat{X}_0}\to w_\ast\ko_{(\widehat{X_0})}$ is the induced map. 
\bigskip

\noindent\textbf{Proof:} 
(1) We apply Proposition \ref{proposition6.14} to $f^\red: X^\red\to T$, to $f^{>1}:X^{>1}\to T$ and to $\widehat{f} : \widehat{X} \to T$, which are flat by Lemma \ref{lemma7.4}.
\medskip

$(i)$ follows then from Proposition \ref{prop1} (1)(ii), saying that $X^\red_t$ is reduced, while $(ii)$ follows from Proposition \ref{prop1} (2)(ii), noting that $\overline{X}_t$ is normal, hence $-\delta(\overline{X}_t)=\sharp\{\text{isolated points of }\overline{X}_t\}=\sharp\{\text{isolated points of } X_t\}=r_1(X,x)$.
\medskip

Applying Proposition \ref{proposition6.14} to $f$ and $f^{>1}$ gives the first equality and applying it again to $f^{>1}$ and $\overline{f^{>1}}$ we deduce $(iii)$ from $(ii)$.
\medskip

To prove $(iv)$ we  apply Proposition \ref{proposition6.14} to $\widehat{f}$  and get $\delta(X_0)-\delta(X_t)=\delta(\widehat{X}_0)-\delta(\widehat{X}_t)$ for $ t\in T$.

Since $\widehat{X}_t$ is weakly normal by Proposition \ref{prop1} (3), we get from Lemma \ref{lemma8.2} for $t\neq 0$ 
\[
\begin{array}{lcl}
\delta(\widehat{X}_t) & = & \delta(\widehat{X}_t)^{>0}-\varepsilon(\widehat{X}_t) \\
& = & r'(\widehat{X}_t)^{>0}-\sharp \text{ isolated points of }\widehat{X}_t \\
& = & r'(X_t)-r_1(X,x).
\end{array}
\]
To determine $\delta(\widehat{X}_0)$, consider the weak normalization $ w:(\widehat{X_0}) \to  \widehat{X}_0$ and the normalization $\mu: (\overline{X_0})\to \widehat{X}_0$ of $\widehat{X}_0$, which factors as
\[
\mu: (\overline{X_0})\to(\widehat{X_0})\xrightarrow{w}\widehat{X}_0.
\]
We have $\delta(\widehat{X}_0)=\delta(\widehat{X}_0)^\red-\varepsilon(\widehat{X}_0)$ and from the morphisms
\[
\ko_{\widehat{X}_0}\xrightarrow{w^\sharp} w_\ast\ko_{(\widehat{X_0})}\to \mu_\ast\ko_{(\overline{X_0})},
\]
we get $\delta(\widehat{X}_0)^\red = \delta(\widehat{X_0}) + \dim_\C w_\ast\ko_{(\widehat{X_0})}/ \ko_{(\widehat{X}_0)^\red} =r'(X_0)+\dim_\C\Coker(w^\sharp)$ by Lemma \ref{lemma8.2}. 
\medskip

(2) From Proposition \ref{proposition6.14} we get $\varepsilon(X_0)-\varepsilon (X_t)=\varepsilon(X_0^{>1})-\varepsilon(X_t^{>1})$ since $\Coker(\varphi_t)=0$ for $t\in T$ in this case. Since $X^{>1}$ is reduced, $X^{>1}_t$ is reduced by Proposition \ref{prop1} (1)(iv). As it has no $0$--dimensional part, we get $\varepsilon(X_t^{>1})=0$.\hfill\qed

\bigskip

As an immediate consequence of Theorem \ref{theorem6.15} we get:

\begin{Corollary}
	\begin{enumerate}
		\item [(1)] The functions $t\mapsto\delta(X_t)$ and $t\mapsto\varepsilon(X_t)$ are upper semicontinuous on $T$.
		\item [(2)] If $r_1(X,x)=0$, then $(\overline{X}_0, \bar{x})$ is reduced and $\delta(X_0)-\delta(X_t)=\delta(\overline{X}_0)\geq 0$, with $\delta(\overline{X}_0)=0$ iff $\overline{X}_0$ is normal.
	\end{enumerate}
\end{Corollary}
\bigskip

From Theorem \ref{theorem6.15} we deduce also easily $\delta$--constant criteria for the simultaneous normalization of $1$--parametric families of isolated non--normal singularities.

\begin{Theorem} \label{Theorem6.17} Let $f: (X,x)\sr (\C,0)$ be flat with $(X_0,x)$ an isolated non-normal singularity and $\dim (X_0,x) \geq 1$.
	Then the following are equivalent.
	\begin{itemize}
		\item[(1)] $f$ is $\delta$--constant,
		\item[(2)] $f^{red}$ is $\delta$--constant,
		\item[(3)] $f^{>1}$ is $\delta$--constant,
		\item[(4)] $f^{>1}$ is equinormalizable.
	\end{itemize}
\end{Theorem}

\noindent \textbf{Proof}:   $(1)\td (2)$ resp. $(1)\td (3)$ follows immediately  from Theorem \ref{theorem6.15}.
Moreover, $f$ is $\delta$--constant if and only if $\delta(\gt{X}_0^{>1}) = 0$ and this is equivalent to $\gt{X}_0^{>1}$ normal, since $\gt{X}_0^{>1}$ is reduced by Proposition \ref{lemma5} (1). Hence $(4) \Sr (1)$. Conversely, if $f$ is $\delta$--constant, then  $\gt{X}_0^{>1}$ is normal. It follows that $\gt{{f}^{>1}}$ is normal since it is flat and $\gt{X}_t^{>1} \mtn X_t^{>1}$ is the normalization by Proposition \ref{prop1} (2), which implies that $f^{>1}$ is equinormalizable.  \hfill\qed
\bigskip

\begin{Theorem} \label{theorem1} Let $f : (X,x) \mtn (\C,0)$ be flat with $(X_0,x)$ an INNS of dimension $\geq 1$.
	Then $f$ is equinormalizable if and only if $f$ is $\delta$--constant and $r_1(X,x) = 0$.
\end{Theorem}

\noindent\textbf{Proof}:  
\rm {The necessity of $r_1(X,x) = 0$ for equinormalizability follows from Lemma \ref{lemma1.1} (2). If $r_1(X,x)=0$ then $f^{red}=f^{>1}$ and the fibres of $f$ are generically reduced by Proposition \ref{prop1} (1). The result follows from Corollary \ref{Theorem6.17} and the fact that $f$ is equinormalizable if and only if $f^{red}$ is equinormalizable by Lemma \ref{lemma5.5}.} \hfill\qed
\bigskip

\begin{Corollary}\label{cor6.18}
	Let $f : (X,x) \mtn (\C,0)$ be active with $(X_0,x)$ a reduced INNS of dimension $\geq 1$.
	Then $f$ is equinormalizable iff it is $\delta$--constant.
\end{Corollary}

\noindent\textbf{Proof}: 
\rm By Lemma  \ref{lemma1} $(X,x)$ is reduced and $f$ is flat. Then $r_1(X,x) = 0$ by Proposition \ref{lemma5} (iv) and the result follows from Theorem \ref {theorem1}. \hfill\qed
\bigskip

Note that $\dim (X_0,x) = 0$ implies that $f: (X,x) \mtn (\C,0)$ is equinormalizable if and only if $f$ is an isomorphism (cf. Remark \ref{remark1.1} (2)).
\bigskip

\begin{Theorem}\label{theo6.16}
	Let $f: X\to T$ be a flat morphism of complex spaces with $T$ a $1$--dimensional complex manifold, such that the non--normal locus of $f$ is finite over $T$. Then the following are equivalent:
	\begin{enumerate}
		\item [(i)] $f$ is equinormalizable.
		\item [(ii)] $\delta(X_t)$ is constant on $T$ and the $1$--dimensional part $X^1$ of $X$ is smooth and does not meet the higher dimensional part $X^{>1}$.
	\end{enumerate}
\end{Theorem}

\noindent\textbf{Proof:} Let $\nu:\overline{X}\to X$ be the normalization and $\bar{f}=f\circ \nu$. Then $f$ is equinormalizable iff $\bar{f}$ is normal (which is a local condition on $\overline{X}$), and if $\nu:\overline{X}_t\to X_t$ is bimeromorphic for $t\in f(X)$.

We have by Proposition \ref{Prop6.1} that $\NNor(X)\subset \NNor(f)$ and that $\NNor(X_t)=\NNor(f)\cap X_t$. Hence 
\[
\begin{array}{c}
\overline{X}\smallsetminus(\nu^{-1}\NNor(f))\to X\smallsetminus\NNor(f),\text{ and }\\[0.5ex]
\overline{X}_t\smallsetminus (\nu^{-1}\NNor(X_t)) \to \overline{X}_t\smallsetminus\NNor(X_t)
\end{array}
\]
are isomorphisms. \\
It follows that $\overline{X}_t\to X_t$ is bimeromorphic iff $\dim(\NNor X_t, x)<\dim(X_t, x)$ for all $x\in X_t$, which is a local condition on $X$. Hence  we have to show the equivalence of (i) ad (ii) only for the germs $f:(X,x)\to (T, f(x))$, for each $x\in X$.

If $x\in X^{>1}$ then $\dim(X_{f(x)}, x)\geq 1$ and we get from Theorem \ref{theorem1} that $f:(X,x)\to (T, f(x))$ is equinormalizable iff $f$ is $\delta$--constant and $x\not\in X^1\cap X^{>1}$. In particular, $X^1\cap X^>=\emptyset$.

If $x\in X^1$ then $(X^1,x)=(X,x)$, $\dim(X_{f(x)}, x)=0$ and $f:(X,x)\to (T, f(x))$ is equinormalizable iff $f$ is an isomorphism. This implies that $X$ is smooth at $x$ and $\delta(X_t)=-1$ for $t$ in a neighbourhood of $f(x)$. 
Conversely, if $(X^1, x)=(X,x)$ is smooth, then $ (X,x) = (\overline X,\nu^{-1}(x))$ and $ f = \bar f$ is normal at $x$. Hence 
$f$ is regular on $U\smallsetminus\{x\}$ for a neighbourhood $U$ of $x$ in $X$. Since $f$ is flat, $\dim_\C \ko_{X_{f(y),y}} =1$  and hence $\delta(X_{f(y)}, y)=-1$ for $y\in U$.\hfill\qed
\bigskip

\begin{Remark}
	\rm (1) Theorem \ref{theo6.16} shows that the equidimensionality assumption by Chiang-Hsieh and Lipman in Theorem \ref{Theo10} and by Koll$\acute{\text{a}}$r in Theorem \ref{Theo-ProjMorph} is not necessary, at least for $1$--parametric families of INNS. The only {\em geometric} necessary condition for equinormalizability is that of \ref{theo6.16} (ii).
	
	(2) Theorem \ref{theo6.16} applies of course to proper analytic, and hence algebraic, morphisms $f$ with $\NNor(f)$ finite over $T$. Koll$\acute{\text{a}}$r's result in Theorem \ref{Theo-ProjMorph} does however not imply Theorem \ref{theo6.16} since in our case $f:X\to T$ need not be algebraic. The paper \cite{BF79} contains examples of flat analytic morphisms $f:(X,x)\to (\C, 0)$ with all fibres (homogeneous) singularities of dimension $n$, for any $n\geq 0$, such that $\ko_{X,x}$ is not the analytic ring of some algebraic $\C$--scheme. If $n\geq 1$ (resp. $\geq 2$) then $(X,x)$ can even be chosen reduced (resp. normal).
\end{Remark}
\bigskip

\section{Connected components of the Milnor fibre}

The main purpose of this section is to study topological properties of the general fibre of a flat morphism $f: (X,x) \to (\C,0)$ with special fibre $(X_0,x)$ an INNS of arbitrary dimension. We determine the number of connected components of the Milnor fibre $X_t=f^{-1}(t), \ t\in T\smallsetminus \{0\}$, for a good representative $f:X\to T$. Some of these properties have already been proved in \cite{BrG90} for generically reduced curves. 
\medskip

We continue to use the notations of Section 5 and in particular consider the induced morphisms $f^{red}, \bar {f}, \widehat{f}$ and the invariants $\varepsilon(X_t), \delta(X_t), r'(X_t)$ and $\mu(X_t)$.
In addition, we denote for any topological space $Y$ by
\[
b_i(Y):=\dim_\C H^i(Y,\C),
\]
the \textbf{i--th-Betti-number} and by
\[
\chi (Y):=\sum\limits_{i\geq 0}(-1)^ib_i(Y),
\]
the \textbf{topological Euler characteristic}. $b_0(Y)$ is the number of connected components of $Y$ and $b_i(Y)=0$ for $i>n$, if $Y$ is a Stein complex space of dimension $n$. In our situation, for a good representative, $X$ is a Stein space and the fibres are also Stein spaces.\\

Before we consider families of INNS, we start with an arbitrary morphism $f: (X,x) \to (\C,0)$ and prove a connectedness result for the Milnor fibre (Theorem \ref{Prop7.4} and Corollary \ref{Cor8.4}), which seems to be new.

For  $\dim(X,x) = 2$ this was proved in \cite{BuG80}, by proving coherence and local freeness of a certain hypercohomology group. The method can easily be generalized to higher dimensions if the singular set of $f$ is finite over  $T$ and, with more effort, also to our more general situation. In a personal correspondence Helmut A. Hamm \cite{Ha16} and J. Fernández de Bobadilla \cite{Bo16} proposed easier topological proofs for $(X,x)$ irreducible. Hamm's proof goes by induction on the dimension and uses a weak Lefschetz  theorem (he obtains also some information about fundamental groups), while Bobadilla proposed to use the monodromy. The final general result, as presented in Proposition \ref{Prop7.4}, is based on Bobadilla's idea and resulted from discussions of the author with Hamm. 

We need also a general fibration theorem due to L\^{e} D\~{u}ng Tr$\acute{\text{a}}$ng (\cite{Le77}, see also \cite{Le76} for a detailed account), saying that for an arbitrary morphism of complex germs $f: (X,x) \to (\C,0)$ and for a good representative $f: X \to T$, 
              $$f: X^* := X \setminus X_0 \ \to \ T \setminus \{0\} := T^*$$ 
is a topological fibre bundle. 
The geometric monodromy $\rho$ of this fibre bundle permutes the irreducible components of the Milnor fibre $F:=X_t, \ t \in T^*$. 
\medskip

\begin{Lemma}\label{Lem8.1}  Let $f: X \to T$ be a good representative  of  $f: (X,x) \to (\C,0)$.  
\begin{enumerate}
\item [(1)] If $f$ is flat, then \\
$\Sing(f) = \Sing(X_0) \cup \Sing(X),  \ \Sing(X_t) = X_t \cap \Sing(f)$, $ t\in T$.
\item [(2)] If the monodromy of the fibre bundle $f: X^* \to T^*$ acts trivially on the irreducible components of $F$, then there is a bijection between the irreducible components of $F$ and those of $X$ that are not contained in $X_0$.
\end{enumerate}
\end{Lemma}

\noindent\textbf{Proof}: (1) Since $f$ is flat,
$\Sing(X_t) = X_t \cap \Sing(f)$ for $t\in T$ by \cite{GLS07}, Theorem I.1.115 (c). $(X,z)$ is smooth if $f$ is smooth at $z$ by Proposition \ref{Prop6.1}, hence $\Sing(X_0) \cup \Sing(X) \subset \Sing(f)$.  

We claim that  $X_t \setminus \Sing(X)$ is smooth for $t \neq 0$. Otherwise $\Sing(f)$ would contain a germ of a curve $(C,0) \subset (X,0)$, finite over $T$ and not contained in $\Sing(X)$. At $z \in C \smallsetminus {0}$, $X$ is smooth and $f$ has a singularity. Hence for $t$ close to $f(z)$ the fibre $X_{t}$ is smooth at $z'$ close to $z$ (cf. \cite {GLS07}, Lemma I.2.4), contradicting that $C$ is a curve, which is finite over $T$. This proves (1).

(2) Since the statement is purely set-theoretically, we may assume that $X$ is reduced. If $f(X) = 0$, the statement is trivially true and if $X$ contains irreducible components that are mapped to $0$ by $f$, we omit them as they do not affect the Milnor fibre. We may therefore assume that $f$ is active and hence flat (since $X$ is reduced). For any reduced complex space $X$, the irreducible components are the topological closure of the connected components of $X \setminus \Sing(X)$. Moreover, $X$ is irreducible iff $X \setminus A$ is connected for any subspace $A$ of codimension 1. 

It follows that the irreducible components of $X$ that are not contained in $X_0$, correspond to the connected components of $X^* \setminus \Sing(X)$, and from (1) we get that the irreducible components of $F$ correspond to the connected components of $F \setminus \Sing(X)$. The monodromy acts on $X^* \setminus \Sing(X)$ and induces an action $h$ on the integer homology. It gives rise to an exact sequence, the Wang-sequence of the fibre bundle $X^* \to T^*$ (see \cite {Mi68}, Lemma 8.4),

\[
\xymatrix{
H_0(F \setminus \Sing(X)) \ar[r]^{h-id} & H_0(F \setminus  \Sing(X)) \ar[r] & H_0(X^* \setminus \Sing(X)) \to 0. 
}
\]
Since $h$ acts trivially on $H_0(F \setminus \Sing(X))$, the second arrow is an isomorphism, which proves the claim.
\hfill\qed \\

\begin{Theorem}\label{Prop7.4} {\bf (Bobadilla, Greuel, Hamm 2017)}
Let $f:(X,x)\to (\C, 0)$ be a morphism of complex germs and $f: X \to T$ a good representative. 
\begin{enumerate}
\item [(1)] Let $(X,x)$  be irreducible and $\Red(X_0) \neq \emptyset$, i.e. there exist points $y \in X_0$ arbitrary close to $x$ such that  $(X_0,y)$ is reduced. Then the Milnor fibre of $f$ is irreducible.
\item [(2)] Let $(X,x)$  be reducible with irreducible components $(X^i,x), i=1, \ldots r,$ and assume that the intersection graph $G(f)$ is connected.
Then the Milnor fibre of $f$ is connected.


\end{enumerate}
Here $G(f)$ is the graph with vertices $i=1, \ldots r,$ and we join $i \neq j$ by an edge iff $ Red(X_0) \cap X^i \cap X^j \neq \emptyset$, i.e. there exist points $y \in X_0 \cap X^i \cap X^j$ arbitrary close to $x$ such that  $(X_0,y)$ is reduced. 
\end{Theorem}

\noindent\textbf{Proof}: Let $f^{red}$ denote the restriction of $f$ to the reduction $X^{red}$ of $X$. Since $(f^{red})^{-1}(0)$ is reduced at $y$ if $X_0$ is reduced at $y$ and since the statements about the Milnor fibre concern only the reduced structure, we may assume that $X$ is reduced.

(1) If $f(X) = 0$, the statement is trivially true, hence we may assume that f is flat.
The monodromy $\rho$ of the fibre bundle  $X^* \to T^*$ permutes the irreducible components of the Milnor fibre and we denote by $r$ a common multiple of the cardinalities of the orbits of $\rho$. In order to trivialize the monodromy we
consider the base change  $\varphi: T \to T $, $\varphi(t) = t^r$,
\[
\xymatrix{
X_T :=  X\times_T T\ar[rr]^{\pi} \ \ \ar[d]^g && X \ar[d]^f\\
T\ar[rr]_\varphi & &T
}
\]

We claim that $X_T$ is irreducible. Since $\Red(X_0) \neq \emptyset$ also $\Reg(X_0) \neq \emptyset$ and there are smooth points of $X_0$ arbitrary close to 0.
If $Z$ is an irreducible component of $X_T$, then $\pi (Z) = X$ since $X$ is irreducible. It follows that $g^{-1}(0)$ is contained in every irreducible component of $X_T$. Since $g^{-1}(0)$ is isomorphic to $X_0$ it has smooth points $z$ arbitrary close to $(x,0)$. By Lemma \ref{Lem8.1}  $X_T$ is smooth at $z$ and $z$ cannot be in the intersection of two components of $X_T$, i.e. $X_T$ is irreducible.

Since the monodromy of $g$ acts trivially on the irreducible components of the Milnor fibre of $g$ by the choice of the base change, we can apply Lemma \ref{Lem8.1} and get that the Milnor fibre of $g$ and hence the Milnor fibre of $f$ is irreducible.

(2) We proceed by induction on $r \geq 2$. Let $f_i$ denote the restriction of $f$ to $X^i$.
The assumption implies that $Red(X_0) \cap X^i \neq \emptyset$ for $i=1, \ldots r,$ hence the Milnor fibre of  $f_i$ is connected by (1). If $r=2$ then $ \mbox {Red} (X_0)  \cap X^1 \cap X^2 \neq \emptyset$ since $G(f)$ is connected and Proposotion \ref{Prop8.4} implies that the Milnor fibre of $f$ is connected. 

If $r > 2$ let $i$ be a leaf of a spanning tree of $G(f)$ such that the graph  $G'(f)$  obtained from $G(f)$ by deleting $i$ remains connected.
After a renumeration we may assume that $i = 1$. Let $f'$ be the restriction of $f$ to $X' := X^2 \cup \ldots \cup X^r$. Then $G(f') = G'(f)$ is connected and hence the Milnor fibre of $f'$ is connected by induction hypothesis. Since the Milnor fibre of $f_1$ is connected the Milnor fibre of $f$ is connected by Proposotion \ref{Prop8.4}.
\hfill\qed 
\bigskip

\begin{Proposition}\label{Prop8.4}
Let $(X,x)$ be the union of two reduced subgerms $(X^1,x)$  and $(X^2,x)$ and let $f_i$ be the restriction of $f$ to $X^i, i=1,2$. Assume that  the Milnor fibres of $f_i$ are connected and that $ \mbox {Red} (X_0)  \cap X^1 \cap X^2 \neq \emptyset$.
Then the Milnor fibre of $f$ is connected.
\end{Proposition}

\noindent\textbf{Proof}: Note that we do not assume that $X^1$ or $X^2$ are irreducible. 

If $X$ contains irreducible components that are mapped to $0$ by $f$, we omit them as they do not affect the Milnor fibre. Hence we may assume that $f$ is flat. Then the $f_i$  are flat too and we have for the fibres  $X^i_t$ of $f_i$,
$$ X_t = X^1_t \cup X^2_t, \  t \in T \ \text {(as sets).}$$

 We prove now that $f(X^1  \cap X^2) \neq \{0\}$. 
 Choose $y \in \mbox {Red} (X_0)  \cap X^1 \cap X^2$  arbitrary close to $x$ and consider the commutative diagram with exact rows and with vertical arrows induced by multiplication with $f$ ($\ko :=\ko_{\C^n, y }$, $\ko_{X,y} = \ko/I$ and $\ko_{X^i,y}= \ko/I_i$) 

\[
\xymatrix{
0\ar[r] & \ko/I\ar[d]\ar[r] & \ko/I_1\oplus\ko/I_2\ar[d]\ar[r] & \ko/(I_1+I_2)\ar[d]^f\ar[r] & 0\\
0\ar[r] & \ko/I\ar[r] & \ko/I_1\oplus\ko/I_2\ar[r] & \ko/(I_1+I_2)\ar[r] & 0\ .
}
\]
The two first vertical arrows are injective (by flatness) and hence we get from the snake lemma the exact sequence
\[
\xymatrix{ 0\ar[r]  & Ker(f)\ar[r]  & \ko/I+\langle f \rangle \ar[r] & \ko/I_1+\langle f \rangle\oplus\ko/I_2+\langle f \rangle,}
\]
i.e. $Ker(f) = (I_1+\langle f \rangle) \cap (I_2+\langle f \rangle) / (I+\langle f \rangle)$. The numerator defines $(X^1_0,y) \cup (X^2_0,y)$ and the denominator $(X_0,y) $, which coincide as sets. Therefore the radicals of the ideals coincide and since $(X_0,y) $ is reduced, $Ker(f) = 0$. It follows that $f|(X^1  \cap X^2,y)$ is flat, hence open, showing that $f(X^1  \cap X^2) \neq \{0\}$. 

Then $\exists \ z \in X^1  \cap X^2$ such that $f(z) = t \neq 0$ and hence $z \in X^1_t \cap X^2_t$. Since $X^i_t$ is connected, it follows that $X_t$ is connected.
\hfill\qed 
\bigskip

Theorem \ref{Prop7.4} has the following immediate corollary:
\medskip

\begin{Corollary}\label{Cor8.4} Let $f:(X,x)\to (\C, 0)$ be a morphism of complex germs with $(X_0,x)$ reduced. Then the Milnor fibre of $f$ is connected.
\end{Corollary}
\medskip

\begin{Example}\label{Ex8.5} \rm
The following examples are meant to illustrate the previous results. F denotes the Milnor fibre of $f$ and in (1)--(4) $(X,0) = (\mathbb{C}^2,0)$; 
\begin{enumerate}
\item [(1)] $f(x,y):=x^2$: Here $F$ is not connected, $X_0$ is nowhere reduced. 
\item [(2)] $f(x,y):=xy$: $X_0$ is reduced and $F$ is irreducible. 
\item [(3)] $f(x,y):=xy^2$: $F$ is irreducible. Note that $X_0$ is not reduced but has reduced points arbitrarily close to 0.

\item [(4)] $f(x,y):= x^2y^3$: Here $F$ is irreducible as one can see by using the unramified covering $F\to\mathbb{C}^*$, $(x,y)\mapsto y$. Theorem \ref{Prop7.4} is not applicable as $X_0$ is nowhere reduced. However, one can see that the acts monodromy is trivially on $H_0(F)$, hence the irreducibility of $F$ follows also from Lemma \ref{Lem8.1}.

\item [(5)] $X=\{xy=0\}\subset \mathbb{C}^3$, $f(x,y,z)=z$. $X$ is reducible, $X_0$ is reduced and $F$ is connected but not irreducible. Further examples with reduced $X_0$ and connected $F$ with several irreducible components are given in Example \ref{Ex7.5} (1), (2).
\item [(6)] $X=X^1 \cup X^2$ with $X^1=\{x=y=0\}, \ X^2 =\{u=v=0\} \subset \mathbb{C}^4$ and $f(x,y,u,v)=x+y+u+v$. We have $X_0 =  \{x=y=u+v=0\} \cup \{u=v=x+y=0\}$ with an embedded component $\{y=v=x+u=u^2=0\}$ at $0$. Here $X_0 \cap X^1 \cap X^2 = \{0\}$, thus the assumptions of Theorem \ref{Prop7.4} are not satisfied and $F$ consists of 2 disjoint lines. For another example with 3 irreducible components that do not satisfy the assumptions of Theorem \ref{Prop7.4} and with disconnected Milnor fibre see Example \ref{Ex7.5} (4).
\end{enumerate}
Examples (1) and (6) show that the assumption of Theorem \ref{Prop7.4} can in general not be weakened.

\end{Example}

We consider now the behavior of the $\mu$--invariant in a flat family of INNS of arbitrary dimension.\\

Since $\delta-r'=\mu-\delta$, we get immediately from Theorem \ref{theorem6.15} (iv):

\begin{Proposition}\label{Prop7.8} Let $f:(X,x)\to (\C, 0)$ be flat with $(X_0, x)$ an isolated non--normal singularity. For a good representative $f:X \to T$ and $t \in T\smallsetminus\{0\}$ we have
\[
\mu(X_0)-\mu(X_t)= \delta(X_0)-\delta(X_t)+r_1(X,x)+\dim_\C w_\ast\ko_{\widehat{(X_0)}}/ \ko_{(\widehat{X}_0)^\red}-\varepsilon(\widehat{X}_0).
\]
\end{Proposition}

This result was stated in \cite{BrG90}, Lemma 2.4.1, for families of curves without the term $-\varepsilon(\widehat{X}_0)$, which was mistakenly omitted. We will show by examples at the end of this section, that this term may appear and is responsible for $\mu$ being in general not semicontinuous and that $\mu$--constant does not imply $\delta$--constant. Semicontinuity of $\mu$ follows if $\varepsilon(\widehat{X}_0)=0$, i.e. if $\widehat{X}_0$ is reduced. Moreover, in this case $(\mu - \delta)$--constant characterizes simultaneous weak normalizability of $f$:\\

\begin{Theorem}\label{Theorem7.9}
Let $f:(X,x)\to (\C,0)$ be flat, $\dim(X,x) \geq 2$ and $(X_0,x)$ an INNS. If $(\widehat{X}_0, \widehat{x})$ is reduced (equivalently depth $(\widehat{X}, \widehat{x}) \geq 2$) then $r_1(X,x) =0$ and the following holds for a good representative $f:X\to T$.
\begin{enumerate}
\item [(1)] $\mu(X_0)-\mu(X_t)=\delta(X_0)-\delta(X_t)+\dim_\C w_\ast\ko_{\widehat{(X_0)}}/\ko_{\widehat{X}_0}, t\in T\smallsetminus\{0\}$,\\ 
in particular $\mu$ is upper semicontinuous.
\item [(2)] $\mu(X_t)$ is constant iff $\delta(X_t)$ and $r'(X_t)$ are constant for $t\in T$. 
\item [(3)] $\mu(X_t)-\delta(X_t)$ is constant for $t\in T$ iff $\omega: \widehat{X}\to X$ is a simultaneous weak normalization of $f$.
\item [(4)] $(\widehat{X}_0, \widehat{x})$ is reduced in the following cases:
\begin{enumerate}
\item [(i)] $(X_0, x)$ is reduced.
\item [(ii)] $X$ is locally irreducible.
\end{enumerate}
\end{enumerate}
\end{Theorem}

\noindent\textbf{Proof}: 
Since  $\widehat{f}$ is flat we get from  Proposition \ref{lemma5} (iv)  that $(\widehat{X}_0, \widehat{x})$ is reduced iff depth $(\widehat{X}, \widehat{x})\geq 2$ and  $r_1(\widehat{X}, \widehat{x}) = r_1(X,x)=0$. Hence (1) follows from Proposition \ref{Prop7.8}.
\medskip

(2) If $\delta(X_t)$ and $r'(X_t)$ are constant, then $\mu(X_t)$ is constant by the definition of $\mu$. If $\mu(X_t)$ is constant, then $\delta(X_t)$ is constant by (1) and by the semicontinuity of $\delta$.
\medskip

(3) By (1), $\mu(X_t)-\delta(X_t)$ is constant iff $w_\ast \ko_{\widehat{(X_0)}}=\ko_{\widehat{X}_0}$, i.e. iff $\widehat{X}_0$ is weakly normal. Then $w: \widehat{X}_0\to X_0$ is the weak normalization of $X_0$, and hence $\widehat{X}_t \to X_t$ is the weak normalization for all $t$.
\medskip

(4) (i) If $(X_0, x)$ is reduced, then depth $(X,x) \geq 2$ by Proposition \ref{lemma5} (iv) and hence depth $(\widehat{X}, \widehat{x}) \geq 2$ (by \cite{BF93}). 

(ii) If $X$ is locally irreducible, then $r_1(X,x) = 0$ and the weak normalization coincides with the normalization $\overline{X}$. Hence $\widehat{X}_0=\overline{X}_0$ is reduced by Proposition \ref{lemma5} (i).\hfill\qed
\bigskip

The most drastic difference between families of reduced and non--reduced INNS is, that in the non--reduced case the Milnor fibre $X_t$ may have several connected components. Proposition \ref{Prop7.5} gives a precise answer about the number of connected components. Using a primary decomposition and computations such as in Example \ref{Ex5.6}, this number can be effectively computed.
\medskip

Part (1) of the following result was proved for isolated curve singularities in \cite{BrG90}. 

\begin{Proposition}\label{Prop7.5} Let $f:(X,x)\to (\C, 0)$ be active, $\dim(X,x) \geq 2$, $(X_0,x)$ an isolated non-normal singularity and  $X_t$ the Milnor fibre of $f$.
\begin{enumerate}
\item [(1)] If $r_1(X,x)=0$, then $b_0(X_t)=b_0(X\smallsetminus\{x\}) $.
\item [(2)] If $r_1(X,x)\neq 0$, then $b_0(X_t)=b_0(X^{>1} \smallsetminus\{x\}) + \varepsilon(X^1_0,x)$, \\
with $\varepsilon(X^1_0,x) = \dim_\C\ko_{X^1,x}/\langle f \rangle = \sharp\{\text{isolated points of }X_t, t\neq 0\} = r_1(X,x)$.
\end{enumerate}
\end{Proposition}

\noindent\textbf{Proof}: Since $X$ and $X^\red$ as well as $X_t$ and $X_t^\red$ are topologically the same spaces, we may assume that $X$ is reduced and that $f$ is flat. 

(1) Let $r_1(X,x)=0$. 
If $(X,x)$ is irreducible, then  $X_t$ is connected for $t \neq 0$ by Proposition \ref{Prop7.4} since $\dim(X_0,x) \geq 1$ and $X_0 \smallsetminus \{x\}$ is reduced.

Now let $X$ be reducible with $X_i, i=1, \ldots ,r,$ the irreducible components, $f_i$ the restriction of $f$ to $X_i$, and $X_{i,t}$ the Milnor fibre of $f_i$. Then $X_{t}$ is the union of the $X_{i,t}$ and  $X_{i,t}$  is connected by the previous step. We claim that  $X_{i,t} \cap X_{j,t} \neq \emptyset$ for $i \neq j$ iff $X_i \cap X_j  \neq \{x\}$ . Note that $X_i \cap X_j \cap X_0 = \{x\}$ since $X$ is normal at $y \in X_0 \smallsetminus \{x\}$. Therefore $X_i \cap X_j  \neq \{x\}$ iff  $f(X_i \cap X_j)$ is at least 1-dimensional, i.e. $X_{i,t} \cap X_{j,t} = (f | X_i \cap X_j)^{-1} (t) \neq \emptyset$. This proves the claim.

To see (2), let $X^1$ be the union of the $1$--dimensional irreducible components of $X$. By Remark \ref{rem2.9} $r_1(X,x)=\sharp\{\text{isolated points of }X_t\} = \varepsilon (X^1_0)$. 
The result follows from (1), since the number of connected components of $X_t$ is the number of its positive dimensional components plus the number of its isolated points. \hfill\qed \\

Apart from the number of connected components, we cannot say much about the topology of the  Milnor fibre for an arbitrary INNS. This is different for families of generically reduced curves to which we switch in the next section.

\section{Topology of families of generically reduced curves}

We consider families of generically reduced curve singularities $f:(X,x)\to (\C, 0)$ with $(X,x)$ a surface singularity, i.e. a 2--dimensional complex germ. We have
\[
(X^\red, x)=(X^2,x)\cup(X^1, x)
\]
with $(X^i, x)$ the pure $i$--dimensional part of $(X^\red, x)$. $(X,x)$ is pure $2$--dimensional if $(X^\red, x)=(X^2,x)$, i.e. if $r_1(X,x)=0$ but $(X,x)$ may not be reduced.\\

Let us now investigate $\mu$--constant families. For families of reduced curves $\mu$--constant along a section is equivalent to topological triviality by Theorem \ref{Theo3.9}. We show by an example below that this is in general not true for families of non--reduced curve singularities.
Nevertheless, the nice thing about $\mu$ is, that it determines the Euler characteristic of the Milnor fibre in a family of generically reduced curves.\\

\begin{Theorem}\label{Theorem7.10}Let $f:(X,x)\to (\C,0)$ be flat with $(X_0, x)$ an isolated curve singularity and $f:X\to T$ a good representative. Then the following holds:
\begin{enumerate}
\item [(1)] $\delta(X_0)-\delta(X_t)=\delta(X_0^\red)-\delta(X_t^\red)=\delta(\overline{X}_0^2) \geq 0$ for $t\in T\smallsetminus \{0\}$.
\item [(2)] $\mu(X_0)-\mu(X_t)=\mu(X_0^\red)-\mu(X_t^\red)=1-\chi(X_t)$ for $t\in T$.
\end{enumerate}
\end{Theorem}

\noindent\textbf{Proof}: (1) is a special case of Theorem \ref{theorem6.15}.

(2) The equality $\mu(X_0)-\mu(X_t)=1-\chi(X_t)$ was proved in \cite{BrG90}, Satz 3.1.2. The proof is completely topological using the reducedness of the Milnor fibre of $\bar f$. Since $\chi(X_t)=\chi(X_t^\red)$ for all $t$, the second equality follows (it follows of course also easily from (1) and the definition of $\mu$).\hfill\qed\\

\begin{Corollary}\label{Cor7.12}
With the assumption of Theorem \ref{Theorem7.10} we have for $t\in T$:
\begin{enumerate}
\item [(1)] $\mu(X_t)$ is constant iff $\chi(X_t)$ is constant.
\item [(2)] The following are equivalent:
\begin{enumerate}
\item [(i)] 	$\mu(X_t)-b_0(X_t)$ is constant,
\item [(ii)] 	$b_1(X_t)=0$,
\item [(iii)]   each connected component of $X_t$ is contractible.
\end{enumerate}	
\end{enumerate}
\end{Corollary}

\noindent\textbf{Proof}: The equivalence of (iii) to the other statements in (2) follows as in \cite{BuG80}, Theorem 4.2.4.\hfill\qed\\

Even if $\mu$ lacks, if $(X_0,x)$ is not reduced, some of the well--known properties (like semicontinuity),  the above results show that $\mu$ always controls the Euler characteristic of the Milnor fibre of an arbitrary isolated curve singularity. 
This is somewhat surprising as $(X_0,x)$ may have a rather complicated embedded component.\\

In the following theorem  we use the notion of weak simultaneous resolution for families with non-reduced fibres. 
If $X$ is a non-reduced space we say that $\pi:\widetilde{X}\to X$ is a \textbf{resolution of singularities} if $\widetilde{X}$ is reduced and if the induced map $\pi:\widetilde{X} \to X^\red$ is a resolution of singularities (in the classical sense). Then the notion of very weak (resp.  weak, resp.  strong) simultaneous resolution from Definition \ref{Def2.1} carries over verbatim.\\

For investigating topological triviality, we have to assume that $X_t, t\in T$, has only one singular point along a section.

\begin{Theorem}\label{Theorem7.14}
Let $f:(X,x)\to (\C,0)$ be a deformation with section of the isolated curve singularity $(X_0, x)$ and $f:X\to T$ a good representative with section $\sigma:T\to X$, such that $X_t\smallsetminus\sigma(t)$ is smooth for $t\in T$. Then the following are equivalent:
\begin{enumerate}
\item [(i)] $(X,x)$ is pure $2$--dimensional and $f:X\to T$ admits a weak simultaneous resolution.

\item [(ii)] $(X,x)$ is pure $2$--dimensional and $\delta(X_t, \sigma(t))$ and $r(X,\sigma(t))$ are constant.

\item [(iii)] $\mu(X_t)$ is constant and $X_t$ is connected for $t\in T$.

\item [(iv)] $X_t$ is connected and $H^1(X_t, \C)=0$ for $t\in T$.

\item [(v)] $X_0$ and $X_t$ are embedded topologically equivalent for $t\in T$.

\item [(vi)] $f:X\to T$ is topologically trivial.
\end{enumerate}
\end{Theorem}

\noindent\textbf{Proof}: We prove the following implications

\[
\xymatrix{(i)\ar@{=>}[d] & (ii)\ar@{=>}[l] & (iii)\ar@{=>}[l] \ar@{<=>}[d]\\
(vi)\ar@{=>}[rr]&& (iv) &(v)\ar@{<=>}[l]
}
\]

(ii)$\Rightarrow$(i): Since $r_1(X,x)=0, \   \delta(X_t)=$ constant implies very weak simultaneous resolution by Theorem \ref{theorem1}. Since $r(X_t, \sigma(t))=\#\nu^{-1}(\sigma(t))$ is constant, $\nu^{-1}(\sigma(T))^\red\cong\nu^{-1}(x)^\red\times T$ over  $T$, and we get (i).
\medskip

(i) $\Rightarrow$ (vi): the proof is the same as in \cite{BuG80}, (4) $\Rightarrow$ (6).
\medskip

(vi) $\Rightarrow$ (iv): trivial
\medskip

(iv) $\Leftrightarrow$ (v): (v) $\Rightarrow$ (iv) is trivial, while the proof for (iv) $\Rightarrow$ (v) is the same as in \cite{BuG80}, (3) $\Rightarrow$ (5): it only needs the connectedness of the Milnor fibre but not the reducedness of $(X_0, x)$.
\medskip

(iv) $\Leftrightarrow$ (iii): follows from Corollary 7.11 (2).
\medskip

(iii) $\Rightarrow$ (ii): Since (iii) $\Rightarrow$ (iv) was already proved, we can use that there is a homeomorphism $h:(B, X_0, x)\xrightarrow{\approx} (B, X_t, \sigma(t))$ for some embedding $X\subset B\times T$, with $f$ the restricton of the second projection. Hence $h(X_0\smallsetminus \{x\})=X_t\smallsetminus \sigma(t)$ and $X_0\smallsetminus\{x\}$ and $X_t\smallsetminus\sigma (t)$ have the same number of connected components. Since $X_t$ is connected and $X_t\smallsetminus\sigma (t)$ is smooth, the number of connected components of $X_t\smallsetminus\sigma(t)$ is equal to the number of (global) irreducible components of $X_t$, $t\in T$, and all irreducible components of $X_t$ meet at $\sigma(t)$. Therefore $r(X_t, \sigma (t))$ is constant. Together with $\mu(X_t, \sigma (t))=$ constant we get $\delta(X_t, \sigma(t))=$ constant. Moreover, $X_t$ connected implies that $X_t$ has no isolated points, i.e. $r_1(X,x)=0$, and hence (ii) follows. \hfill\qed \\

\begin{Remark}\label{Rem7.14} \rm
The equivalence of (i), (ii) and (vi) was proved by L\^{e} C\^{o}ng-Tr\`{i}nh in \cite{LC15}, Theorem 4.5, for $(X,x)$ pure dimensional and reduced, using Corollary \ref{Cor7.15} below. Note however, that the constancy of the Milnor number is not sufficient to guarantee (ii), see Example \ref{Ex7.5} (4). 
\end{Remark}
\bigskip

We need now a local Bertini theorem for a generic linear projection $f:(X,0)\to (\C, 0)$, due to Flenner \cite{Fl77}.  $f$ is called a \textbf{generic linear projection} if it is the restriction to $(X,0)\subset(\C^n, 0)$ of a linear projection
\[
\pi:(\C^n,0)\to (\C,0)\ ,\ (x_1, \ldots, x_n)\mapsto (\lambda_1x_1+ \ldots +\lambda_nx_n),
\]
with $\lambda = (\lambda_1,\ldots,\lambda_n) \in \C^{n}$ generic. $f_\lambda^{-1}(0)$ is called a {\bf generic hyperplane section}.
\medskip

\begin{Lemma} \label{Lem7.4} Let $(X,0)\subset(\C^n, 0)$ be a complex germ of positive dimension and $f:(X,0)\to (\C, 0)$ a generic linear projection. Then the following holds for a good representative $f:X\to T$.
\begin{enumerate}
\item [(i)] $f$ is active on $X$.
\item [(ii)] Let $(X,0)$ be reduced. Then $f$ is flat on $X$ and for all $t\in T$ the fibre $X_t$ is smooth at $y$ iff $y$ is a smooth point of $X$, i.e. $\Sing(f) = \Sing(X)$. 
\end{enumerate}
\end{Lemma} 
\medskip

\noindent\textbf{Proof}:  (i) $f$ is active iff no irreducible component of $(X,0)$ is contained in the linear subspace $\pi^{-1}(0)$, which can be achieved for a generic $\pi$ since $\dim(X,0)>0$.

(ii) The flatness is an immediate consequence of (i). That $f^{-1}(0) \cap Reg(X)$ is smooth is proved in  \cite{Fl77}, (4.1) Satz. $f^{-1}(t) \cap Reg(X)$ is smooth for $t \neq0$ by Lemma \ref{Lem8.1} (1). Conversely, if $(X_t,y)$ is smooth then $(X,y)$ is smooth by Proposition \ref{Prop6.1} (2).
\hfill\qed\\

From Theorem \ref{Theorem7.14} we deduce

\begin{Corollary}\label{Cor7.14}
Let $(X,x)$ be a reduced, pure dimensional surface singularity with smooth 1-dimensional singular locus. A good representative $f: X \to T$ of a generic linear projection makes $X$ a topologically trivial family of curves, iff the normalization $\overline{X}$ of $X$ becomes an analytically trivial family of smooth curves. 
\end{Corollary}

\noindent\textbf{Proof}: 
We can choose a generic $f$ such that $\Sing(X)$ is finite and smooth over $T$, i.e. there is a section $\sigma :T \to \Sing(X)$ of $f$. Moreover, $f$ is flat on $X$ and $f^{-1}(t) \smallsetminus \sigma (t)$ is smooth by Lemma \ref{Lem7.4}. Now we can apply Theorem \ref{Theorem7.14}.\hfill\qed
\bigskip

The following consequence is due to J. Fernández de Bobadilla, J. Snoussi and M. Spivakovsky (\cite{BSS16}, Theorem 4.4).

\begin{Corollary}\label{Cor7.15} {\bf (Bobadilla, Snoussi, Spivakovsky)}
If a reduced, pure dimensional surface singularity $X$ with smooth 1-dimensional singular locus is a topologically trivial family of curves, then $\overline{X}$ is smooth.
\end{Corollary}
\medskip

As another corollary we can characterize strong simultaneous resolutions for families of generically reduced curves, which follows from Theorem \ref{Theorem7.14} and Lemma \ref{Lem7.16}. Let $mt$ denote the Hilbert-Samuel multiplicity. 

\begin{Corollary}\label{Cor7.16}
With the assumptions of Theorem \ref {Theorem7.14},
$f:X\to T$ admits a strong simultaneous resolution iff any of the conditions $(i) \ldots (vi)$ hold and if $mt(X_t,\sigma(t))$ is constant for $t\in T$.	
\end{Corollary}
\medskip

\begin{Remark}\label{Rem9.9} \rm
(1) In \cite{BT08} J. J. Nu\~no--Ballesteros and J. N. Tomazella study the Milnor number $\mu(f)$ for a finite morphism $f:(X_0,0) \to (\C,0)$, \ $(X_0,0)$ a reduced curve singularity (in the sense of V. Goryunov, D. Mond and D. van Straten).  They define $m_1(X_0,0) := \mu(f)$ for $f$ a generic linear projection and prove (based on \cite{BuG80}) a "Greuel -- L\^e" type formula
 $$m_1(X_0,0) = \mu(X_0,0) + mt(X_0,0) -1.$$
They show that $m_1(X_t,\sigma(t))$ being constant in a family $f:X \to T$ as above, with $X_0$ reduced, is equivalent to Whitney equisingularity. 
\medskip

(2) If we define $m_1(X_0,0)$ by the above formula for a generically reduced curve singularity and with  $\mu(X_0,0)$ as in Definition \ref{Def5.2}, we get from Corollary \ref{Cor7.16}:

If $f:(X,x)\to (\C,0)$ is a deformation with section $\sigma$ of the generically reduced curve singularity $(X_0, x)$ such that $X_t\smallsetminus\sigma(t)$ is smooth then 
\[
 f:X \to T \text{ admits a strong simultaneous resolution}
  \]
 \[
  \Leftrightarrow  X_t   \text{ is connected and }  
 m_1(X_t,\sigma(t))  \text{ is constant.}
 \]
 It would be interesting to study properties $\mu$ and $m_1$ for functions on a generically reduced curve or, more generally, on an INNS of arbitrary dimension.
\end{Remark}
\medskip

The following lemma is well-known for reduced curve singularities. Using that $mt(X_t,\sigma(t)) = mt(X^{red}_t,\sigma(t))$ by the additivity formula of Serre (cf. \cite{Se65}, V-3), the proof is the same as in the reduced case. 

\begin{Lemma}\label{Lem7.16}
Let $f:(X,x)\to (S,0)$ be a deformation with section $\sigma$ of an isolated curve singularity $(X_0, x)$ with $(S,0)$ reduced, $(X,x)$ pure dimensional and $X_s \smallsetminus \sigma (s)$ smooth for $s \in S$. Then a weak simultaneous resolution along $\sigma (S)$ is strong iff $\mt(X_s, \sigma (s))$, $s\in S$, is constant.
\end{Lemma}

\noindent\textbf{Proof}:
\rm{ Let $\pi:\widetilde{X}\to X$ be a weak simultaneous resolution of $f:X\to S$ and let $(X_0^1, x), \ldots, (X^r_0, x)$ be the branches of $(X_0, x)$. Then $r$ is also the number of irreducible components of $(X_s, \sigma (s))$, $s\in S$, and of $(X,x)$ by Theorem \ref {Theorem7.14}. 

Then $\pi$ can be given as a deformation of the parametrization of $(X_0, 0) \subset (\C,0)$,
	\[
	\begin{array}{c}
	(\widetilde{X}, \widetilde{0})=\overset{r}{\underset{i=1}{\coprod}}(\widetilde{X}^i, \widetilde{0}_i)=\overset{r}{\underset{i=1}{\coprod}}(\C\times S, 0)\overset{\pi}{\twoheadrightarrow} (X,0)\hookrightarrow (\C^n\times S,0)\\(t_i, s)\mapsto(X_{i,1}(t_i, s), \ldots, X_{i,n}(t_i, s), s),
	\end{array}
	\]
	with 
	\[
	X_{i,j}(t_i, s)= x_{i,j}(t_i)+s\cdot a_{i,j}(t_i, s), i=1, \ldots, r, j=1, \ldots, n,
	\]
	such that $t_i\mapsto (x_{i,1}(t_i), \ldots, x_{i,n}(t_i))$ is a parametrization of $(X_0^i, 0)$.
	
	Then $W=\sigma(S)=\{0\}\times S$ and $\widetilde{W}=\pi^{-1}(W)$ is the disjoint union of the $\widetilde{W}_i$ defined by the ideal $I_i=\langle X_{i,1}(t_i, s), \ldots, X_{i,n}(t_i, s)\rangle$ in $\ko_{\widetilde{X}^i, \widetilde{0}_i}=\C\{t_i, x\}$. $\pi^{-1}(0)=\{\widetilde{0}_1, \ldots, \widetilde{0}_r\}$ is defined by the ideals $\langle x_{i,1}(t_i), \ldots, x_{i,n}(t_i)\rangle=t^{m_i}\C\{t_i\}$ in $\C\{t_i\}$ where $m_i=\min\{\ord(x_ {i,1}), \ldots, \ord(x_{i,n})\}$ is the multiplicity of $(X_0^i, 0)$. Hence $\widetilde{W}\cong\pi^{-1}(0)\times S$ iff $I_i=t^{m_i}\C\{t_i, s\}$ i.e. iff $m_i=\min\{\ord_{t_i}X_{i,1}$\\$(t_i, s),\ldots,\ord_{t_i}X_{i,n}(t_i, s)\}$. But this is equivalent to $m_i=\min\{\ord X_{i,1}(-, s),$\\$ \ldots, \ord X_{i,n}(-, s)\}=\mt(X_s^i, 0)$ for each fixed $s\in S$ sufficiently close to $0$. 
	
	Since $\mt(X_s, 0)=\mt(X^{red}_s,0) =\sum\limits^r_{i=1}\mt(X_s^i, 0)$ and since the multiplicity is semicontinuous, $\widetilde{W}\cong \pi^1(0)\times S$ iff $\mt(X_s, \sigma (s))$ is constant. This proves the lemma.}\hfill\qed
\bigskip

We illustrate the previous results by several examples. 

\begin{Example}\label{Ex7.5} \rm
In example (1), (2), (4) $f:X\to \C$ is the projection to the $t$--axis. 
\begin{enumerate}
	\item [(1)] Let $X \subset \C^4$ be the reduced surface consiting of 3 planes and  defined by the ideal $\langle xz, xy, yz+zt\rangle$ with depth $(X,0)=2$. $X_0$ is the  reduced curve singularity consisting of the 3 coordinate axes in $ \C^3$, while $X_t$ consists of 3 lines meeting in two ordinary nodes. This family is $\mu$--constant and $\delta$--constant with $\mu(X_t)=\delta(X_t)=r'(X_t)=2$.
	
	\item [(2)] Let $X\subset\C^4$ be defined by $\langle xz, xy, yz+yt+zt\rangle$, again reduced, pure $2$--dimensional of depth $2$. $X$ consists of a plane and a surface with an $A_1$--singularity, meeting in the $t$--axis. We have $X_0$ as in (1) but for $t\neq 0$ the Milnor fibre $X_t$ is the union of two smooth curves, meeting in an $A_1$--singularity in $0$ (we smoothed away one $A_1$--singularity in (1)), hence $\mu(X_t)=\delta(X_t)=r'(X_t)=1$. The family is $(\mu-\delta)$--constant. $X$ and all fibres are weakly normal, and $\id:X\to X$ is the weak normalization of $f$, in accordance with Theorem \ref{Theorem7.9} (3).
	
	\item [(3)] Consider $X\subset \C^4$, defined by $\langle x,y\rangle\cap\langle u,v\rangle$ and let $f(x,y,u,v)=x+y+u+v$. $X$ is the $1$--point union of 2 planes with depth $(X,0)=1$. $X_0$ consists of 2 lines through $0$ with an embedded component of length $\varepsilon(X_0)=1$ and $X_t$ consists of 2 disjoint lines for $t \neq 0$. We get $\delta(X_0)=0, r'(X_0)=1, \mu(X_0)=-1$ and $\delta(X_t)=r'(X_t)=\mu(X_t)=0$. In this example $X$ is reduced and pure 2--dimensional with $\mu(X_t)$ not upper semicontinuous. 
	
	\item [(4)] Let $X\subset \C^6$ be the reduced surface with 3 components given by $\langle x^2y+y^2x+t(x^2+y^4)\rangle\cap\langle w-t, v, y,x\rangle\cap \langle w,u,y,x\rangle$. The first ideal of the intersection describes the deformation from a $D_4$--singularity to an $A_3$--singularity, the other 2 ideals describe 2 planes that meet in $0$. $X$ has depth $1$, $X_0$ consists of a $D_4$--singularity and 2 lines with an embedded point of length $\varepsilon(X_0)=1$. The general fibre $X_t$ is a reduced curve, consisting of an $A_3$--singularity at $0$, a line through $0$ and another disjoint line. We get $\delta(X_0)=r'(X_0)=\mu(X_0)=4$ and $\delta(X_t)=3, r'(X_t)=2$, $\mu(X_t)=4$. This family is $\mu$--constant along the trivial section (the $t$--axis) but not $\delta$--constant.
		
	The computation for these examples were done with {\sc Singular} \cite{DGPS15}, similar to the computations in Example \ref{Ex5.6}.
	\end{enumerate}
\end{Example}
\bigskip

\section{Comments, open problems and conjectures}

\subsection{Equisingularitiy for families of generically reduced curves}\label{subsec8.1}

The theory of equisingularity for familes of plane curves is well understood since many years. Part of it has been reviewed in Section 1 with emphasis on topological equisingularity and simultaneous resolution. 

In Section 2 we considered the case of families of reduced but not necessarily plane curve singularities, based on \cite{BuG80}. This had been continued in \cite{BGG80}, who investigated stronger conditions, e.g. Whitney regularity and equisaturation. Equisingular families of plane curve singularities with embedded points have been studied by A. Nobile \cite {No95}. Families of generically reduced curves have been considered in \cite{BrG90}, in \cite{BSS16} and in this paper but several questions remain still open, which we are going to discuss in this section.
\bigskip

{\bf General assumption:} Let $(X,x)$ be a surface singularity, $f:(X,x)\to (\C, 0)$ flat with $(X_0, x)=(f^{-1}(0), x)$ an isolated curve singularity and $f:X\to T$ a good representative with generically reduced curves as fibres. We keep also the other notations from Section 7.\\

If $(X_0, x)$ is reduced (equivalent to $(X,x)$ Cohen--Macaulay; then $(X,x)$ is pure dimensional and reduced) and if $\sigma:T\to X$ is a section of $f$, such that $X_t \smallsetminus \sigma(t)$ is smooth, we have the following implications, shown in the es-diagram below (diagram taken from \cite{BGG80}). 
\bigskip

\noindent
$
{\scriptsize
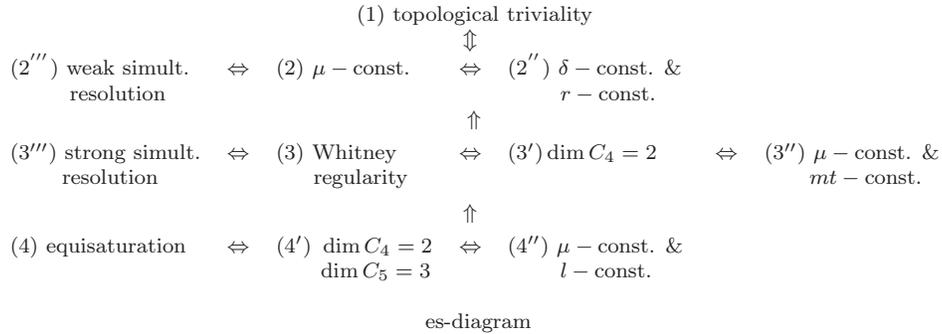
\begin{array}{lclclcl}
\multicolumn{7}{c}{\text{(1) topological triviality}} \\
& & & \Updownarrow & &&\\
(2^{'''})\   \text{weak simult.}  & \Leftrightarrow & (2)\ \mu-\text{const.} & \Leftrightarrow & (2^{''})\ \delta-\text{const. \& }&&\\
\ \ \ \ \ \ \ \ \text{resolution}&&&& \ \ \ \ \ \ \ r-\text{const.}\\
& & & \ \Uparrow & & \\[1.0ex]
(3''')\   \text{strong simult.} & \Leftrightarrow & (3)\ \text{Whitney} & \Leftrightarrow & (3') \dim C_4=2 & \Leftrightarrow & (3'')\ \mu-\text{const. \&}\\
\ \ \ \ \ \ \ \text{resolution} & & \ \ \ \ \  \text{regularity} & & & &\ \ \ \ \ \ mt-\text{const.}\\[1.0ex]
& & &  \Uparrow & & \\[1.0ex]
(4)\ \text{equisaturation}  & \Leftrightarrow & (4')\ \dim C_4=2  & \Leftrightarrow & (4'')\ \mu-\text{const. \&}\\
& & \ \ \ \ \ \ \dim C_5=3 & & \ \ \ \ \ \ \ l-\text{const.}\\ \\

\multicolumn{7}{c}{\text{ es-diagram}}
\end{array}
}$ \\
\bigskip

Moreover, we consider\\ 
{\scriptsize $(3^{iv})$ Zariski’s discrimininant criterion}.\\

We have $l(X_t, \sigma(t))=\delta(X'_t, 0)-\delta(X_t,\sigma(t))$ with $X'_t$ the plane curve obtained by a generic projection of $(X_t, \sigma (t))$ to $(\C^2, 0)$. $C_4(X_t, \sigma(t))$ denotes the cone of limits of tangent vectors at $X_t\smallsetminus \sigma (t)$ and $C_5(X_t, \sigma(t))$ the cone of limits of secant vectors.\\

The implications (1) $\Leftrightarrow (2) \Leftrightarrow (2'') \Leftrightarrow (2''')$ are due to \cite{BuG80}, $(3')  \Leftrightarrow (3''), (4')  \Leftrightarrow (4'')$ and $(4') \Rightarrow (3)$ to \cite{BGG80}, $(3) \Rightarrow (1)$ to Thom-Mather, $(3)  \Leftrightarrow (3''')$ to Teissier and $(3)  \Leftrightarrow (3')$ and $(4)  \Leftrightarrow (4')$ to Stutz (reference see \cite{BGG80}).\\

\noindent\textbf{Conjecture 1:} Topological triviality means the existence of a homeomorphism of pairs $h:(X, \sigma(T))\xrightarrow{\approx} (X_0\times T, \{x\}\times T)$ over $T$. In \cite{BuG80} this was shown to be equivalent to the existence of homeomorphisms $h_t:(B_\varepsilon, X_0, 0)\xrightarrow{\approx} (B_\varepsilon, X_t, \sigma(t))$ for each $t\in T$. We conjecture that the homeomorphisms $h_t, t\in T$, can be glued to one homeomorphism, providing embedded topological triviality. The question is whether there exists a homeomorphism of triples
\[
h:(B_\varepsilon\times T, X, \sigma(T))\xrightarrow{\approx}(B_\varepsilon\times T, C_0\times T, \{x\}\times T)
\]
over $T$.

This question was communicated to me by O. Nogueira and J. Snoussi. It is purely topological, independent from $X$ and $X_0$ being reduced or not.\\

The proof of non-embedded topological triviality in \cite{BuG80} used, for $X_0\subset\C^n, n\geq 3$, a theorem by Lickorish about the unknotting of cones in high codimension. A proof of the conjecture would require a generalization about the simultaneous unknotting of cones in a $1$--parametric family or a direct construction using vector fields as in \cite{Ti77}.\\

\noindent\textbf{Remark} (added in proof): Fern\'andez de Bobadilla communicated to me that Conjecture 1 has a positive answer. In section 2 of his paper \cite{Bo13} he introduces the notion of ''cuts'' and a technique for topological trivialization. Using the same arguments as in his generalization of Timourian's result Conjecture 1 follows, according to Bobadilla.\\

\noindent\textbf{Problem 2:} In \cite{BrG90} and in this paper we prove (1) $\Leftrightarrow (2'')$ $\Leftrightarrow (2''')$ and $(3'') \Leftrightarrow (3''')$ for families of generically reduced curves. We also show (1) $\Rightarrow  (2)$, while $(2) \Rightarrow$ (1) holds only if the general fibre $X_t$ is connected.

It is open, which of the remaining implications of the es-diagram above hold for families of generically reduced curves, resp. which additional assumptions have to be made.  \\

\noindent\textbf{Problem 3:} Zariski's discrimininant criterion for equisingularity of plane curves is also not fully understood in general for generically reduced curves.

Let $(X,x)$ be a pure dimensional surface singularity and $\pi: (X,x)\to(\C^2, 0)$ a finite morphism. The image $\pi(\Sing(\pi))$ of the singular locus of $\pi$ is called the discriminant of $\pi$. 

Zariski's discriminant criterion $(3^{iv})$ says that for a generic linear projection (the restriction to $(X,x)$ of a linear projection from an ambient space), the reduced discriminant is smooth.\\

The following was proved in \cite{BuG80}, Theorem 6.2.7. Let $f:(X,x)\to (\C,0)$ be flat with $(X_0, x)$ a reduced complete intersection curve, and let $\pi:(X,x)\to (\C, 0)$ be a generic projection. 
Then $\pi=(\pi_1, f):(X,x)\to (\C^2,0)$ satisfies Zariskis discriminant criterion iff a good representative $f:X\to T$ admits a section $\sigma$ such that $X_t\smallsetminus\sigma(t)$ is smooth and $\mu(X_t, \sigma(t))$ and $\mt(X_t, \sigma (t))$ are constant for $t\in T$. 

This result was complemented in \cite{BGG80}, Theoreme II.5, where the authors show that Zariski's discriminant criterion is equivalent to Whitney regularity of $X$ (i.e. $\Sing(X)$ is smooth and the pair ($X\smallsetminus\Sing (X),\Sing (X)$) satisfies the Whitney conditions (a) and (b)). 
\medskip

In \cite{BSS16}, Theorem 3.1, the authors generalize this to the following:

Let $(X,x)$ be reduced and irreducible with $1$--dimensional smooth singular locus. Then Zariski's discriminant criterion is equivalent to Whitney regularity.

\medskip

It remains still open whether the irreducibility and reducedness in this theorem is necessary in our general setting. More generally, which implications between $(3^{iv}), (3), (3')$ and ($\widetilde{3}''$) hold under the assumptions of Theorem 7.13 of this paper, where ($\widetilde{3}''$) means that $\delta(X_t, \sigma(t))\ ,\ r(X_t, \sigma(t))$ and $\mt(X, \sigma(t))$ are constant.\\

We like to add that for hypersurfaces of arbitrary dimension  \cite [Theorem 3.2] {Li00} shows that a Zariski-equisingular stratification (defined by constancy of ''dimensionality type'') is indeed a Whitney stratification.\\

\noindent\textbf{Problem 4:} Study properties of the Milnor number $\mu(f)$ (introduced by Goryunov, Mond and van Straten) and the invariant $m_1(X_0)$ (introduced by Nu\~no--Ballesteros and Tomazella)
for functions $f$ on a generically reduced curve or, more generally, on an INNS of arbitrary dimension (see Remark \ref{Rem9.9}). In particular their behaviour under deformations and their relation to other equisingularity conditions.\\

Some of the questions in Problem 1, 2 and 4 should have an easy answer for families of generically reduced curves, given the results of the present paper. The invariants from the es-diagram and from Problem 4 may be studied for families of higher-dimensional INNS and related to equisingularity questions. We can however not expect, that the invariants $\mu, \delta, r', mt$ are meaningful enough to study the topology (e.g. the higher homology groups) of the Milnor fiber in such families.\\

\subsection{Deformation of the normalization and of the parametrization}

Let $(C,0)\subset (\C^n,0)$ be a reduced curve singularity and consider the diagram
\[
\xymatrix{
(\overline{C}, \overline{0})\ar[d]_n\ar[dr]^\varphi &\\
(C,0)\ar@{^{(}->}[r]^j & (\C^n, 0)
}
\]
with $j$ the given embedding, $n$ the normalization and $\varphi$ the parametrization of $(C,0)$. We compare deformations of $n$ and of $\varphi$. By definition, a deformation of a morphism includes the deformation of the source and the target (cf. Appendix). Deformations of $n$ are hence at least as complicated as deformations of $(C,0)$. On the other hand, deformations of $\varphi$ are very simple, as $(\overline{C}, 0)$ and $(\C^n, 0)$ are both smooth. They can be explicitly described as in Section 4, just with $n$ components insteady of $2$.

For plane curves, however, deformations of $\varphi$ and of $n$ are up to isomorphism the same thing as we explained in Remark \ref{Remark4.4}, and the same holds for deformations with section (cf. \cite{GLS07}, Proposition 2.23). This identification allowed us in Theorem \ref{Theo4.1} to relate the base space of the semiuniversal deformation of $\varphi$ with the $\delta$--constant stratum in the seminuniversal deformation of $(C,0)$.\\

Consider now for $(C,0)\subset(\C^n,0)$ a $1$--parametric deformation of $\varphi$,
\[
(\overline{X}, \overline{x})\cong(\overline{C}\times \C,0)\xrightarrow{\Phi} (\C^n\times \C, 0)\xrightarrow{\pr}(\C,0),
\]
and let $(X,x) \subset (\C^n\times \C,0)$ be the image of $\Phi$, which is analytic since $\Phi$ is finite. If we endow $(X,x)$ with its reduced structure, then the restriction of $\pr$ to $(X,x)$,
\[
f: (X,x)\to (\C, 0),
\]
is flat. If $(X,x)$ is not Cohen Macaulay, then $(X_0, x)$ has an embedded component and $\varphi$ induces a deformation of $(X_0, x)$ but not of $(C,0)\cong(X_0^\red, x)$. Therefore we cannot expect a similar relation between the semiuniversal deformations of $\varphi$ and of $(C,0)$ as for plane curves.
\medskip

Nevertheless, $(C,0)$ has a semiuniversal deformation, which we denote by $\Psi:(\kd_C,0)\to (B_C,0)$. Consider the restriction $\Psi_\Delta$ of $\Psi$ to the $\delta$--constant stratum $\Delta\subset B_C$ of $\Psi$ and let $\Psi_{\overline{\Delta}}:\kd_{\overline{\Delta}}\to \overline{\Delta}$ be the pullback of $\Psi_\Delta$ via the normalization map $\overline{\Delta}\to \Delta$. Then $\Psi_{\overline{\Delta}}$ is flat and a $\delta$--constant deformation of the reduced curve singularity $(C,0)$ over a normal base space. 
By Theorem \ref{Theorem2.4} the normalization $\overline{\kd}_{\overline{\Delta}}\xrightarrow{m} \kd_{\overline{\Delta}}$ is a simultaneous normalization of $\Psi_{\overline{\Delta}}$. Since the deformation $\Psi_{\overline{\Delta}}$ can be embedded, we get a deformation of the normalization of $(C,0)$
\[
(\overline{\kd}_{\overline{\Delta}}, m^{-1}(0))\to (\kd_{\overline{\Delta}}, 0)\hookrightarrow(\C^n \times\overline{\Delta}, 0)\to(\overline{\Delta}, 0),
\]
and the composition $\Phi:(\overline{\kd}_{\overline{\Delta}}, m^{-1} (0) \to (\C^n\times \overline{\Delta}, 0)\to (\overline{\Delta}, 0)$ is a deformation of the parmaetrization $\varphi$. Hence $\Phi$ can be induced from the seminuniversal deformation $\kd_\varphi\to \kb_\varphi$ of $\varphi$ by some morphism $\beta:\overline{\Delta}\to \kb_\varphi$. $\beta$ is only unique up to first order, but nevertheless we may suggest (in view of Theorem \ref{Theo4.1}):\\

\noindent\textbf{Problem 5:} What can be said about $\beta$, in particular about its image in $\kb_\varphi$? In general probably not much since deformations of $\varphi$ are unobstructed, while deformations of $(C,0)$ can be arbitrary complicated (it is still an open problem whether there exist non--smooth rigid reduced curve singularities). It might nevertheless be interesting to study special classes of curves, such as complete intersections or reduced curves in $(\C^3, 0)$. For both classes the semiuniversal deformation can be explicitly described and has a smooth base space (for the latter case see \cite{BrG90}, Proposition 6.3).
\medskip

Note that in any case we have an explicit description of the semiuniversal deformation of $\varphi$.:

Any deformation of $\varphi :(\overline{C}, \overline{0})\to (\C^n,0)$ is given as in Section 4 (with $n$ components instead of 2). Also the semiuniversal deformation of $\varphi$ is analogous to the case of plane curves, which is explicitly described in \cite{GLS07}, Proposition II.2.27. We have to use a basis of $T^1_{(\overline{C}, \overline{0})\to (\C^n, 0)}$ which is the same as $M^0_\varphi$ in \cite{GLS07} but now with $n$ components (see \cite{GrL08}, Theorem 1.3).
\bigskip

For an isolated curve singularity $(C,0)\subset(\C^n, 0)$ we can define equisingularity by any of the es--conditions, i.e. conditions of the es--diagram, including Zariski's discriminant criterion ($3^{iv}$) above. So far, mainly $1$--parametric deformations have been considered, even for reduced curves. The following problem appears natural.\\

\noindent\textbf{Problem 6:} For a (generically) reduced curve singularity $(C,0)$ study deformations that satisfy any of the es--conditions over a reduced base space of arbitrary dimension. In particular does a versal es--deformation (to be understood in the category of reduced complex spaces) exist? What can be said about the base space if it exists (e.g. smoothness, irreducibility)? Note that for reduced plane curve singularities the restriction to the $\mu$--constant stratum of the semiuniversal deformation is a seminuniversal Zariski--equisingular deformation.
May be the case of complete intersections or of reduced curves in $(\C^3,0)$ is more feasible.

\subsection{Simultaneous normalization}

Our main new result about simultaneous normalization is Theorem \ref{theorem1} (resp. its global version Theorem \ref{theo6.16}). It says that a flat morphism $f:(X,x)\to (\C, 0)$ with  $(X_0,x)$ a positive dimensional INNS is equinormalizable iff $f$ is $\delta$--constant and without isolated points in the fibres. We conjecture that this result extends to familes over a normal base space.
\bigskip

\noindent\textbf{Conjecture 7:} Let $f:(X,x)\to (S,s)$ be flat with $(S,s)$ normal such that the non--normal locus $\NNor(f)$ finite over $S$ and $\dim(X_0, x)\geq 1$. Then $f$ is equinormalizable iff $f$ ist $\delta$--constant and without isolated points in the fibres.
\medskip

A proof of this conjecture would be an important step to answer some of the questions posed in Problem 6. The conjecture was proved by L\^{e} C\^ong Tr\`inh(\cite{LC15}, Theorem 3.6) for families of isolated curve singulrities under the assumption that $(S,s)$ is smooth and $(X,x)$ reduced and pure dimensional with normalization $(\overline{X}, \overline{x})$ being Cohen--Macaulay.
\bigskip

We may ask  whether in our local situation the functor of simultaneous normalizations is representable, as this was proved by Koll$\acute{\text{a}}$r \ref{Theo6.13} for projective morphisms. It would probably be of great help to prove Conjecture 7.\\

\noindent\textbf{Problem 8:} Let $f:(X,x)\to (S,s)$ be a morphism of complex germs and consider the functor $\SimNor(f)$ from complex germs to sets, associating to a germ $(T,t)$ the set
\[
\{\varphi:(T,t)\to (S,s)|\varphi^\ast f \text{ admits a simultaneous normalization over } (T,t)\}.
\]
Show that this functor is representable if the fibres of $f$ are generically reduced without isolated points and if $\NNor(f)$ is finite over $(S,s)$.
May be one has to assume that $(S,s)$ is normal or weakly normal.
\bigskip

Concerning projective morphisms Chiang--Hsieh and Lipman pose the question why the Hilbert polynomial, a global object, should be involved in the characterization of simultaneous normalization, which is a purely local phenomenon. Moreover, they ask if there are local invariants, somehow related to the Hilbert polynomial, that characterize equinormalizability at a point. This question makes sense also for local analytic morphisms when no Hilbert polynomial is available.
\bigskip

\noindent\textbf{Problem 9:} Let $f:(X,x) \to (S,s)$ by a morphism of complex germs. Does there exist numerical invariants of a complex space such that their constancy along the fibres characterizes 
the existence of a simultaneous normalization of $f$?
\medskip

One may think of certain multiplicities associated to the non--normal locus, to hyperplane sections of various dimensions or to strata of a stratification.
We have seen that our $\delta$--invariant gives a positive answer if the non-normal locus of $f$ is finite over $(S,s)$ and if $(S,s) = (\C,0)$ (and we conjecture it for $(S,s)$ normal).

\subsection{Equisingularity in positive characteristic}

Equisingular deformations of plane algebroid isolated curve singularities $R=K[[x,y]]/\langle  f\rangle$, $K$ an algebraically closed field of characteristic $p\geq 0$, have been studied in \cite{CGL07a} and \cite{CGL07b}, using Zariski's definition of equisingularity requiring equimultiplicity along sections through successive blowing up points. The aim was to prove smoothness of the equisingularity stratum in the semiuniversal deformation of the equation as well as of the parametrization. We do not go into details here but like to mention that the parametric approach works very much like in characteristic $0$, while deformations of the equation may depend on the characteristic for small $p$. In particular, the base space of the semiuniversal $\es$--deformation of the parametrization is always smooth, while this holds for deformations of the equation only in ''good'' characteristic. Furthermore, the classification of isolated hypersurface singularities has been started in \cite{GrK90} and in \cite{GrN16}. Otherwise not very much is known about equisingularity in positive characteristic. In particular, the conditions at the beginning of Section 8.1 have not been systematically studied for algebroid curve singularities. Therefore one might consider the following problem.
\bigskip

\noindent\textbf{Problem 10:} Let $K$ be an algebraically closed field of characteristic $\geq 0$ and $C$ an isolated algebroid curve singularity, i.e. $C=\Spec(R)$ with $R=K[[x_1, \ldots, x_n]]/I$ a $1$--dimensional analytic ring such $\Spec(R)\smallsetminus\frak{m}$ is regular ($\frak{m}$ the maximal ideal of $R$).  The problem is to study the es-conditions at the beginning of Section \ref{subsec8.1} and their relations for $1$--parametric deformations of $R$, whenever they make sense.
\medskip

The algebraic framework for this, in even greater generality, has been developed in \cite{CL06} for $R$ reduced. The authors prove there that for deformations over a normal base space the family is equinormalizable iff it is $\delta$--constant (even for perfect fields). \\

Finally, the following problem seems to be of general interest, e.g. in the context of generalized Jacobians.
\bigskip

\noindent\textbf{Problem 11:} Let $C$ be an isolated algebroid plane curve singularity, defined over an algebraically closed field of positive characteristic. Is the normalization of the $\delta$--constant stratum in the semiuniversal deformation of $C$ smooth? More generally, does Theorem \ref{Theo4.1} hold in this situation? \\
\smallskip

\section{Appendix: Deformation categories and functors}

For the convenience of the reader we collect notations and definitions concerning deformations of complex spaces, resp. of germs and good representatives, to be used throughout this article. As a general reference we refer to  \cite{GLS07}.
\medskip

\begin{enumerate}
\item [(1)] Let $\Phi:X\to S$ and $\varphi:T\to S$ be morphisms of complex spaces\footnote{We always assume that a complex space is separated and has a countable basis of the topology.}. We denote by $X_T$ the fibre product of $X$ with $T$ over $S$, giving rise to a Cartesian diagram
\[
X\times_ST =
\xymatrix{X_T\ar[r]\ar[d]_{\varphi^\ast\Phi} & X\ar[d]^\Phi\\
T\ar[r]_\varphi & S\ 
} \raisebox{-8ex}{,}
\]
with $\varphi^\ast\Phi$ called the \textbf{pull back} of $\Phi$. If $s\in S$, then $X_s:=\Phi^{-1}(s)$, denotes the \textbf{fibre} of $\Phi$ over $s$. $S$ is called the \textbf{base space} or the \textbf{parameter space} of $\Phi$. Analogous notations are used for morphisms of complex (space) germs.

\item [(2)] Let $(X,x)$ and $(S,s)$ be complex germs. A \textbf{deformation} of $(X,x)$ over $(S,s)$ consists of a flat morphism $\Phi:(\kx, x)\to (S,s)$ of complex germs together with a isomorphism $i:(X,x)\overset{\cong}{\rightarrow}(\kx_s, x)$. Hence it is given by a diagram as in (1) with $\Phi$ flat, but for germs.
We denote a deformation by $(i, \Phi)$ or just by $\Phi$.

\item [(3)] Let $f:(X,x)\to (T, t)$ be a morphism of complex germs. A \textbf{deformation of the morphism} $f$ over a complex germ $(S,s)$ is a Cartesian diagram
\[
\xymatrix{
(X,x)\ar[d]_f \ar@{^{(}->}[r]^i \ar@{}[dr]|{\Box}& (\kx, x)\ar[d]^F\ar@/^2pc/[dd]^\Phi\\
(T, t)\ar@{^{(}->}^j [r]^i\ar[d] \ar@{}[dr]|{\Box}& (\kt, t)\ar[d]^p\\
\{pt\}\ar@{^{(}->}[r] & (S,s)
}
\]

with $i$ and $j$ closed embeddings, $p$ and $\Phi$ flat (hence deformations of $(T,t)$ resp. of $(X,x)$ over $(S,s)$) and $\{pt\}$ the reduced point. We denote it by $(i,j,F,p)$ or just by $(F,p)$. A \textbf{deformation with (compatible) sections} of $f$ is a deformation $(F, p)$ together with sections $\sigma:(S,s)\to (\kt, t)$ of $p$ and $\overline{\sigma}:(S,s)\to (\kx, x)$ of $\Phi$ such that $\sigma=F\circ \overline{\sigma}$.
\item [(4)] A \textbf{morphism of two deformations} of $f$, $(i,j,F,p)$ over $(S,s)$ and $(i', j', F', p')$ over $(S', s')$ is given by morphism $\Psi_1:(\kx, x)\to (\kx', x')$, $\Psi_2:(\kt, t)\to (\kt', t')$ and $\varphi: (S,0)\to (S', s')$ such that the obvious diagram commutes. Forgetting the upper square of the diagram in (3) we get a morphism of deformations of $(T,t)$. If $(S,s)=(S', s')$ we always assume $\varphi=\id$. A morphism of deformations with sections is given by morphisms $\Psi_1$ and $\Psi_2$ that commute with the sections.
\item [(5)] Deformations of $f$ form a \textbf{category of deformations} of $f$, denoted by $\Def_f$. $\Def_f (S, s)$ denotes the subcategory of $\Def_f$ consisting of deformations of $f$ over $(S,s)$ with $\varphi$ being the identity on $(S,s)$. The category of deformations of $(X,x)$ are denoted by $\Def_{(X,x)}$. The diagram in (3) shows that we have forgetful functors $\Def_f\to \Def_{(X,x)}$ and $\Def_f\to \Def_{(T,t)}$. 

The categories of deformations with section are denoted by $\Def_f^{sec}$ and $\Def_f^{sec}(S,s)$ and similar for $(X,x)$ instead of $f$.

\item [(6)] Let $\varphi: (S',s')\to (S,s)$ be a morphism and $(F, p)$ a deformation of $f$ over $(S,s)$. The pull backs of $F$ and $p$ give rise to a deformation of $f$ over $(S',s')$, the \textbf{induced deformation} by $\mathbf{\varphi}$. If we denote by $\underline{\Def}_f (S,s)$ the set of \textbf{isomorphism classes of deformations} of $f$ over $(S,s)$, then we get a functor
\[
\begin{tabular}{p{0.8cm}p{5.76cm}p{0.1cm}p{4.5cm}}
$\underline{\Def}_f:$ &(category of complex space germs) & $\to$ & (category of sets),
\end{tabular}
\]
mapping $(S,s)$ to $\underline{\Def}_f(S,s)$ and a morphism of complex germs to the induced deformation. This is called the \textbf{deformation functor} of $f$. Similar notations are used for deformations of $(X,x)$ and for deformations with section.

\item [(7)]We denote by $T_\varepsilon$ the complex space consisting of one point with local ring $\C[\varepsilon]=\C\oplus\varepsilon\C, \varepsilon^2=0$. Deformations over $T_\varepsilon$ are called \textbf{infinitesimal deformations}. For isomorphism classes of deformations over $T_\varepsilon$ we use the notion
\[
T^1_{(X, x)}\vcentcolon= \underline{\Def}_{(X,x)} (T_\varepsilon)\ ,\ T^1_f:= \underline{\Def}_f(T_\varepsilon)\ ,
\]
and $T^{1,\sec}:=\underline{\Def}^{\sec}(T_\varepsilon)$ for deformations with section. Note that these sets are complex vector spaces (cf. \cite{GLS07}, section II.1.4. and Appendix C).

\item [(8)] Let $(X,x)$ be a complex subgerm of $(\C^n,0)$. A deformation $\Phi:(\kx,x)\to (S,s)$ of $(X,x)$ with section $\sigma$ is isomorphic to an \textbf{embedded deformation with trivial section}. That is, there exists a commutative diagram
\[
\xymatrix{(X,x)\ar@{^{(}->}[r]^i\ar[d] & (\kx, x)\ar@{^{(}->}[r]\  \ar[d]^\Phi & (\C^n\times S, (0,s))\ar[dl]^{pr}\\
\{pt\}\ar@{^{(}->}[r] & (S,s)\ar@/^1pc/^\sigma[u] &
}
\]
with $(\kx, x)$ embedded in $(\C^n\times S, (0,s))$, $pr$ the secod projection and $\sigma(S,s)=\{0\}\times (S,s)$. (cf. \cite{GLS07}, Prop. II.1.5., $\Phi$ need not be flat).

\item [(9)] Let $\Phi: (\kx, x)\to (S,s)$ be a morphism with fibre isomorphic to $(X,x)\subset(\C^n, 0)$ and with $(S,s)\subset (\C^k, 0)$. Choose an embedding of $\Phi$ as in (8) and let $B\subset \C^n$ resp. $D\subset \C^k$ be open balls around $0$ of radius $\varepsilon$ resp. $\eta$ and $X \subset B$ a closed subspace of $B$ representing $(X,x)$. \\
If $S\subset D$ resp. $\kx\subset B\times S$ are closed subspaces representing $(S,s)$ resp. $(\kx, x)$, and if $0<\eta (\varepsilon)<<\varepsilon$ are sufficiently small (depending on the required properties of $\Phi$), then $\Phi=pr|\kx:\kx\to S$ is called a \textbf{good representative} of $\Phi:(\kx,x)\to (S,s)$. By identifying $B \times \{t\}$ with $B$, we usually consider the fibres $\kx_t$,  $t \in S,$ of $\Phi$ as subsets of $B$ and identify the central fibre $\kx_s$ with $X$. Talking about good representatives, we allow $ \varepsilon$ and $\eta (\varepsilon)$ to shrink during an argument, if necessary.\\
$B$ is called a {\bf Milnor ball} if $\varepsilon$ is sufficiently small and the general fibre $\kx_t$ is called the \textbf{Milnor fibre} of $\Phi$.
\end{enumerate}

\end{document}